\documentclass[12pt]{article}
\usepackage{amsmath,amssymb}
\usepackage{graphicx}
\newcommand{\eqdef}{\stackrel{\text{def}}{=}}
\newcommand{\n}{\nonumber\\}
\newcommand{\bm}{\boldsymbol}
\newcommand{\ignore}[1]{}
\numberwithin{equation}{section}
\newcommand{\Romannumeral}[1]{\uppercase\expandafter{\romannumeral#1}}
\newcommand{\I}{\text{\Romannumeral{1}}}
\newcommand{\II}{\text{\Romannumeral{2}}}
\newcommand{\norm}[1]{|\!|#1|\!|}
\newcommand{\bin}[2]{(\!(#1,#2)\!)}
\newcommand{\bbin}[2]{\bigl(\!\bigl(#1,#2\bigr)\!\bigr)}

\allowdisplaybreaks[4]

\begin{document}

\baselineskip=20pt

\newcommand{\preprint}{
\vspace*{-20mm}\begin{flushleft}\end{flushleft}
   \begin{flushright}\normalsize \sf
    DPSU-16-1\\
  \end{flushright}}
\newcommand{\Title}[1]{{\baselineskip=26pt
  \begin{center} \Large \bf #1 \\ \ \\ \end{center}}}
\newcommand{\Author}{\begin{center}
  \large \bf Satoru Odake and Ryu Sasaki \end{center}}
\newcommand{\Address}{\begin{center}
     Faculty of Science, Shinshu University,\\
     Matsumoto 390-8621, Japan
   \end{center}}
\newcommand{\Accepted}[1]{\begin{center}
  {\large \sf #1}\\ \vspace{1mm}{\small \sf Accepted for Publication}
  \end{center}}

\preprint
\thispagestyle{empty}

\Title{Orthogonal Polynomials from Hermitian Matrices $\II$}

\Author

\Address
\vspace{1cm}

\begin{abstract}
This is the second part of the project `unified theory of classical orthogonal
polynomials of a discrete variable derived from the eigenvalue problems of
hermitian matrices.' In a previous paper, orthogonal polynomials having
Jackson integral measures were not included, since such measures cannot be
obtained from single infinite dimensional hermitian matrices.
Here we show that Jackson integral measures for the polynomials of the big
$q$-Jacobi family are the consequence of the recovery of self-adjointness of
the unbounded Jacobi matrices governing the difference equations of these
polynomials. The recovery of self-adjointness is achieved in an extended
$\ell^2$ Hilbert space on which a direct sum of two unbounded Jacobi matrices
acts as a Hamiltonian or a difference Schr\"odinger operator for an infinite
dimensional eigenvalue problem. The polynomial appearing in the upper/lower
end of Jackson integral constitutes the eigenvector of each of the two
unbounded Jacobi matrix of the direct sum.
We also point out that the orthogonal vectors involving the $q$-Meixner
($q$-Charlier) polynomials do not form a {\em complete basis\/} of the $\ell^2$
Hilbert space, based on the fact that the dual $q$-Meixner polynomials
introduced in a previous paper fail to satisfy the orthogonality relation.
The complete set of eigenvectors involving the $q$-Meixner polynomials is
obtained by constructing the duals of the dual $q$-Meixner polynomials
which require the two component Hamiltonian formulation.
An alternative solution method based on the closure relation, the Heisenberg
operator solution, is applied to the polynomials of the big $q$-Jacobi family
and their duals and $q$-Meixner ($q$-Charlier) polynomials.
\end{abstract}

\tableofcontents


\section{Introduction}
\label{sec:intro}

In a previous paper \cite{os12}, to be referred to I and its equation as
(I.2.3) etc. hereafter, a unified theory of classical orthogonal polynomials
of a discrete variable has been presented.
The classical orthogonal polynomials are polynomials satisfying the three
term recurrence relation and second order differential/difference equations.
On top of the well known Hermite, Laguerre, Jacobi and Bessel polynomials
\cite{szego} which satisfy second order differential equations,
the rest of about 40
classical orthogonal polynomials satisfying second order difference equations
are classified according to Askey scheme of hypergeometric orthogonal
polynomials \cite{askey,ismail,gasper,koeswart,koeswart2}.
In contrast to the ordinary orthogonal polynomials, which satisfy the three
term recurrence relation only, unified understanding of various properties
of the classical orthogonal polynomials is possible by considering them as
the main part of the eigenfunctions (eigenvectors) of self-adjoint second
order differential/difference operators which are exactly solvable.
The Schr\"odinger operators  in quantum mechanics are the typical examples
of such second order self-adjoint differential operators.
Thus we called the second order self-adjoint difference operators `discrete
Schr\"odinger operators' or the Hamiltonians in `discrete quantum mechanics'
\cite{os13,os14,os24}.
Various concepts and methods accumulated since the birth of quantum mechanics
are now available for unified understanding of classical orthogonal
polynomials \cite{crum,infhul,susyqm,os7}.
For example, the orthogonality is the consequence of the self-adjointness
and the orthogonality measures are provided by the square of the lowest
(the ground state) eigenfunctions (eigenvectors) which have no zeros due to
the oscillation theorem.

Orthogonal polynomials of a discrete variable \cite{nikiforov,koeswart}
have orthogonality measures concentrated on discrete points, either finite
or infinite in number.
The eigenvalue problems governing the classical orthogonal polynomials of
a discrete variable are based on a special class of hermitian matrices,
which are real symmetric tri-diagonal (Jacobi) matrices \eqref{Hdef} of
finite or infinite dimensions.
Since the spectrum of a Jacobi matrix is simple, the orthogonality of
eigenvectors (eigenpolynomials) is guaranteed except for unbounded infinite
dimensional ones, for which the self-adjointness could be broken.
The lowest eigenvectors (ground state vectors) of Jacobi matrices (to be
called the Hamiltonian hereafter) satisfy zero mode equations \eqref{phi0eq},
which can be solved easily.
By similarity transforming the Hamiltonians in terms of the ground state
vectors, one obtains difference operators which are {\em upper triangular\/}
in certain bases called the {\em sinusoidal coordinates\/} $\eta(x)$
\cite{os12,os14,os7}. The eigenvalues $\mathcal{E}(n)$ can be easily read off
as the coefficients of the highest degree term in the sinusoidal coordinates.
By solving the same equations in a different way, one obtains the
{\em dual polynomials\/} in the eigenvalue $\mathcal{E}(n)$
\cite{leonard}--\cite{terw3}.
All the Jacobi matrices corresponding to the classical orthogonal polynomials
are shown to have symmetries called {\em shape invariance\/} \cite{os12} and
the {\em closure relations\/} \cite{os7} which in turn provide universal
Rodrigues formula \eqref{univrod} and the coefficients of the three term
recurrence relations through the solutions of the Heisenberg equations
\cite{os7} for the sinusoidal coordinates (I.4.52)-(I.4.53).
These are the main results reported in I.

The purpose of the present paper is to rectify two shortcomings of I.
The first is that among many proposed dual orthogonal polynomials in I,
the dual $q$-Meixner and the dual $q$-Charlier polynomials do not satisfy
the orthogonality relations.
Secondly, those polynomials having the Jackson integral orthogonal measures
\cite{askey,gasper,koeswart}, {\em i.e.\/} the big $q$-Jacobi polynomial
(b$q$J) and related polynomials, are not included in I.
These two points have the same root.
The dual orthogonality relation is equivalent to the completeness relation
of the original polynomial.
That is, the $q$-Meixner ($q$-Charlier) polynomials  with the orthogonality
weight function \cite{koeswart} {\em do not form a complete set of basis of
the corresponding $\ell^2$ Hilbert space\/}\, \cite{atakishi1}.
There is no direct path to restore the breakdown of the completeness relation.
On the other hand, the failure of the orthogonality relation can be traced
to the breakdown of the self-adjointness.
Another Jacobi matrix (Hamiltonian) is introduced so that in the extended
Hilbert space the direct sum of the two Hamiltonians recover the
self-adjointness.
The resulting orthogonality relation in the {\em two component Hamiltonian
formalism takes the form of Jackson integral measure.\/}
By constructing the duals of the eigenpolynomials of the two component
Hamiltonian system, the remaining part of the complete set of eigenvectors
belonging to the $q$-Meixner ($q$-Charlier) polynomials is obtained
\cite{atakishi1}. The problem of identifying various dual polynomials,
{\em e.g.} big $q$-Jacobi and $q$-Meixner etc, has been tackled extensively
by Atakishiyev and Klimyk \cite{atakishi1}--\cite{atakishi4}
by methods different from ours.

This paper is organised as follows.
In section \ref{formulation} the essence of the program `orthogonal
polynomials from hermitian matrices' is reviewed with an eye to the possible
breakdown of the self-adjointness.
Several basic concepts, notions and notation are introduced and explained
in section \ref{polys}.
They are the sinusoidal coordinates, shape invariance, universal Rodrigues
formula, universal normalisation condition, dual polynomials, etc.
The lack of orthogonality of the naive dual $q$-Meixner polynomial is
examined in section \ref{sec:naive} and its recovery in the framework of
two component Hamiltonian formalism is presented in section \ref{sec:binham}.
The dual $q$-Charlier case is treated similarly in section \ref{sec:dqChar}.
The orthogonal polynomials having Jackson integral measures are discussed
in section \ref{sec:Jackson}.
Starting with the general structure of the two component Hamiltonian systems
in section \ref{sec:genstr}, the most generic case of the big $q$-Jacobi
polynomial is discussed in some detail in section \ref{sec:bqJ}.
The interesting features of the dual big $q$-Jacobi polynomial are explored
in section \ref{sec:dbqJ} \cite{atakishi1}--\cite{atakishi4}.
Restrictions of the parameters of the big $q$-Jacobi polynomial together with
certain limiting procedures produce other polynomials in the same family.
The big $q$-Laguerre and its dual, the Al-Salam-Carlitz $\I$ and the discrete
$q$-Hermite $\I$ are discussed in section \ref{sec:lim_bqJ}.
In \S\,\ref{sec:compqm} another set of infinitely many orthogonal polynomials
related with the $q$-Meixner ($q$-Charlier) polynomials is shown to
constitute the remaining part of the complete set of orthonormal bases
involving the $q$-Meixner ($q$-Charlier) polynomials.
This set is obtained by constructing the duals of the dual $q$-Meixner
($q$-Charlier) polynomials.
The discrete $q$-Hermite $\II$ is discussed in section \ref{sec:dqHII},
since it has also Jackson integral measure but it is defined on the full
integer lattice $x\in\mathbb{Z}$ in contrast to the half line integer lattice
$x\in\mathbb{Z}_{\ge0}$ in other examples.
Another polynomial defined on the full integer lattice, the $q$-Laguerre
polynomial is discussed in \S\ref{sec:qLag} by pursuing the well known
connection between the even degree Hermite and the Laguerre polynomials.
The universal Rodrigues formulas for those having Jackson integral measures are
presented in \S\,\ref{sec:uniRod}.
One interesting application of the orthogonal polynomials of a discrete
variable, the stationary stochastic processes called `Birth and Death (BD)
processes,' is discussed in section \ref{sec:bdp}.
The polynomials discussed in \S\,\ref{sec:dqMeix}--\S\,\ref{sec:dqHII} provide
new types of exactly solvable birth and death processes.
This is a supplementary sequence to \cite{bdproc}.
As explained in \S\,\ref{sec:bdp}, the exact solvability of BD processes hangs
on the completeness of the corresponding eigenpolynomials. The necessary
modifications to the solutions of the BD process corresponding to the
$q$-Meixner ($q$-Charlier) polynomials are mentioned at the end of
\S\,\ref{sec:bdp}.
The issue of the identification of the duals of the classical orthogonal
polynomials attracted many researchers \cite{atakishi1}--\cite{atakishi4}.
Most of them are, however, obtained {\em automatically by the interchange
$x\leftrightarrow n$} \eqref{duality1} from the original polynomials
{\em satisfying the universal normalisation condition} \eqref{Pzero}.
In section \ref{sec:nnorm} we propose new $q$-hypergeometric expressions
for four classical polynomials, the little $q$-Jacobi, little $q$-Laguerre,
Al-Salam-Carlitz $\II$ and the alternative $q$-Charlier ($q$-Bessel),
which satisfy the universal condition \eqref{Pzero}.
The final section is for a summary and comments.
In Appendix, an alternative solution method for the classical orthogonal
polynomials is applied to those polynomials discussed in the main text.
That is an algebraic method based on the symmetry properties of classical
orthogonal polynomials called closure relations \cite{os7,os12,os14}.
For the dual big $q$-Jacobi family and the $q$-Meixner ($q$-Charlier),
the same data govern the ordinary sector and the supplementary sector,
which is necessary for the completeness.

\section{Orthogonal Polynomials from Hermitian Matrices}
\label{sec:rdQM}

\subsection{Formulation}
\label{formulation}

Let us recapitulate the essence of the discrete quantum mechanics with real
shifts developed in I.
The Hamiltonian $\mathcal{H}=(\mathcal{H}_{x,y})$ is a {\em tri-diagonal real
symmetric\/} (Jacobi) matrix and its rows and columns are indexed by
integers $x$ and $y$, which take values in $\{0,1,\ldots,N\}$ (finite) or
$\mathbb{Z}_{\geq 0}$ (semi-infinite) or $\mathbb{Z}$ (full infinite).
Since the finite dimensional case is fully developed in I,
we concentrate here on the semi-infinite  ($\mathbb{Z}_{\geq 0}$) case.
The full infinite case will be discussed separately in \S\,\ref{sec:dqHII}.
Let us assume that the spectrum consists of discrete eigenvalues only 
and is bounded from below.
By adding a scalar matrix to the Hamiltonian, the lowest eigenvalue is
adjusted to be zero. This makes the Hamiltonian {\em positive semi-definite}.
Since the eigenvector corresponding to the zero eigenvalue has definite sign,
{\em i.e.\/} all the components are positive or negative, the Hamiltonian
$\mathcal{H}$ has the following tri-diagonal form ($x,y\in\mathbb{Z}_{\ge0}$)
\begin{equation}
  \mathcal{H}_{x,y}\eqdef
  -\sqrt{B(x)D(x+1)}\,\delta_{x+1,y}-\sqrt{B(x-1)D(x)}\,\delta_{x-1,y}
  +\bigl(B(x)+D(x)\bigr)\delta_{x,y},
  \label{Hdef}
\end{equation}
in which the potential functions $B(x)$ and $D(x)$ are real and positive
but vanish at the boundary
\begin{equation}
  B(x)>0\ \ (x\geq 0),\quad D(x)>0\ \ (x\geq 1),\ \ D(0)=0,
  \label{BDcond}
\end{equation}
and the Hamiltonian \eqref{Hdef} is real symmetric, $\mathcal{H}_{x,y}
=\mathcal{H}_{y,x}$.
Reflecting the positive semi-definiteness, the Hamiltonian \eqref{Hdef}
can be expressed in a factorised form:
\begin{align}
  &\mathcal{H}=\mathcal{A}^{\dagger}\mathcal{A},\qquad
  \mathcal{A}=(\mathcal{A}_{x,y}),
  \ \ \mathcal{A}^{\dagger}=((\mathcal{A}^{\dagger})_{x,y})
  =(\mathcal{A}_{y,x}),
  \label{factor}\\
  &\mathcal{A}_{x,y}\eqdef
  \sqrt{B(x)}\,\delta_{x,y}-\sqrt{D(x+1)}\,\delta_{x+1,y},\quad
  (\mathcal{A}^{\dagger})_{x,y}=
  \sqrt{B(x)}\,\delta_{x,y}-\sqrt{D(x)}\,\delta_{x-1,y}.
\end{align}
Here $\mathcal{A}$ ($\mathcal{A}^\dagger$) is an upper (lower) triangular
matrix with the diagonal and the super(sub)-diagonal entries only.
For simplicity in notation, we write $\mathcal{H}$, $\mathcal{A}$ and
$\mathcal{A}^{\dagger}$ as follows:
\begin{align}
  &e^{\pm\partial}=((e^{\pm\partial})_{x,y}),\quad
  (e^{\pm\partial})_{x,y}\eqdef\delta_{x\pm 1,y},\quad
  (e^{\partial})^{\dagger}=e^{-\partial},
  \label{partdef}\\
  &\mathcal{H}=-\sqrt{B(x)D(x+1)}\,e^{\partial}
  -\sqrt{B(x-1)D(x)}\,e^{-\partial}+B(x)+D(x)\n
  &\phantom{\mathcal{H}}=-\sqrt{B(x)}\,e^{\partial}\sqrt{D(x)}
  -\sqrt{D(x)}\,e^{-\partial}\sqrt{B(x)}+B(x)+D(x),
  \label{Hdef2}\\
  &\mathcal{A}=\sqrt{B(x)}-e^{\partial}\sqrt{D(x)},\quad
  \mathcal{A}^{\dagger}=\sqrt{B(x)}-\sqrt{D(x)}\,e^{-\partial}.
  \label{A,Ad}
\end{align}
(We suppress the unit matrix $\bm{1}=(\delta_{x,y})$:
$(B(x)+D(x))\bm{1}$ in \eqref{Hdef2}, $\sqrt{B(x)}\,\bm{1}$ in \eqref{A,Ad}.)
Note that the self-adjointness of the Hamiltonian \eqref{Hdef} is trivial
for finite systems, but non-trivial for infinite systems.

The Hamiltonian \eqref{Hdef} is a linear operator on the real $\ell^2$
Hilbert space with the inner product of two real vectors $f=(f(x))$ and
$g=(g(x))$ defined by
\begin{equation}
  (f,g)=\sum_{x=0}^{\infty}f(x)g(x),\quad|\!|f|\!|^2\eqdef(f,f)<\infty,
  \label{inpro}
\end{equation}
where the infinite sum is defined by the limit of $N$-truncated inner
product $(f,g)_N$:
\begin{equation}
  (f,g)\eqdef\lim_{N\to\infty}(f,g)_N,\quad
  (f,g)_N\eqdef\sum_{x=0}^Nf(x)g(x).
  \label{inproN}
\end{equation}

The Schr\"odinger equation is the eigenvalue problem for the hermitian
matrix $\mathcal{H}$,
\begin{equation}
  \mathcal{H}\phi_n(x)=\mathcal{E}(n)\phi_n(x)\quad
  (n=0,1,\ldots),\quad
  0=\mathcal{E}(0)<\mathcal{E}(1)<\cdots,
  \label{schreq}
\end{equation}
where the eigenvector $\phi_n=(\phi_n(x))$ is, by definition, of finite
norm, $|\!|\phi_n|\!|<\infty$.
Let us recall the fact that the spectrum of a Jacobi matrix is simple.
The ground state eigenvector, which is chosen positive $\phi_0(x)>0$
($x\in\mathbb{Z}_{\geq 0}$), satisfies the zero mode equation:
\begin{equation}
  \mathcal{A}\phi_0=0\ \Rightarrow\mathcal{H}\phi_0=0,\quad
  \sqrt{B(x)}\,\phi_0(x)=\sqrt{D(x+1)}\,\phi_0(x+1)
  \ \ (x\in\mathbb{Z}_{\geq 0}),
  \label{phi0eq}
\end{equation}
and it is easily obtained (convention: $\prod_{k=n}^{n-1}*=1$):
\begin{equation}
  \phi_0(x)=\phi_0(0)\prod_{y=0}^{x-1}\sqrt{\frac{B(y)}{D(y+1)}}
  \ \ (x\in\mathbb{Z}_{\geq 0}).
  \label{phi0}
\end{equation}
The self-adjointness of the Hamiltonian and the non-degeneracy of the
spectrum \eqref{schreq} imply that the eigenvectors are mutually orthogonal:
\begin{equation}
  (\phi_n,\phi_m)=\frac{\delta_{nm}}{d_n^2}\quad
  (n,m=0,1,\ldots),
  \label{ortho}
\end{equation}
where $d_n^2$ ($d_n>0$) is a normalisation constant.
It should be emphasised that $\phi_0(x)^2$ can be analytically continued to
the entire complex $x$-plane as a meromorphic function.
Likewise $d_n^2$ can also be analytically continued to the entire complex
$n$-plane as a meromorphic function.
They vanish on the negative integer lattice
\begin{equation}
  \phi_0(x)^2=0\ \ (x\in\mathbb{Z}_{<0}),\quad
  d_n^2=0\ \ (n\in\mathbb{Z}_{<0}).
  \label{mero}
\end{equation}

Now let us consider the self-adjointness of the above tri-diagonal infinite
dimensional Hamiltonian \eqref{Hdef}.
In general the self-adjointness of the Hamiltonian $\mathcal{H}$ means the
equality of the following two quantities,
\begin{equation}
  (f,\mathcal{H}g)=(\mathcal{H}f,g)\ \ \Bigl(\Leftrightarrow 
  \lim_{N\to\infty}\bigl((f,\mathcal{H}g)_N-(\mathcal{H}f,g)_N\bigr)=0\Bigr),
  \label{sadj}
\end{equation}
for an arbitrary choice of two vectors $f$ and $g$, $\norm{f},\norm{g}<\infty$.
Of course the summations of both sides should be absolutely convergent.

The tri-diagonality provides a simple criterion of the self-adjointness by
considering the action of the Hamiltonian $\mathcal{H}$ on the vector $f(x)$,
which has a factorised form $f(x)=\phi_0(x)\check{\mathcal{P}}(x)$:
\begin{align}
  &\quad\mathcal{H}f(x)=\sum_{y=0}^{\infty}\mathcal{H}_{x,y}f(y)\n
  &=\bigl(B(x)+D(x)\bigr)f(x)
  -\sqrt{B(x)D(x+1)}\,f(x+1)-\sqrt{B(x-1)D(x)}\,f(x-1)\n
  &=\phi_0(x)\Bigl(B(x)\bigl(\check{\mathcal{P}}(x)
  -\check{\mathcal{P}}(x+1)\bigr)
  +D(x)\bigl(\check{\mathcal{P}}(x)-\check{\mathcal{P}}(x-1)\bigr)\Bigr)
  =\phi_0(x)\widetilde{\mathcal{H}}\check{\mathcal{P}}(x).
  \label{Hf}
\end{align}
Here the similarity transformed Hamiltonian $\widetilde{\mathcal{H}}$ is
\begin{equation} 
  \widetilde{\mathcal{H}}
  \eqdef\phi_0(x)^{-1}\circ\mathcal{H}\circ\phi_0(x)
  =B(x)(1-e^{\partial})+D(x)(1-e^{-\partial}).
  \label{Ht}
\end{equation}
To verify the equality \eqref{sadj}, we use the $N$-truncated inner product
\eqref{inproN}.
For two vectors $f(x)=\phi_0(x)\check{\mathcal{P}}(x)$ and
$g(x)=\phi_0(x)\check{\mathcal{Q}}(x)$, we have
\begin{align}
  (f,\mathcal{H}g)_N&=\sum_{x=0}^N\phi_0(x)\check{\mathcal{P}}(x)\phi_0(x)
  \Bigl(B(x)\bigl(\check{\mathcal{Q}}(x)-\check{\mathcal{Q}}(x+1)\bigr)
  +D(x)\bigl(\check{\mathcal{Q}}(x)-\check{\mathcal{Q}}(x-1)\bigr)\Bigr)\n
  &=\sum_{x=0}^N\phi_0(x)^2\bigl(B(x)+D(x)\bigr)\check{\mathcal{P}}(x)
  \check{\mathcal{Q}}(x)
  -\sum_{x=0}^N\phi_0(x)^2B(x)\check{\mathcal{P}}(x)\check{\mathcal{Q}}(x+1)\n
  &\qquad
  -\sum_{x=0}^N\phi_0(x)^2D(x)\check{\mathcal{P}}(x)\check{\mathcal{Q}}(x-1),
  \label{fHgN}
\end{align}
and the last term is rewritten with the help of the zero mode equation
\eqref{phi0eq}
\begin{align}
  &\quad\sum_{x=0}^N\phi_0(x)^2D(x)\check{\mathcal{P}}(x)
  \check{\mathcal{Q}}(x-1)
  =\sum_{x=1}^N\phi_0(x)^2D(x)\check{\mathcal{P}}(x)\check{\mathcal{Q}}(x-1)\n
  &=\sum_{x=0}^{N-1}\phi_0(x+1)^2D(x+1)\check{\mathcal{P}}(x+1)
  \check{\mathcal{Q}}(x)
  =\sum_{x=0}^{N-1}\phi_0(x)^2B(x)\check{\mathcal{P}}(x+1)
  \check{\mathcal{Q}}(x).
  \label{fHgNlast}
\end{align}
The expression of $(\mathcal{H}f,g)_N$ can be obtained by exchanging
$\check{\mathcal{P}}$ and $\check{\mathcal{Q}}$.
The criterion of the self-adjointness of the Hamiltonian with respect to
the inner product $(\cdot,\cdot)$ \eqref{inpro}--\eqref{inproN} is the
vanishing of the difference in $N\to\infty$ limit,
\begin{align}
  0&=\lim_{N\to\infty}\bigl((f,\mathcal{H}g)_N-(\mathcal{H}f,g)_N\bigr)\n
  &=\lim_{N\to\infty}\phi_0(N)^2B(N)\bigl(\check{\mathcal{P}}(N+1)
  \check{\mathcal{Q}}(N)
  -\check{\mathcal{P}}(N)\check{\mathcal{Q}}(N+1)\bigr).
  \label{criter1}
\end{align}

\subsection{Hypergeometric orthogonal polynomials and their duals}
\label{polys}

Each specific theory to be discussed hereafter depends on a certain set of
parameters, to be denoted symbolically by $\bm{\lambda}$.
Various quantities are expressed like, 
$\mathcal{H}(\bm{\lambda})=\mathcal{A}(\bm{\lambda})^\dagger
\mathcal{A}(\bm{\lambda})$, $\mathcal{E}(n;\bm{\lambda})$,
$\phi_n(x;\bm{\lambda})$, and $P_n(\eta(x;\bm{\lambda});\bm{\lambda})$, etc.
The parameter $q$ is $0<q<1$ and $q^{\bm{\lambda}}$ stands for
$q^{(\lambda_1,\lambda_2,\dots)}=(q^{\lambda_1},q^{\lambda_2},\ldots)$.
As shown in I, for special choices of the potential functions
$B(x;\bm{\lambda})$ and $D(x;\bm{\lambda})$, the eigenvalue problem
\eqref{schreq} is exactly solvable and various hypergeometric orthogonal
polynomials $\{P_n\}$, the ($q$-)Meixner, (dual) ($q$-)Hahn, ($q$-)Racah,
etc. are obtained as the main part of the eigenvector
\begin{align}
  &\phi_n(x;\bm{\lambda})
  =\phi_0(x;\bm{\lambda})\check{P}_n(x;\bm{\lambda}),\quad
  \check{P}_n(x;\bm{\lambda})\eqdef
  P_n\bigl(\eta(x;\bm{\lambda});\bm{\lambda}\bigr),
  \label{factsol}\\
  &\widetilde{\mathcal{H}}(\bm{\lambda})\check{P}_n(x;\bm{\lambda})
  =\mathcal{E}(n;\bm{\lambda})\check{P}_n(x;\bm{\lambda})
  \ \ (n=0,1,\ldots),
  \label{polysol}
\end{align}
in which $\phi_0(x;\bm{\lambda})$ is the ground state eigenvector
\eqref{phi0} and $\eta(x;\bm{\lambda})$ is a certain function of $x$,
called the {\em sinusoidal coordinate}.
In other words, the similarity transformed Hamiltonian
$\widetilde{\mathcal{H}}$ is {\em triangular} in the basis spanned by
$1,\eta(x;\bm{\lambda}),\ldots,\eta(x;\bm{\lambda})^n,\ldots$,
\begin{equation}
  \widetilde{\mathcal{H}}(\bm{\lambda})\eta(x;\bm{\lambda})^n
  =\mathcal{E}(n;\bm{\lambda})\eta(x;\bm{\lambda})^n+
  \text{lower degrees in}\ \eta(x;\bm{\lambda}),
  \label{lowtri}
\end{equation}
which leads to the {\em solvability}.
The eigenvalue $\mathcal{E}(n;\bm{\lambda})$ can be easily obtained as the
coefficient of the highest degree monomial in $\eta(x;\bm{\lambda})$.
Five different types of the sinusoidal coordinates are known \cite{os12,os14}:
\begin{equation}
  \eta(x;\bm{\lambda})=x,\ \ x(x+d),\ \ 1-q^x,\ \ q^{-x}-1,
  \ \ (q^{-x}-1)(1-dq^x),
  \label{listsinu}
\end{equation}
in which $d$ is a real parameter.
They all satisfy the universal boundary condition
\begin{equation}
  \eta(0;\bm{\lambda})=0.
  \label{etazerocond}
\end{equation}
The eigenpolynomials are chosen to satisfy the {\em universal boundary
condition}
\begin{equation}
  \check{P}_n(0;\bm{\lambda})
  =P_n\bigl(\eta(0;\bm{\lambda});\bm{\lambda}\bigr)=P_n(0;\bm{\lambda})=1
  \ \ (n=0,1,2,\ldots).
  \label{Pzero}
\end{equation}
They satisfy the orthogonality condition $(n,m=0,1,\ldots)$
\begin{equation}
  (\phi_n,\phi_m)
  =\sum_{x=0}^{\infty}\phi_0(x;\bm{\lambda})^2
  P_n\bigl(\eta(x;\bm{\lambda});\bm{\lambda}\bigr)
  P_m\bigl(\eta(x;\bm{\lambda});\bm{\lambda}\bigr)
  =\frac{\delta_{nm}}{d_n(\bm{\lambda})^2},
  \label{ortho2}
\end{equation}
in which the orthogonality weight $\phi_0(x;\bm{\lambda})^2$ is expressed
in terms of ($q$-)shifted factorials \cite{os12}.

The above triangularity \eqref{lowtri} does not necessarily imply the
existence of explicit expressions of the corresponding eigenpolynomials,
like the hypergeometric form.
A stronger condition governing the parameter dependence of the Hamiltonian
or the potential functions, called {\em shape invariance\/}
(see section \Romannumeral{4} of I for more details) plays that role.
The shape invariance provides the explicit forms of all eigenvalues
$\{\mathcal{E}(n;\bm{\lambda})\}$ \eqref{spectrumform} and the universal
Rodrigues type formula of the eigenpolynomials
$\{P_n(\eta(x;\bm{\lambda});\bm{\lambda})\}$ \eqref{univrod}.

The shape invariance condition is
\begin{equation}
  \mathcal{A}(\bm{\lambda})\mathcal{A}(\bm{\lambda})^{\dagger}
  =\kappa\mathcal{A}(\bm{\lambda}+\bm{\delta})^{\dagger}
  \mathcal{A}(\bm{\lambda}+\bm{\delta})
  +\mathcal{E}(1;\bm{\lambda}),
  \label{shapeinv1}
\end{equation}
in which $\bm{\delta}$ denotes the shift of the parameters, $\kappa$ is
a positive constant and $\mathcal{E}(1;\bm{\lambda})$ is the eigenvalue
of the first excited state $\mathcal{E}(1)>0$ with the explicit parameter
dependence. In terms of the potential functions the above condition means 
the following two relations:
\begin{align}
  B(x+1;\bm{\lambda})D(x+1;\bm{\lambda})
  &=\kappa^2
  B(x;\bm{\lambda}+\bm{\delta})D(x+1;\bm{\lambda}+\bm{\delta}),
  \label{shapeinv1cond1}\\
  B(x;\bm{\lambda})+D(x+1;\bm{\lambda})
  &=\kappa\bigl(
  B(x;\bm{\lambda}+\bm{\delta})+D(x;\bm{\lambda}+\bm{\delta})\bigr)
  +\mathcal{E}(1;\bm{\lambda}).
  \label{shapeinv1cond2}
\end{align}
By using the discrete counterpart of Crum's theorem \cite{crum,os15,os22},
the spectrum is simply generated by $\mathcal{E}(1,\bm{\lambda})$:
\begin{equation}
  \mathcal{E}(n;\bm{\lambda})
  =\sum_{s=0}^{n-1}\kappa^s\mathcal{E}(1;\bm{\lambda}+s\bm{\delta}),
  \label{spectrumform}
\end{equation}
and the corresponding eigenvectors are generated from the known form
of the ground state eigenvector $\phi_0(x;\bm{\lambda})$ \eqref{phi0}
together with the multiple action of the successive
$\mathcal{A}(\bm{\lambda})^\dagger$ operator:
\begin{equation}
  \phi_n(x;\bm{\lambda})\propto
  \mathcal{A}(\bm{\lambda})^{\dagger}
  \mathcal{A}(\bm{\lambda}+\bm{\delta})^{\dagger}\cdots
  \mathcal{A}\bigl(\bm{\lambda}+(n-1)\bm{\delta}\bigr)^{\dagger}
  \phi_0(x;\bm{\lambda}+n\bm{\delta}).
  \label{eigvecform}
\end{equation}
Let us introduce an auxiliary function $\varphi(x)$
defined on the entire complex $x$-plane (I.4.12), (I.4.23),
\begin{equation}
  \varphi(x;\bm{\lambda})\eqdef
  \frac{\eta(x+1;\bm{\lambda})-\eta(x;\bm{\lambda})}{\eta(1;\bm{\lambda})}.
  \label{varphidef}
\end{equation}
On the non-negative integer lattice it satisfies
\begin{equation}
 \varphi(x;\bm{\lambda})= \sqrt{\frac{B(0;\bm{\lambda})}{B(x;\bm{\lambda})}}
  \frac{\phi_0(x;\bm{\lambda}+\bm{\delta})}{\phi_0(x;\bm{\lambda})}.
\end{equation}
Then we have
\begin{align}
  \phi_0(x;\bm{\lambda}+\bm{\delta})^{-1}\circ
  \mathcal{A}(\bm{\lambda})\circ\phi_0(x;\bm{\lambda})
  &=\sqrt{B(0;\bm{\lambda})}\,\varphi(x;\bm{\lambda})^{-1}(1-e^{\partial}),
  \label{preF}\\
  \phi_0(x;\bm{\lambda})^{-1}\circ
  \mathcal{A}(\bm{\lambda})^{\dagger}\circ\phi_0(x;\bm{\lambda}+\bm{\delta})
  &=\frac1{\sqrt{B(0;\bm{\lambda})}}
  \bigl(B(x;\bm{\lambda})-D(x;\bm{\lambda})e^{-\partial}\bigr)
  \varphi(x;\bm{\lambda})
  \label{preB}\\
  &=\sqrt{B(0;\bm{\lambda})}\,\phi_0(x;\bm{\lambda})^{-2}\circ(1-e^{-\partial})
  \varphi(x;\bm{\lambda})^{-1}\circ\phi_0(x;\bm{\lambda}+\bm{\delta})^2.
  \nonumber
\end{align}

In terms of the auxiliary function $\varphi(x;\bm{\lambda})$, 
let us introduce the forward and backward shift operators
$\mathcal{F}(\bm{\lambda})$ (I.4.18) and $\mathcal{B}(\bm{\lambda})$ (I.4.19),
\begin{equation}
  \mathcal{F}(\bm{\lambda})
  =B(0;\bm{\lambda})\varphi(x;\bm{\lambda})^{-1}(1-e^{\partial}),\quad
  \mathcal{B}(\bm{\lambda})
  =\frac1{B(0;\bm{\lambda})}
  \bigl(B(x;\bm{\lambda})-D(x;\bm{\lambda})e^{-\partial}\bigr)
  \varphi(x;\bm{\lambda}).
  \label{FB}
\end{equation}
The following forward and backward shift 
relations
\begin{equation}
  \mathcal{F}(\bm{\lambda})\check{P}_n(x;\bm{\lambda})
  =\mathcal{E}(n;\bm{\lambda})\check{P}_{n-1}(x;\bm{\lambda}+\bm{\delta}),
  \quad
  \mathcal{B}(\bm{\lambda})\check{P}_{n-1}(x;\bm{\lambda}+\bm{\delta})
  =\check{P}_n(x;\bm{\lambda}),
  \label{FP=}
\end{equation}
are the consequences of the shape invariance for the polynomials satisfying the
universal boundary condition \eqref{Pzero}.
Starting from $P_0(\eta;\bm{\lambda})\eqdef 1$, $\check{P}_n(x;\bm{\lambda})$
can be written as
\begin{equation}
  \check{P}_n(x;\bm{\lambda})=\mathcal{B}(\bm{\lambda})
  \mathcal{B}(\bm{\lambda}+\bm{\delta})\cdots
  \mathcal{B}\bigl(\bm{\lambda}+(n-1)\bm{\delta}\bigr)\cdot 1.
  \label{RodrigB}
\end{equation}
{\em The universal Rodrigues formula\/} for all the hypergeometric orthogonal
polynomials (of a discrete variable) satisfying the universal boundary
conditions $P_n(0;\bm{\lambda})=1$ \eqref{Pzero},
$\phi_0(0;\bm{\lambda})=1$ (I.2.19) and
$\phi_0(x)^2=0$ ($x\in\mathbb{Z}_{<0}$) \eqref{mero} reads:
\begin{align}
  P_n\bigl(\eta(x;\bm{\lambda});\bm{\lambda}\bigr)
  &=\phi_0(x;\bm{\lambda})^{-2}\,
  \mathcal{D}(\bm{\lambda})
  \mathcal{D}(\bm{\lambda}+\bm{\delta})\cdots
  \mathcal{D}\bigl(\bm{\lambda}+(n-1)\bm{\delta}\bigr)\cdot
  \phi_0(x;\bm{\lambda}+n\bm{\delta})^2,
  \label{univrod}\\
  \mathcal{D}(\bm{\lambda})&\eqdef
  (1-e^{-\partial})\varphi(x;\bm{\lambda})^{-1}.
  \label{Dlamdef}
\end{align}
For the orthogonal polynomials having Jackson integral measures, both of
the conditions $P_n(0;\bm{\lambda})=1$ \eqref{Pzero} and
$\phi_0(0;\bm{\lambda})=1$ are not satisfied.
Thus the universal Rodrigues formula gets certain modification as shown
in section \ref{sec:uniRod}.
It should be emphasised that the conventional Rodrigues formulas
\cite{koeswart,koeswart2} were derived one by one for each polynomial.
 
The dual polynomial arises naturally as the solution of the original
eigenvalue problem \eqref{schreq} or \eqref{polysol} obtained in a
different way.
The similarity transformed eigenvalue problem
$\widetilde{\mathcal{H}}\,v(x)=\mathcal{E}\,v(x)$ \eqref{polysol} can be
rewritten into an explicit matrix form with the change of the notation
$v(x)\to{}^t(Q_0,Q_1,\ldots,Q_x,\ldots)$
\begin{equation}
  \sum_{y}\widetilde{\mathcal{H}}_{x,y}Q_y=\mathcal{E}Q_x
  \ \ (x=0,1,\ldots).
\end{equation}
Because of the tri-diagonality of $\widetilde{\mathcal{H}}$, it is
in fact a three term recurrence relation for $\{Q_x\}$ as polynomials
in $\mathcal{E}$:
\begin{equation}
  \mathcal{E}Q_x(\mathcal{E})
  =B(x)\bigl(Q_x(\mathcal{E})-Q_{x+1}(\mathcal{E})\bigr)
  +D(x)\bigl(Q_x(\mathcal{E})-Q_{x-1}(\mathcal{E})\bigr)
  \ \ (x=0,1,\ldots).
  \label{dual3term}
\end{equation}
Starting with the boundary (initial) condition
\begin{equation}
  Q_0=1,
  \label{Qzero}
\end{equation}
we determine $Q_x(\mathcal{E})$ as a degree $x$ polynomial in $\mathcal{E}$.
It is easy to see
\begin{equation}
  Q_x(0)=1\ \ (x=0,1,\ldots).
  \label{Qxzero}
\end{equation}
Note that the boundary condition $D(0)=0$ guarantees that $Q_{-1}$ decouples
from the above three term recurrence relation \eqref{dual3term}.
When $\mathcal{E}$ is replaced by the actual value of the $n$-th eigenvalue
$\mathcal{E}(n)$ \eqref{lowtri} in $Q_x(\mathcal{E})$, we obtain the
explicit form of the eigenvector
\begin{equation}
  \sum_{y}\widetilde{\mathcal{H}}_{x,y}Q_y\bigl(\mathcal{E}(n)\bigr)
  =\mathcal{E}(n)Q_x\bigl(\mathcal{E}(n)\bigr)
  \ \ (x=0,1,\ldots).
\end{equation}

The two expressions (polynomials) for the eigenvectors of the problem
\eqref{schreq} belonging to the eigenvalue $\mathcal{E}(n)$, $P_n(\eta(x))$
and $Q_x(\mathcal{E}(n))$ are in fact equal on the integer lattice points
\begin{equation}
  P_n\bigl(\eta(x;\bm{\lambda});\bm{\lambda}\bigr)
  =Q_x\bigl(\mathcal{E}(n;\bm{\lambda});\bm{\lambda}\bigr)
  \ \ (n=0,1,\ldots\,;\ x=0,1,\ldots),
  \label{Duality}
\end{equation}
due to the simplicity of the spectrum of Jacobi matrices and the boundary
conditions \eqref{Pzero} and \eqref{Qxzero}.
When $n\in\mathbb{Z}_{\ge0}$ is fixed, the hypergeometric expression of
the polynomial $P_n(\eta(x))$, say (I.5.28) for the Hahn polynomial,
or (I.5.73) for the $q$-Racah polynomial, etc.,
\begin{alignat*}{3}
  \text{Hahn}&:
  \ P_n\bigl(\eta(x);\bm{\lambda}\bigr)={}_3F_2\Bigl(
  \genfrac{}{}{0pt}{}{-n,\,n+a+b-1,\,-x}{a,\,-N}\!\Bigm|\!1\Bigr),
  &\bm{\lambda}&=(a,b,N),
  &\qquad&(\text{I}.5.28)\\
  \quad\ \text{$q$-Racah}&:
  \ P_n\bigl(\eta(x;\bm{\lambda});\bm{\lambda}\bigr)={}_4\phi_3\Bigl(
  \genfrac{}{}{0pt}{}{q^{-n},\,\tilde{d}q^n,\,q^{-x},\,dq^x}
  {a,\,b,\,c}\!\Bigm|\!q\,;q\Bigr),
  &\quad q^{\bm{\lambda}}&=(a,b,c,d),
  &\qquad&(\text{I}.5.73)
\end{alignat*}
is a degree $n$ polynomial in $\eta(x;\bm{\lambda})$ ($x$ for the Hahn and
$(q^{-x}-1)(1-{d}q^x)$ for the $q$-Racah), whereas when
$x\in\mathbb{Z}_{\ge0}$ is fixed, it is a degree $x$ polynomial in
$\mathcal{E}(n;\bm{\lambda})$ ($n(n+a+b-1)$ for the Hahn and
$(q^{-n}-1)(1-\tilde{d}q^n)$, $\tilde{d}\eqdef abcd^{-1}q^{-1}$ for the
$q$-Racah).
When one wants to express the dual polynomial of $P_n(\eta(x))$,
{\em i.e.\/} $Q_x(\mathcal{E}(n))$, as a degree $n$ polynomial on the
$x$ lattice, it is
\begin{equation}
  P^{\text{d}}_n\bigl(\eta^{\text{d}}(x;\bm{\lambda});\bm{\lambda}\bigr)
  \eqdef Q_n\bigl(\mathcal{E}(x;\bm{\lambda});\bm{\lambda}\bigr).
\end{equation}
Its explicit form is obtained by the interchange $x\leftrightarrow n$ in
the hypergeometric expression, for example:
\begin{alignat*}{2}
  \text{dual Hahn}&:
  \ P^{\text{d}}_n\bigl(\eta^{\text{d}}(x;\bm{\lambda});\bm{\lambda}\bigr)
  ={}_3F_2\Bigl(
  \genfrac{}{}{0pt}{}{-x,\,x+a+b-1,\,-n}{a,\,-N}\!\Bigm|\!1\Bigr),
  &\quad\eta^{\text{d}}(x;\bm{\lambda})&\eqdef\mathcal{E}(x;\bm{\lambda}),\\
  \text{dual $q$-Racah}&:
  \ P^{\text{d}}_n\bigl(\eta^{\text{d}}(x;\bm{\lambda});\bm{\lambda}\bigr)
  ={}_4\phi_3\Bigl(
  \genfrac{}{}{0pt}{}{q^{-x},\,\tilde{d}q^x,\,q^{-n},\,dq^n}
  {a,\,b,\,c}\!\Bigm|\!q\,;q\Bigr),
  &\quad\mathcal{E}^{\text{d}}(n;\bm{\lambda})&\eqdef\eta(n;\bm{\lambda}).
\end{alignat*}
Here is the summary of the dual correspondence:
  \begin{align}
  &x\leftrightarrow n,\quad
  \eta(x)\leftrightarrow\mathcal{E}(n),\quad
  \eta(0)=0\leftrightarrow\mathcal{E}(0)=0,
  \label{duality1}\\
  &B(x)\leftrightarrow -A_n,\quad
  D(x)\leftrightarrow -C_n,\quad D(0)=0=C_0,\quad
  \frac{\phi_0(x)}{\phi_0(0)}\leftrightarrow\frac{d_n}{d_0},
  \label{duality2}
\end{align}
in which $A_n$ and $C_n$ are the coefficients of the three term recurrence
relation of the orthogonal polynomials $P_n(\eta(x))$, \eqref{3rel}.

The orthogonality relation of the dual polynomials is related to the
{\em normalised} orthogonality relation of the original polynomials:
\begin{align}
  (\hat{\phi}_n,\hat{\phi}_m)
  &=\sum_{x=0}^{\infty}\phi_0(x;\bm{\lambda})^2d_nd_m
  P_n\bigl(\eta(x;\bm{\lambda});\bm{\lambda}\bigr)
  P_m\bigl(\eta(x;\bm{\lambda});\bm{\lambda}\bigr)=\delta_{nm}
  \ \ (n,m=0,1,\ldots),
  \label{normortho}\\
  &\quad\hat{\phi}_n(x;\bm{\lambda})\eqdef
  \phi_0(x;\bm{\lambda})d_n(\bm{\lambda})
  P_n\bigl(\eta(x;\bm{\lambda});\bm{\lambda}\bigr).
\end{align}
By multiplying $\hat{\phi}_n(y)$ to \eqref{normortho} and summing over $n$ 
produces, under certain conditions, the completeness relation
\begin{align}
  \sum_{n=0}^{\infty}\hat{\phi}_n(x)\hat{\phi}_n(y)
  &=\sum_{n=0}^{\infty}\phi_0(x)\phi_0(y)d_n^2
  P_n\bigl(\eta(x)\bigr)P_n\bigl(\eta(y)\bigr)=\delta_{xy}
  \ \ (x,y=0,1,\ldots)\n
  &=\sum_{n=0}^{\infty}\phi_0(x)\phi_0(y)d_n^2
  Q_x\bigl(\mathcal{E}(n)\bigr)Q_y\bigl(\mathcal{E}(n)\bigr),
  \label{dualortho}
\end{align}
or the dual orthogonality relation
\begin{align}
  &\quad\sum_{n=0}^{\infty}d_n(\bm{\lambda})^2
  P_n\bigl(\eta(x;\bm{\lambda});\bm{\lambda}\bigr)
  P_n\bigl(\eta(y;\bm{\lambda});\bm{\lambda}\bigr)\n
  &=\sum_{n=0}^{\infty}d_n(\bm{\lambda})^2
  Q_x\bigl(\mathcal{E}(n;\bm{\lambda});\bm{\lambda}\bigr)
  Q_y\bigl(\mathcal{E}(n;\bm{\lambda});\bm{\lambda}\bigr)
  =\frac{\delta_{xy}}{\phi_0(x;\bm{\lambda})^2}.
  \label{Qnortho}
\end{align}
As is well known, the above completeness relation \eqref{dualortho} simply
means that any square summable vector $f(x)$, $\norm{f}<\infty$ in the
$\ell^2$ Hilbert space can be expanded by $\{\hat{\phi}_n\}$:
\begin{equation*}
  f(x)=\sum_{n=0}^{\infty}(\hat{\phi}_n,f)\hat{\phi}_n(x),
\end{equation*}
and in particular $f(x)=\delta_{xy}$ leads to \eqref{dualortho}.
For the orthogonal polynomials defined on a finite lattice, for example,
the ($q$-)Racah polynomials, the completeness relation is an automatic
consequence of the original orthogonality relation \eqref{normortho}.
However, for some orthogonal polynomials defined over non-negative integers,
{\em e.g.\/} the $q$-Meixner and the $q$-Charlier polynomials, the
completeness relation \eqref{dualortho} or the dual orthogonality relation
\eqref{Qnortho} fails \cite{atakishi1}.
In the upcoming section we will show that the failure is due to the breakdown
of the self-adjointness of the naive dual polynomial system, and present
the prescription to derive the correct orthogonality relation, which will
be termed as the two component Hamiltonian formalism.

Before closing this section, let us emphasise again that the simple dual
correspondence relation \eqref{Duality} is the consequence of the universal
boundary conditions \eqref{etazerocond}, \eqref{Pzero}, \eqref{Qxzero}.
As will be shown in subsequent sections, the universal boundary conditions
cannot be imposed for those polynomials having Jackson integral measures.
For them, {\em e.g.\/} the big $q$-Jacobi family, and the $q$-Meixner
($q$-Charlier), the dual correspondence relation contains an extra constant
factor reflecting the boundary conditions.
The duality does not exist for polynomials defined on the full integer
lattice, {\em i.e.\/} the discrete $q$-Hermite $\II$ and the $q$-Laguerre
to be discussed in section \ref{sec:dqHII}.

\section{Dual $q$-Meixner and Two Component Hamiltonian}
\label{sec:dqMeix}

\subsection{Naive dual $q$-Meixner polynomials}
\label{sec:naive}

Here we recapitulate the basic data of the $q$-Meixner polynomials (to be
abbreviated as $q$M hereafter) presented in I:
\begin{align}
  &q^{\bm{\lambda}}=(b,c),\quad
  \bm{\delta}=(1,-1),\quad \kappa=q,\quad
  0<b<q^{-1}\ \text{or}\,\ 0<-b<c^{-1}q^{-1},\quad c>0,
  \label{qMpara}\\
  &B(x;\bm{\lambda})=cq^x(1-bq^{x+1}),\quad
  D(x;\bm{\lambda})=(1-q^x)(1+bcq^x),
  \label{qMB&D}\\
  &\mathcal{E}(n)=1-q^n,\quad\eta(x)=q^{-x}-1,
  \label{qMsin}\\
  &\check{P}_n(x;\bm{\lambda})
  =P_n\bigl(\eta(x);\bm{\lambda}\bigr)={}_2\phi_1\Bigl(
  \genfrac{}{}{0pt}{}{q^{-n},\,q^{-x}}{bq}\!\Bigm|\!q\,;-c^{-1}q^{n+1}\Bigr)
  =M_n(q^{-x};b,c;q),
  \label{qMpol}\\
  &\phi_0(x;\bm{\lambda})^2
  =c^xq^{\frac12x(x-1)}\frac{(bq;q)_x}{(q,-bcq;q)_x}
  \quad\bigl(\Leftarrow\phi_0(0;\bm{\lambda})=1\bigr),
  \label{qMphi0}\\
  &d_n(\bm{\lambda})^2
  =\frac{q^n(bq;q)_n}{(q,-c^{-1}q;q)_n}\times d_0(\bm{\lambda})^2,\quad
  d_0(\bm{\lambda})^2=\frac{(-bcq;q)_{\infty}}{(-c;q)_{\infty}},
  \label{qMdn}\\[4pt]
  &A_n(\bm{\lambda})=-cq^{-2n-1}(1-bq^{n+1}),\quad
  C_n(\bm{\lambda})=-q^{-2n}(1-q^n)(q^n+c).
  \label{3coef}
\end{align}
The last entry \eqref{3coef} is the coefficients of the three term
recurrence relation of the $q$M ($n=0,1,\ldots$),
\begin{equation}
  \eta(x)\check{P}_n(x;\bm{\lambda})
  =A_n(\bm{\lambda})\check{P}_{n+1}(x;\bm{\lambda})
  -\bigl(A_n(\bm{\lambda})+C_n(\bm{\lambda})\bigr)
  \check{P}_n(x;\bm{\lambda})
  +C_n(\bm{\lambda})\check{P}_n(x;\bm{\lambda}).
  \label{3rel}
\end{equation}
Note that the potentials are positive for the extended range of $b$ as
in \eqref{qMpara}.

Based on these data, we reported that the ``dual $q$-Meixner polynomial''
is given by
\begin{equation}
  P^{\text{d}}_n\bigl(\eta^{\text{d}}(x);\bm{\lambda}\bigr)
  ={}_2\phi_1\Bigl(
  \genfrac{}{}{0pt}{}{q^{-n},\,q^{-x}}{bq}\!\Bigm|\!q\,;-c^{-1}q^{x+1}\Bigr),
  \label{dualqmeixner0}
\end{equation}
which satisfies the difference equation ($n=0,1,\ldots$)
\begin{align}
  &\bigl(B^{\text{d}}(x;\bm{\lambda})(1-e^{\partial})
  +D^{\text{d}}(x;\bm{\lambda})(1-e^{-\partial})\bigr)
  P^{\text{d}}_n\bigl(\eta^{\text{d}}(x);\bm{\lambda}\bigr)
  =\mathcal{E}^{\text{d}}(n)
  P^{\text{d}}_n\bigl(\eta^{\text{d}}(x);\bm{\lambda}\bigr),
  \label{dualdif}\\
  &\qquad
  B^{\text{d}}(x;\bm{\lambda})=cq^{-2x-1}(1-bq^{x+1}),\quad
  D^{\text{d}}(x;\bm{\lambda})=q^{-2x}(1-q^x)(q^x+c),
  \label{dualqmeixner1}\\
  &\qquad
  \mathcal{E}^{\text{d}}(n)=q^{-n}-1,\quad\eta^{\text{d}}(x)=1-q^x,\\
  &\qquad
  \phi^{\text{d}}_0(x;\bm{\lambda})^2
  =\phi^{\text{d}}_0(0;\bm{\lambda})^2
  \frac{q^x(bq;q)_x}{(q,-c^{-1}q;q)_x},
  \label{dualqmeixner2}\\
  &\qquad A^{\text{d}}_n(\bm{\lambda})=-cq^n(1-bq^{n+1}),\quad
  C^{\text{d}}_n(\bm{\lambda})=-(1-q^n)(1+bcq^n).
  \label{dualqmeixner3}
\end{align}
The difference equation \eqref{dualdif} is just rewriting of the above three
term recurrence relation \eqref{3rel} in terms of the dual correspondence
\eqref{duality1}--\eqref{duality2}.
But the claim of this ``dual $q$-Meixner polynomial'' is not totally correct,
since the completeness relation or the dual orthogonality relation
\eqref{Qnortho}
\begin{equation}
  \sum_{n=0}^{\infty} d_n(\bm{\lambda})^2
  P_n\bigl(\eta(x);\bm{\lambda}\bigr)P_n\bigl(\eta(y);\bm{\lambda}\bigr)
  =\frac{\delta_{xy}}{\phi_0(x;\bm{\lambda})^2},
  \label{dortho-qM}
\end{equation}
{\em does not hold}.
The non-orthogonality can be easily seen for certain parameter ranges.
For $x=0$, $y=1$, we obtain $P_n(0;\bm{\lambda})=1$,
$P_n(\eta(1);\bm{\lambda})=1+c^{-1}(q^n-1)/(1-bq)$.
In the parameter range $A\eqdef 1-{c}^{-1}(1-bq)^{-1}>0$,
the left hand side of \eqref{dortho-qM} does not vanish,
l.h.s.$=\sum\limits_{n=0}^{\infty}d_n(\bm{\lambda})^2
\bigl(A(1-q^n)+q^n\bigr)>0$.

This is due to the {\em breakdown of the self-adjointness} of the naive
dual Hamiltonian
\begin{align}
  &\mathcal{H}^{\text{d}}(\bm{\lambda})
  =\mathcal{A}^{\text{d}}(\bm{\lambda})^\dagger
  \mathcal{A}^{\text{d}}(\bm{\lambda}),\quad
  \mathcal{A}^{\text{d}}(\bm{\lambda})
  \eqdef\sqrt{B^{\text{d}}(x;\bm{\lambda})}
  -e^{\partial}\sqrt{D^{\text{d}}(x;\bm{\lambda})},
  \label{dqMHd1}\\
  &\mathcal{H}^{\text{d}}(\bm{\lambda})\phi^{\text{d}}_n(x;\bm{\lambda})
  =\mathcal{E}^{\text{d}}(n)\phi^{\text{d}}_n(x;\bm{\lambda})
  \ \ (n=0,1,\ldots),\\
  &\phi^{\text{d}}_n(x;\bm{\lambda})
  \eqdef\phi^{\text{d}}_0(x;\bm{\lambda})
  P^{\text{d}}_n\bigl(\eta^{\text{d}}(x);\bm{\lambda}\bigr).
  \label{dqMHd3}
\end{align}
The generic element of the corresponding Hilbert space has the form
$f(x)=\phi^{\text{d}}_0(x;\bm{\lambda})\mathcal{P}(\eta^{\text{d}}(x))$,
$\norm{f}<\infty$, with a smooth function $\mathcal{P}$.
At large $N$, we have
\begin{align*}
  &\mathcal{P}\bigl(\eta^{\text{d}}(N)\bigr)\simeq
  \mathcal{P}(1)-\mathcal{P}'(1)q^N+O(q^{2N}),\\
  &\ \ \Rightarrow
  \mathcal{P}\bigl(\eta^{\text{d}}(N+1)\bigr)
  \mathcal{Q}\bigl(\eta^{\text{d}}(N)\bigr)
  -\mathcal{P}\bigl(\eta^{\text{d}}(N)\bigr)
  \mathcal{Q}\bigl(\eta^{\text{d}}(N+1)\bigr)
  \simeq(1-q)\text{W}[\mathcal{P},\mathcal{Q}](1)q^N+O(q^{2N}),\\
  &\phi^{\text{d}}_0(N)^2\simeq
  q^N\bigl(1+O(q^N)\bigr)\times\text{const},\quad
  B^{\text{d}}(N)\simeq cq^{-2N-1}\bigl(1+O(q^N)\bigr),
\end{align*}
in which $\text{W}[f,g](x)=f(x)g'(x)-f'(x)g(x)$ is the Wronskian.
These show that the right hand side of the criterion formula \eqref{criter1}
is finite in the $N\to\infty$ limit.
The self-adjointness is broken mainly by the very strong increase
($\sim\!q^{-2x}$) of the potential functions $B^{\text{d}}(x)$,
$D^{\text{d}}(x)$ at $x\to \infty$ \eqref{dualqmeixner1} and the almost
constant behaviour of the polynomial
$P^{\text{d}}_n(\eta^{\text{d}}(x);\bm{\lambda})$ at $x\to\infty$.

\subsection{Two component Hamiltonian system}
\label{sec:binham}

A clue for the recovery of the self-adjointness of the naive dual $q$-Meixner
(d$q$M) Hamiltonian system \eqref{dualqmeixner1},
\eqref{dqMHd1}--\eqref{dqMHd3} is the fact that it is accompanied by
{\em another Hamiltonian\/} sharing the same eigenvalues and the corresponding
eigenvectors which are obtained by rescaling the sinusoidal coordinate.
A close look at the potential functions \eqref{dualqmeixner1} of the naive
Hamiltonian shows that the rescaling of the variable
\begin{equation}
  q^x\to-cq^x,
  \label{rescale}
\end{equation}
provides another valid Hamiltonian. For later discussion, let us denote
the original naive Hamiltonian and its potential functions etc., by
superscript $(+)$ and those of the rescaled one by $(-)$:
\begin{align}
  &\mathcal{H}^{(\pm)}={\mathcal{A}^{(\pm)}}^{\dagger}\mathcal{A}^{(\pm)},
  \qquad \mathcal{A}^{(\pm)}\eqdef
  \sqrt{B^{(\pm)}(x)}-e^{\partial}\sqrt{D^{(\pm)}(x)},
  \label{dqMpm1}\\
  &B^{(+)}(x)=cq^{-2x-1}(1-bq^{x+1}),
  \qquad\ D^{(+)}(x)=q^{-2x}(1-q^x)(q^x+c),
  \label{dqMpm2}\\
  &\phantom{B^{(+)}(x)\ }
  \big\Downarrow\qquad\qquad\quad\ q^x\to-cq^x\ \ \qquad\qquad\big\Downarrow\n
  &B^{(-)}(x)=c^{-1}q^{-2x-1}(1+bcq^{x+1}),\quad
  D^{(-)}(x)=c^{-1}q^{-2x}(1-q^{x})(1+cq^x),
  \label{dqMpm3}\\
  &\check{P}^{(+)}_n(x)={}_2\phi_1\Bigl(
  \genfrac{}{}{0pt}{}{q^{-n},\,q^{-x}}{bq}\!\Bigm|\!q\,;-c^{-1}q^{x+1}\Bigr),
  \ \ \check{P}^{(-)}_n(x)={}_2\phi_1\Bigl(
  \genfrac{}{}{0pt}{}{q^{-n},\,-c^{-1}q^{-x}}{bq}\!\Bigm|\!q\,;q^{x+1}\Bigr),
  \label{dqMpm4}\\
  &\mathcal{H}^{(\pm)}\phi^{(\pm)}_n(x)
  =\mathcal{E}^{\text{d}}(n)\phi^{(\pm)}_n(x),\quad
  \phi^{(\pm)}_n(x)=\phi^{(\pm)}_0(x)\check{P}^{(\pm)}_n(x)
  \ \ (n=0,1,\ldots),
  \label{dqMpm5}\\
  &\mathcal{A}^{(\pm)}\phi^{(\pm)}_0(x)=0,\quad
  \phi^{(\pm)}_0(x)=\phi^{(\pm)}_0(0)
  \prod_{y=0}^{x-1}\sqrt{\frac{B^{(\pm)}(y)}{D^{(\pm)}(y+1)}},
  \label{dqMpm6}\\
  &\widetilde{\mathcal{H}}^{(\pm)}=\phi^{(\pm)}_0(x)^{-1}\circ
  \mathcal{H}^{(\pm)}\circ\phi^{(\pm)}_0(x)
  =B^{(\pm)}(x)(1-e^{\partial})+D^{(\pm)}(x)(1-e^{-\partial}).
  \label{dqMpm7}
\end{align}
Here we have suppressed the parameter $\bm{\lambda}$ dependence for
simplicity. Except for the ratio of $\phi^{(\pm)}_0(0)$ \eqref{dqMpm6},
every quantity in the $(-)$ system is determined from that of the $(+)$
system by rescaling \eqref{rescale}. The positivity of $B^{(-)}(x)$,
$D^{(-)}(x)$ and the boundary condition of $D^{(-)}(x)$ are obviously
satisfied.
The triangularity of $\widetilde{\mathcal{H}}^{(-)}$ \eqref{dqMpm7} is
trivial, except for the vanishing of $q^{-x}$ term in
$\widetilde{\mathcal{H}}^{(-)}q^x$.
This is inherited from $\widetilde{\mathcal{H}}^{(+)}$, as the cancellation
is achieved by the coefficients of $q^{-2x}$ terms in $B^{(+)}(x)$,
$D^{(+)}(x)$. The situation is unchanged by the rescaling \eqref{rescale}.
The eigenvalues of $\widetilde{\mathcal{H}}^{(\pm)}$ \eqref{dqMpm7} are the
same, since they are determined by the $(q^x)^0$ terms in the potential
functions, which are unchanged by the rescaling \eqref{rescale}.
It is easy to verify that both Hamiltonians $\mathcal{H}^{(\pm)}$ are shape
invariant.

Now we have two closely related Hamiltonian systems $\mathcal{H}^{(\pm)}$,
each of which is not self-adjoint separately.
We will demonstrate that the self-adjointness can be recovered in the
combined system, to be called two component Hamiltonian system, by
adjusting the ratio of $\phi^{(\pm)}_0(0)$ \eqref{dqMpm6}.
Since the rescaling factor is always negative, {\em e.g.} $q^x\to -cq^x$,
this approach leads naturally to the orthogonality measure of the
polynomials of {\em Jackson integral\/} type.
The Jackson integral is defined by \cite{koeswart}
\begin{equation*}
  \int_0^{\alpha}\!\!d_qy\,f(y)\eqdef(1-q)\alpha
  \sum_{k=0}^{\infty}f(\alpha q^k)q^k,\quad
  \int_{\alpha}^{\beta}\!\!d_qy\,f(y)\eqdef
  \int_0^{\beta}\!\!d_qy\,f(y)-\int_0^{\alpha}\!\!d_qy\,f(y).
\end{equation*}

\bigskip
Let us introduce appropriate notation for generic two component Hamiltonian
systems.
The vector $\bm{f}$ and its inner product are
\begin{align}
  \bm{f}(x)&=\genfrac{(}{)}{0pt}{}{f^{(+)}(x)}{f^{(-)}(x)}
  \quad(x\in\mathbb{Z}_{\geq 0}),
  \label{vec2}\\
  \bin{\bm{f}}{\bm{g}}&=\bigl(f^{(+)},g^{(+)}\bigr)
  +\bigl(f^{(-)},g^{(-)}\bigr)\n
  &=\lim_{N\to\infty}{\bin{\bm{f}}{\bm{g}}}_N
  =\lim_{N\to\infty}\Bigl(\bigl(f^{(+)},g^{(+)}\bigr)_N
  +\bigl(f^{(-)},g^{(-)}\bigr)_N\Bigr).
  \label{binpro}
\end{align}
The Hamiltonian is a direct sum
\begin{equation}
  \underline{\mathcal{H}}=
  \begin{pmatrix}\mathcal{H}^{(+)}&0\\0&\mathcal{H}^{(-)}\end{pmatrix}.
  \label{dsHam}
\end{equation}
The Schr\"odinger equation is
\begin{equation*}
  \underline{\mathcal{H}}\bm{\phi}_n(x)=\mathcal{E}(n)\bm{\phi}_n(x)\quad
  (n=0,1,\ldots),\quad
  0=\mathcal{E}(0)<\mathcal{E}(1)<\cdots,
\end{equation*}
where the eigenvector
$\bm{\phi}_n(x)=\genfrac{(}{)}{0pt}{}{\phi^{(+)}_n(x)}{\phi^{(-)}_n(x)}$ has
a finite norm.
We define $\underline{B}(x)$, $\underline{D}(x)$, $\underline{\mathcal{A}}$,
etc. as the direct sum as follows:
\begin{align}
  &\underline{B}(x)=\begin{pmatrix}B^{(+)}(x)&0\\0&B^{(-)}(x)\end{pmatrix},
  \quad
  \underline{D}(x)=\begin{pmatrix}D^{(+)}(x)&0\\0&D^{(-)}(x)\end{pmatrix},\n
  &\underline{\mathcal{A}}=
  \begin{pmatrix}\mathcal{A}^{(+)}&0\\0&\mathcal{A}^{(-)}\end{pmatrix},\quad
  \underline{\mathcal{A}}^{\dagger}=
  \begin{pmatrix}\mathcal{A}^{(+)\,\dagger}&0\\
  0&\mathcal{A}^{(-)\,\dagger}\end{pmatrix},\n
  &\underline{\widetilde{\mathcal{H}}}=
  \begin{pmatrix}\widetilde{\mathcal{H}}^{(+)}&0\\
  0&\widetilde{\mathcal{H}}^{(-)}\end{pmatrix},\quad
  \underline{\phi_0}(x)
  =\begin{pmatrix}\phi^{(+)}_0(x)&0\\0&\phi^{(-)}_0(x)\end{pmatrix},\quad
  \underline{e^{\pm\partial}\!\!}\ =
  \begin{pmatrix}e^{\pm\partial}&0\\0&e^{\pm\partial}\end{pmatrix},\n[2pt]
  &\underline{\mathcal{H}}
  =\underline{\mathcal{A}}^{\dagger}\underline{\mathcal{A}}
  =-\sqrt{\underline{B}(x)}\,\underline{e^{\partial}}\sqrt{\underline{D}(x)}
  -\sqrt{\underline{D}(x)}
  \,\underline{e^{-\partial}\!\!}\ \sqrt{\underline{B}(x)}
  +\underline{B}(x)+\underline{D}(x),\n[2pt]
  &\underline{\mathcal{A}}=\sqrt{\underline{B}(x)}
  -\underline{e^{\partial}}\sqrt{\underline{D}(x)},\quad
  \underline{\mathcal{A}}^{\dagger}=\sqrt{\underline{B}(x)}
  -\sqrt{\underline{D}(x)}\,\underline{e^{-\partial}\!\!}\ ,\n[2pt]
  &\underline{\widetilde{\mathcal{H}}}
  =\underline{\phi_0}(x)^{-1}\circ\underline{\mathcal{H}}\circ
  \underline{\phi_0}(x)
  =\underline{B}(x)(1-\underline{e^{\partial}})
  +\underline{D}(x)(1-\underline{e^{-\partial}\!\!}\ ).
  \label{Ht2}
\end{align}
Since $B^{(\pm)}(x)$ and $D^{(\pm)}(x)$ are functions of two sinusoidal
coordinates $\eta^{(\pm)}(x)$ of the form
\begin{equation}
  \eta^{(\pm)}(x)=\alpha^{(\pm)}q^x,\quad\alpha^{(\pm)}\in\mathbb{R},\quad
  \alpha^{(+)}\alpha^{(-)}<0,
  \label{2eta}
\end{equation}
these four input potentials can be expressed by two functions
$B^{\text{J}}(\eta)$ and $D^{\text{J}}(\eta)$ as
\begin{equation}
  B^{(\pm)}(x)=B^{\text{J}}\bigl(\eta^{(\pm)}(x)\bigr),\quad
  D^{(\pm)}(x)=D^{\text{J}}\bigl(\eta^{(\pm)}(x)\bigr),
  \label{tilddef}
\end{equation}
in which J stands for `Jackson'.
It should be stressed that the sinusoidal coordinates $\eta^{(\pm)}(x)$ in the
two component Hamiltonian formalism are not required the universal boundary
condition $\eta(0)=0$ \eqref{etazerocond}, which was imposed in \cite{os14}.

For two vectors $\bm{f}$ with
$f^{(\pm)}(x)=\phi^{(\pm)}_0(x)\check{\mathcal{P}}^{(\pm)}(x)$ and
$\bm{g}$ with $g^{(\pm)}(x)=\phi^{(\pm)}_0(x)\check{\mathcal{Q}}^{(\pm)}(x)$,
the criterion of the self-adjointness \eqref{criter1} now reads
\begin{align}
  0&=\lim_{N\to\infty}\bigl({\bin{\bm{f}}{\underline{\mathcal{H}}\,\bm{g}}}_N
  -{\bin{\underline{\mathcal{H}}\,\bm{f}}{\bm{g}}}_N\bigr)\n
  &=\lim_{N\to\infty}\Bigl(\phi^{(+)}_0(N)^2B^{(+)}(N)
  \bigl(\check{\mathcal{P}}^{(+)}(N+1)\check{\mathcal{Q}}^{(+)}(N)
  -\check{\mathcal{P}}^{(+)}(N)\check{\mathcal{Q}}^{(+)}(N+1)\bigr)\n
  &\ \ \,\qquad
  +\phi^{(-)}_0(N)^2B^{(-)}(N)
  \bigl(\check{\mathcal{P}}^{(-)}(N+1)\check{\mathcal{Q}}^{(-)}(N)
  -\check{\mathcal{P}}^{(-)}(N)\check{\mathcal{Q}}^{(-)}(N+1)\bigr)\Bigr).
  \label{fHgN-HfgN2}
\end{align}
For smooth functions $\mathcal{P}$ and $\mathcal{Q}$, we have the following
large $N$ behaviours:
\begin{align}
  &\quad
  \mathcal{P}(\alpha q^N)\simeq\mathcal{P}(0)+\mathcal{P}'(0)\alpha q^N
  +O(q^{2N}),\n
  &\quad
  \mathcal{P}\bigl(\eta^{(\pm)}(N+1)\bigr)
  \mathcal{Q}\bigl(\eta^{(\pm)}(N)\bigr)
  -\mathcal{P}\bigl(\eta^{(\pm)}(N)\bigr)
  \mathcal{Q}\bigl(\eta^{(\pm)}(N+1)\bigr)\n
  &\simeq(1-q)\text{W}[\mathcal{P},\mathcal{Q}](0)\alpha^{(\pm)}q^N
  +O(q^{2N}).
  \label{largeN:PQ-QP}
\end{align}

Now we are in a position to demonstrate the self-adjointness and to derive
the orthogonality relation for the d$q$M.
Two functions $B^{\text{J}}(\eta)$ and $D^{\text{J}}(\eta)$ are
\begin{equation}
  B^{\text{J}}(\eta)\eqdef\eta^{-2}cq^{-1}(1-bq\eta),\quad
  D^{\text{J}}(\eta)\eqdef\eta^{-2}(1-\eta)(\eta+c).
  \label{qM:BJDJ}
\end{equation}
The explicit forms of the ground state eigenfunctions $\phi^{(\pm)}_0(x)$ are
\begin{equation}
  \phi^{(+)}_0(x)^2
  =\phi^{(+)}_0(0)^2\frac{q^x(bq;q)_x}{(q,-c^{-1}q;q)_x},\quad
  \phi^{(-)}_0(x)^2=\phi^{(-)}_0(0)^2\frac{q^x(-bcq;q)_x}{(q,-cq;q)_x}.
  \label{dqnpm0}
\end{equation}
Let us introduce new sinusoidal coordinates $\eta^{(\pm)}(x)$,
as the original and rescaled $q^x$:
\begin{equation*}
  \eta^{(+)}(x)\eqdef q^x>0,\quad\eta^{(-)}(x)\eqdef -cq^x<0.
\end{equation*}
In terms of these sinusoidal coordinates, the above ground state
eigenfunctions \eqref{dqnpm0} have a unified expression:
\begin{align*}
  \phi^{(\pm)}_0(x)^2&=\pm A^{(\pm)}\eta^{(\pm)}(x)
  \frac{(q\eta^{(\pm)}(x),-c^{-1}q\eta^{(\pm)}(x);q)_\infty}
  {(bq\eta^{(\pm)}(x);q)_\infty},\\
  &A^{(+)}\eqdef\frac{\phi^{(+)}_0(0)^2(bq;q)_\infty}
  {(q,-c^{-1}q;q)_\infty},\quad
  A^{(-)}\eqdef\frac{\phi^{(-)}_0(0)^2(-bcq;q)_\infty}{c(q,-cq;q)_\infty},\\
  \phi^{(\pm)}_0(N)^2&\simeq\pm A^{(\pm)}\eta^{(\pm)}(N)
  \bigl(1+O(q^N)\bigr)\ \ (N\to\infty).
\end{align*}
Then the r.h.s.\!\! of the criterion \eqref{fHgN-HfgN2} reads for
$\check{\mathcal{P}}^{(\pm)}(x)=\mathcal{P}(\eta^{(\pm)}(x))$ and
$\check{\mathcal{Q}}^{(\pm)}(x)=\mathcal{Q}(\eta^{(\pm)}(x))$:
\begin{align*}
  &\lim_{N\to\infty}\Bigl(
  A^{(+)}q^Ncq^{-2N-1}(1-q)\text{W}[\mathcal{P},\mathcal{Q}](0)q^N
  +A^{(-)}cq^Nc^{-1}q^{-2N-1}(1-q)
  \text{W}[\mathcal{P},\mathcal{Q}](0)(-c)q^N\Bigr)\\
  &=\bigl(A^{(+)}-A^{(-)}\bigr)cq^{-1}(1-q)
  \text{W}[\mathcal{P},\mathcal{Q}](0),
\end{align*}
implying that the self-adjointness of the two component Hamiltonian system
is achieved by the choice 
\begin{equation*}
  A^{(+)}=A^{(-)}=1,\quad
  \phi^{(+)}_0(0)^2\eqdef\frac{(q,-c^{-1}q;q)_\infty}{(bq;q)_\infty},\quad
  \phi^{(-)}_0(0)^2\eqdef\frac{c(q,-cq;q)_\infty}{(-bcq;q)_\infty}.
\end{equation*}
The orthogonality relation can also be expressed in terms of the Jackson
integral ($n,m=0,1,\ldots$):
\begin{align} 
  \bin{\bm{\phi}_n}{\bm{\phi}_m}
  &=\sum_{\epsilon=\pm}\sum_{x=0}^{\infty}
  \phi^{(\epsilon)}_0(x;\bm{\lambda})^2
  P_n\bigl(\eta^{(\epsilon)}(x;\bm{\lambda});\bm{\lambda}\bigr)
  P_m\bigl(\eta^{(\epsilon)}(x;\bm{\lambda});\bm{\lambda}\bigr)
  =\frac{\delta_{nm}}{d_n(\bm{\lambda})^2}\\
  &=\sum_{k=0}^{\infty}
  \frac{(q^{k+1},-c^{-1}q^{k+1};q)_\infty}{(bq^{k+1};q)_\infty}
  P_n(q^k;\bm{\lambda})P_m(q^k;\bm{\lambda})q^k\n
  &\quad+c\sum_{k=0}^{\infty}
  \frac{(q^{k+1},-cq^{k+1};q)_\infty}{(-bcq^{k+1};q)_\infty}
  P_n(-cq^k;\bm{\lambda})P_m(-cq^k;\bm{\lambda})q^k\n
  &=\frac1{1-q}\int_{-c}^1\!\!d_qy\,
  \frac{(qy,-c^{-1}qy;q)_\infty}{(bqy;q)_\infty}
  P_n(y;\bm{\lambda})P_m(y;\bm{\lambda}),\\
  d_n(\bm{\lambda})^2&=c^nq^{\frac12n(n-1)}
  \frac{(bq;q)_n}{(q,-bcq;q)_n}\times d_0(\bm{\lambda})^2,\quad
  d_0(\bm{\lambda})^2=
  \frac{(bq,-bcq;q)_{\infty}}{(q,-c,-c^{-1}q;q)_{\infty}}.
  \label{qM:dualdn}
\end{align}
Note that $d_n(\bm{\lambda})^2/d_0(\bm{\lambda})^2$ in \eqref{qM:dualdn}
and $\phi_0(x;\bm{\lambda})^2$ in \eqref{qMphi0} have the same form.
Here the polynomial $P_n(y;\bm{\lambda})$ is obtained from the original
expression \eqref{dualqmeixner0} by the replacement $q^x\to y$:
\begin{equation*}
  P_n(y;\bm{\lambda})={}_2\phi_1\Bigl(
  \genfrac{}{}{0pt}{}{q^{-n},\,y^{-1}}{bq}\!\Bigm|\!q\,;-c^{-1}qy\Bigr).
\end{equation*}
The three term recurrence relation of the above dual $q$-Meixner
polynomials is
\begin{equation}
  (1-y)P_n(y)=A^{\text{d}}_nP_{n+1}(y)
  -(A^{\text{d}}_n+C^{\text{d}}_n)P_n(y)
  +C^{\text{d}}_nP_{n-1}(y),
  \label{threey}
\end{equation}
in which $A^{\text{d}}_n$ and $C^{\text{d}}_n$ are given in
\eqref{dualqmeixner3}.
This is related to the big $q$-Laguerre polynomial by the rescaling of the
variable $y$ and redefinition of the parameters.

\subsection{Dual $q$-Charlier polynomials}
\label{sec:dqChar}

The naive dual $q$-Charlier polynomials presented in I do not satisfy the
orthogonality relation, due to the breakdown of the self-adjointness
condition. As for the $q$-Meixner system, the two component Hamiltonian
formulation offers the remedy.

The Hamiltonian system for the $q$-Charlier polynomials ($q$C) is obtained
from that of $q$M by putting $b=0$ and change of the parameter $c\to a$.
The basic data presented in I are:
\begin{align}
  &q^{\bm{\lambda}}=a,\quad
  \bm{\delta}=-1,\quad\kappa=q,\quad a>0,\n
  &B(x;\bm{\lambda})=aq^x,\quad D(x)=1-q^x,\quad
  \mathcal{E}(n)=1-q^n,\quad\eta(x)=q^{-x}-1,
  \label{qCB&D}\\
  &\check{P}_n(x;\bm{\lambda})
  =P_n\bigl(\eta(x);\bm{\lambda}\bigr)={}_2\phi_1\Bigl(
  \genfrac{}{}{0pt}{}{q^{-n},\,q^{-x}}{0}\!\Bigm|\!q\,;-a^{-1}q^{n+1}\Bigr)
  =C(q^{-x};a;q),
  \label{qCp}\\[2pt]
  &\phi_0(x;\bm{\lambda})^2=\frac{a^xq^{\frac12x(x-1)}}{(q;q)_x}
  \quad\bigl(\Leftarrow\phi_0(0;\bm{\lambda})=1\bigr),
  \label{qCphi0}\\
  &d_n(\bm{\lambda})^2=\frac{q^n}{(q,-a^{-1}q;q)_n}
  \times d_0(\bm{\lambda})^2,\quad
  d_0(\bm{\lambda})^2=\frac{1}{(-a;q)_{\infty}},
  \label{qCdn}\\[2pt]
  &A_n(\bm{\lambda})=-aq^{-2n-1},\quad
  C_n(\bm{\lambda})=-q^{-2n}(1-q^n)(q^n+a).
  \label{qC3coef}
\end{align}

The $(+)$ part of the dual $q$-Charlier (d$q$C) system is
\begin{align}
  &\check{P}^{(+)}_n(x;\bm{\lambda})
  =P_n\bigl(\eta^{(+)}(x);\bm{\lambda}\bigr)
  ={}_2\phi_1\Bigl(
  \genfrac{}{}{0pt}{}{q^{-n},\,q^{-x}}{0}\!\Bigm|\!q\,;-a^{-1}q^{x+1}\Bigr),
  \label{dualqchar-1}\\
  &B^{(+)}(x;\bm{\lambda})=aq^{-2x-1},\quad
  D^{(+)}(x;\bm{\lambda})=q^{-2x}(1-q^x)(q^x+a),
  \label{dualqchar1}\\
  &\mathcal{E}(n)=q^{-n}-1,\quad\eta^{(+)}(x)=q^x,\\
  &\phi^{(+)}_0(x;\bm{\lambda})^2
  =(q,-a^{-1}q;q)_{\infty}\frac{q^x}{(q,-a^{-1}q;q)_x}
  =q^x(q^{x+1},-a^{-1}q^{x+1};q)_\infty,
  \label{dualqchar2}\\
  &A_n(\bm{\lambda})=-aq^n,\quad C_n=-(1-q^n).
  \label{dualqchar3}
  \end{align}
The self-adjointness of the $(+)$ part Hamiltonian is broken as in the
d$q$M case.
The quantities of the $(-)$ part are obtained by the replacement
$q^x\to -aq^x$:
\begin{align}
  &\check{P}^{(-)}_n(x;\bm{\lambda})
  =P_n\bigl(\eta^{(-)}(x);\bm{\lambda}\bigr)
  ={}_2\phi_1\Bigl(
  \genfrac{}{}{0pt}{}{q^{-n},\,-a^{-1}q^{-x}}{0}\!\Bigm|\!q\,;q^{x+1}\Bigr),
  \ \ \eta^{(-)}(x)=-aq^x,
  \label{dualqchar-}\\
  &B^{(-)}(x;\bm{\lambda})=a^{-1}q^{-2x-1},\quad
  D^{(-)}(x;\bm{\lambda})=a^{-1}q^{-2x}(1-q^x)(1+aq^x),
  \label{dualqchar1-}\\
  &\phi^{(-)}_0(x;\bm{\lambda})^2
  =a(q,-aq;q)_{\infty}\frac{q^x}{(q,-aq;q)_x}
  =aq^x(q^{x+1},-aq^{x+1};q)_\infty.
  \label{dualqchar2-}
\end{align}
With the ground state eigenfunctions $\phi^{(\pm)}_0(x;\bm{\lambda})^2$
\eqref{dualqchar2}, \eqref{dualqchar2-}, the self-adjointness is achieved
and the orthogonality relation is
\begin{align}
  \bin{\bm{\phi}_n}{\bm{\phi}_m}
  &=\sum_{\epsilon=\pm}\sum_{x=0}^{\infty}
  \phi^{(\epsilon)}_0(x;\bm{\lambda})^2
  P_n\bigl(\eta^{(\epsilon)}(x;\bm{\lambda});\bm{\lambda}\bigr)
  P_m\bigl(\eta^{(\epsilon)}(x;\bm{\lambda});\bm{\lambda}\bigr)
  =\frac{\delta_{nm}}{d_n(\bm{\lambda})^2}\\
  &=\sum_{k=0}^{\infty}{(q^{k+1},-a^{-1}q^{k+1};q)_\infty}
  P_n(q^k;\bm{\lambda})P_m(q^k;\bm{\lambda})q^k\n
  &\quad+a\sum_{k=0}^{\infty}{(q^{k+1},-aq^{k+1};q)_\infty}
  P_n(-aq^k;\bm{\lambda})P_m(-aq^k;\bm{\lambda})q^k\n
  &=\frac1{1-q}\int_{-a}^1\!\!d_qy\,{(qy,-a^{-1}qy;q)_\infty}
  P_n(y;\bm{\lambda})P_m(y;\bm{\lambda}),\\
  d_n(\bm{\lambda})^2&=\frac{a^nq^{\frac12n(n-1)}}{(q;q)_n}
  \times d_0(\bm{\lambda})^2,\quad
  d_0(\bm{\lambda})^2=\frac{1}{(q,-a,-a^{-1}q;q)_{\infty}},
  \label{qC:dualdx}\\
  P_n(y;\bm{\lambda})&={}_2\phi_1\Bigl(
  \genfrac{}{}{0pt}{}{q^{-n},\,y^{-1}}{0}\!\Bigm|\!q\,;-a^{-1}qy\Bigr),
  \label{qC:Pny}
\end{align}
in which $\{P_n(y)\}$ satisfy the three term recurrence relation
\eqref{threey} with $A^{\text{d}}_n$ and $C^{\text{d}}_n$ replaced by
$A_n(\bm{\lambda})$, $C_n$ in \eqref{dualqchar3}.
This is related to the Al-Salam-Carlitz $\I$ polynomial by the redefinition
of the parameter.

\section{Orthogonal Polynomials and Jackson Integral}
\label{sec:Jackson}

In this section we discuss  polynomials having orthogonality measures of
Jackson integral type and their duals.
As exemplified by the d$q$M in the previous section, the Jackson integral
measures are closely related to the two component Hamiltonian formalism.
The discrete $q$-Hermite polynomial II has some exceptional features,
which will be discussed separately in \S\,\ref{sec:dqHII}.

\subsection{General structure}
\label{sec:genstr}

Let us begin with the general structures of the two component Hamiltonian
systems. The most generic potential functions have the following form
(The superscript J stands for `Jackson,' see \eqref{tilddef})
\begin{equation}
  B^{\text{J}}(\eta)=\eta^{-2}p_1(\eta),\quad
  D^{\text{J}}(\eta)=\eta^{-2}p_2(\eta),
  \label{p1p2}
\end{equation}
in which $p_1(\eta)$ is a polynomial in $\eta\propto q^x$ \eqref{2eta} of
at most degree 2, whereas $p_2(\eta)$ is a polynomial of exact degree 2 in
$\eta$. Thus the system has at most six real parameters. Among them, one
corresponding to the overall normalisation of the eigenvalues
$\{\mathcal{E}(n)\}$ can be fixed to 1.
This constrains the degree two terms of $p_1(\eta)$ and $p_2(\eta)$.
Since the rescaling of $\eta$ does not change the form \eqref{p1p2},
one parameter can be reduced.
The boundary condition of $D^{(\pm)}(0)=0$ determines $\alpha^{(\pm)}$ in
$\eta^{(\pm)}(x)$ \eqref{2eta}.
The triangularity of $\widetilde{\mathcal{H}}^{(\pm)}$ provides the third
constraint, requiring
\begin{equation}
  \widetilde{\mathcal{H}}^{(\pm)}\eta=c_1\eta+c_0+0\times \eta^{-1}.
  \label{gentri}
\end{equation}
This gives the constraint on the constant terms of $p_1(\eta)$ and
$p_2(\eta)$, as $0=p_1(0)(1-q)+p_2(0)(1-q^{-1})$.
With this condition, the triangularity of $\widetilde{\mathcal{H}}^{(\pm)}$
with respect to the basis of $\{1,\eta, \eta^{2}\ldots, \eta^{n}\}$ is obvious.
The strong increase $\sim q^{-2x}$ of the potentials at $x\to\infty$ and
the corresponding eigenpolynomials in $q^x$ require the two component
Hamiltonian formalism for self-adjointness.
The positivity of $B^{(\pm)}(x)$ and $D^{(\pm)}(x)$ for
$x\in\mathbb{Z}_{\ge0}$ restricts the ranges of the parameters.
Thus the most generic two component Hamiltonian systems have three
independent parameters.
In contrast, the systems with the least number of parameters correspond to
the case $p_1(\eta)=const.$ and the d$q$C discussed in \S\,\ref{sec:dqChar}
belong to this class.
When these conditions are satisfied, the solvability of the Hamiltonian is
guaranteed independent of the parametrisation.
However, very special parametrisation is required for the shape invariance
\eqref{shapeinv1}--\eqref{shapeinv1cond2} which ensures the existence of
explicit expressions of the polynomials of the Rodrigues type \eqref{univrod}.
Choosing specific parametrisation of the two polynomials $p_1(\eta)$ and
$p_2(\eta)$ in \eqref{p1p2} is tantamount to discussion of the corresponding
specific polynomial system, which will be presented subsequently.
Before closing this subsection, let us mention that any higher degree
generalisation of the above input \eqref{p1p2}, for example,
\begin{equation*}
  B^{\text{J}}(\eta)=\eta^{-3}\bar{p}_1(\eta),\quad
  D^{\text{J}}(\eta)=\eta^{-3}\bar{p}_2(\eta),
\end{equation*}
with cubic polynomials $\bar{p}_1$ and $\bar{p}_2$, cannot work.
The triangularity of $\widetilde{\mathcal{H}}^{(\pm)}$ is broken, since
the two conditions $\widetilde{\mathcal{H}}^{(\pm)}\eta
=c_1\eta+c_0+0\times \eta^{-1}+0\times \eta^{-2}$ and
$\widetilde{\mathcal{H}}^{(\pm)}\eta^2
=c'_2\eta^2+c'_1\eta+c'_0+0\times \eta^{-1}$ are incompatible.

\subsection{Big $q$-Jacobi}
\label{sec:bqJ}

The big  $q$-Jacobi (to be abbreviated as b$q$J hereafter) system has
three real parameters on top of $q$ and it is the most generic member
having the Jackson integral measure. Its data are:
\begin{align}
  &q^{\bm{\lambda}}=(a,b,c),\quad
  \bm{\delta}=(1,1,1),\quad\kappa=q^{-1},\quad
  0<a<q^{-1},\quad 0<b<q^{-1},\quad c<0,
  \label{bqjpararange}\\
  &\mathcal{E}(n;\bm{\lambda})\eqdef(q^{-n}-1)(1-abq^{n+1}),\quad
  \eta^{(+)}(x;\bm{\lambda})\eqdef aq^{x+1},\quad
  \eta^{(-)}(x;\bm{\lambda})\eqdef cq^{x+1},
  \label{bqjeta}\\
  &B^{\text{J}}(\eta;\bm{\lambda})\eqdef\eta^{-2}aq(1-\eta)(b\eta-c),\quad
  D^{\text{J}}(\eta;\bm{\lambda})\eqdef\eta^{-2}(aq-\eta)(\eta-cq).
  \label{BtDtEn}
\end{align}
It is easy to verify the eigenvalues $\mathcal{E}(n;\bm{\lambda})$
\eqref{bqjeta} and the triangularity \eqref{gentri}.
To be more explicit, the potential functions \eqref{tilddef} are
\begin{align}
  &B^{(+)}(x;\bm{\lambda})=-a^{-1}cq^{-2x-1}(1-aq^{x+1})(1-abc^{-1}q^{x+1}),\n
  &D^{(+)}(x;\bm{\lambda})=q^{-2x}(1-q^x)(q^x-a^{-1}c),
  \label{bqJB+}\\
  &B^{(-)}(x;\bm{\lambda})=-ac^{-1}q^{-2x-1}(1-cq^{x+1})(1-bq^{x+1}),\n
  &D^{(-)}(x;\bm{\lambda})=q^{-2x}(1-q^x)(q^x-ac^{-1}).
  \label{bqJB-}
\end{align}

For the parameter range \eqref{bqjpararange}, the positivity and the
boundary conditions are satisfied
\begin{equation*}
  B^{(\pm)}(x;\bm{\lambda})>0\ \ (x\geq 0),\quad
  D^{(\pm)}(x;\bm{\lambda})>0\ \ (x\geq 1),\quad
  D^{(\pm)}(0;\bm{\lambda})=0.
\end{equation*}
It is straightforward to verify that the shape invariance conditions
\eqref{shapeinv1cond1}--\eqref{shapeinv1cond2} are satisfied and that
the same eigenvalues as above \eqref{bqjeta} are obtained by
\eqref{spectrumform}.
The sinusoidal coordinates $\eta^{(\pm)}(x;\bm{\lambda})$ \eqref{bqjeta}
satisfy the following relations
\begin{align}
  \eta^{(\pm)}(x+1;\bm{\lambda})=q\eta^{(\pm)}(x;\bm{\lambda})
  &=\eta^{(\pm)}(x;\bm{\lambda}+\bm{\delta}),
  \label{etaprop}\\
  \eta^{(\pm)}(x+1;\bm{\lambda})+\eta^{(\pm)}(x-1;\bm{\lambda})
  &=(q+q^{-1})\eta^{(\pm)}(x;\bm{\lambda}),
  \label{etaprop1}\\
  \eta^{(\pm)}(x+1;\bm{\lambda})\eta^{(\pm)}(x-1;\bm{\lambda})
  &=\eta^{(\pm)}(x;\bm{\lambda})^2.
  \label{etaprop2}
\end{align}

The big $q$-Jacobi polynomial $P_n(\eta;\bm{\lambda})$
($n\in\mathbb{Z}_{\geq 0}$), which is a degree $n$ polynomial in $\eta$
satisfying the second order difference equation,
\begin{equation}
  B^{\text{J}}(\eta;\bm{\lambda})
  \bigl(P_n(\eta;\bm{\lambda})-P_n(q\eta;\bm{\lambda})\bigr)
  +D^{\text{J}}(\eta;\bm{\lambda})
  \bigl(P_n(\eta;\bm{\lambda})-P_n(q^{-1}\eta;\bm{\lambda})\bigr)
  =\mathcal{E}(n;\bm{\lambda})P_n(\eta;\bm{\lambda}),
  \label{diffeq}
\end{equation}
has a truncated hypergeometric expression \cite{koeswart}
\begin{equation}
  P_n(\eta;\bm{\lambda})\eqdef P_n(\eta;a,b,c;q)
  \eqdef{}_3\phi_2\Bigl(\genfrac{}{}{0pt}{}
  {q^{-n},\,abq^{n+1},\,\eta}{aq,\,cq}\!\Bigm|\!q\,;q\Bigr),\quad
  P_n(1;\bm{\lambda})=1.
  \label{Pn}
\end{equation}
The explicit expressions of the two component eigenpolynomials are
\begin{align}
  \check{P}^{(\pm)}_n(x;\bm{\lambda})
  &\eqdef P_n\bigl(\eta^{(\pm)}(x;\bm{\lambda});\bm{\lambda}\bigr),
  \label{cPnpm}\\
  \check{P}^{(+)}_n(x;\bm{\lambda})
  &={}_3\phi_2\Bigl(\genfrac{}{}{0pt}{}
  {q^{-n},\,abq^{n+1},\,aq^{x+1}}{aq,\,cq}\!\Bigm|\!q\,;q\Bigr),
  \label{Pn+}\\
  \check{P}^{(-)}_n(x;\bm{\lambda})
  &={}_3\phi_2\Bigl(\genfrac{}{}{0pt}{}
  {q^{-n},\,abq^{n+1},\,cq^{x+1}}{aq,\,cq}\!\Bigm|\!q\,;q\Bigr).
  \label{Pn-}
\end{align}
The coefficients of the three term recurrence relation for the big
$q$-Jacobi polynomial $P_n(\eta;\bm{\lambda})$
\begin{equation}
  (1-\eta)P_n(\eta;\bm{\lambda})
  =A_n(\bm{\lambda})P_{n+1}(\eta;\bm{\lambda})-
  \bigl(A_n(\bm{\lambda})+C_n(\bm{\lambda})\bigr)P_n(\eta;\bm{\lambda})
  +C_n(\bm{\lambda})P_{n-1}(\eta;\bm{\lambda}),
  \label{bqj3}
\end{equation}
are
\begin{align} 
  A_n(\bm{\lambda})&\eqdef
  -\frac{(1-aq^{n+1})(1-abq^{n+1})(1-cq^{n+1})}{(1-abq^{2n+1})(1-abq^{2n+2})},
  \label{bqjAn}\\
  C_n(\bm{\lambda})&\eqdef
  ac q^{n+1}\frac{(1-q^n)(1-abc^{-1}q^n)(1-bq^n)}
  {(1-abq^{2n})(1-abq^{2n+1})}.
  \label{bqjCn}
\end{align}
Let us note that the recurrence relation with the initial condition
$P_0(\eta;\bm{\lambda})=1$ implies
$P_n(1;\bm{\lambda})=1$ ($n=1,2,\ldots$) \eqref{Pn}.

{}From \eqref{phi0}, the ground state vectors
$\phi^{(\pm)}_0(x;\bm{\lambda})>0$ are
\begin{align*}
  \phi^{(+)}_0(x;\bm{\lambda})^2&=\phi^{(+)}_0(0;\bm{\lambda})^2\,
  \frac{q^x(aq,abc^{-1}q;q)_x}{(q,ac^{-1}q;q)_x},\n
  \phi^{(-)}_0(x;\bm{\lambda})^2&=\phi^{(-)}_0(0;\bm{\lambda})^2\,
  \frac{q^x(cq,bq;q)_x}{(q,a^{-1}cq;q)_x},
\end{align*}
in which the ratio of $\phi^{(\pm)}_0(0;\bm{\lambda})^2$ is as yet unspecified
but the overall normalisation is immaterial.
They can be rewritten with the newly introduced constant factor
$A^{(\pm)}$ as
\begin{align}
  &\phi^{(\pm)}_0(x;\bm{\lambda})^2
  =\pm A^{(\pm)}\eta^{(\pm)}(x;\bm{\lambda})
  \frac{\bigl(a^{-1}\eta^{(\pm)}(x;\bm{\lambda}),
  c^{-1}\eta^{(\pm)}(x;\bm{\lambda});q\bigr)_{\infty}}
  {\bigl(\eta^{(\pm)}(x;\bm{\lambda}),
  bc^{-1}\eta^{(\pm)}(x;\bm{\lambda});q\bigr)_{\infty}},
  \label{phi02}\\
  &A^{(+)}\eqdef\frac{\phi^{(+)}_0(0;\bm{\lambda})^2(aq,abc^{-1}q;q)_{\infty}}
  {aq(q,ac^{-1}q;q)_{\infty}},\quad
  A^{(-)}\eqdef\frac{\phi^{(-)}_0(0;\bm{\lambda})^2(cq,bq;q)_{\infty}}
  {-cq(q,a^{-1}cq;q)_{\infty}}.\nonumber
\end{align}
The other eigenvectors of the Hamiltonians
$\mathcal{H}^{(\pm)}(\bm{\lambda})$ are
\begin{align}
  &\phi^{(\pm)}_n(x;\bm{\lambda})=\phi^{(\pm)}_0(x;\bm{\lambda})
  \check{P}^{(\pm)}_n(x;\bm{\lambda}),
  \label{phine}\\
  &\mathcal{H}^{(\pm)}(\bm{\lambda})\phi^{(\pm)}_n(x;\bm{\lambda})
  =\mathcal{E}(n;\bm{\lambda})\phi^{(\pm)}_n(x;\bm{\lambda})
  \ \ (n=0,1,\ldots).
  \label{Hephine=}
\end{align}

Now we demonstrate the self-adjointness.
For two smooth functions $\mathcal{P}$ and $\mathcal{Q}$, we take two
vectors $\bm{f}$ and $\bm{g}$ as follows:
\begin{equation*}
  f^{(\pm)}(x)=\phi^{(\pm)}_0(x)
  \mathcal{P}\bigl(\eta^{(\pm)}(x;\bm{\lambda})\bigr),\quad
  g^{(\pm)}(x)=\phi^{(\pm)}_0(x)
  \mathcal{Q}\bigl(\eta^{(\pm)}(x;\bm{\lambda})\bigr).
\end{equation*}
The asymptotic forms of $\phi^{(\pm)}_0(N;\bm{\lambda})^2$ and
$B^{(\pm)}(N;\bm{\lambda})$ at large $N$,
\begin{equation*}
  \phi^{(\pm)}_0(N;\bm{\lambda})^2
  \simeq\pm A^{(\pm)}\eta^{(\pm)}(N;\bm{\lambda})\bigl(1+O(q^N)\bigr),
  \ \ B^{(\pm)}(N;\bm{\lambda})
  \simeq\frac{-acq}{\eta^{(\pm)}(N;\bm{\lambda})^2}\bigl(1+O(q^N)\bigr),
\end{equation*}
and \eqref{largeN:PQ-QP}
show that the criterion for the self-adjointness \eqref{criter1} of
the separate Hamiltonians $\mathcal{H}^{(\pm)}$ does not hold,
\begin{equation*}
  \lim_{N\to\infty}\bigl((f^{(\pm)},\mathcal{H}^{(\pm)}g^{(\pm)})_N
  -(\mathcal{H}^{(\pm)}f^{(\pm)},g^{(\pm)})_N\bigr)
  =\pm A^{(\pm)}(-acq)(1-q)\text{W}[\mathcal{P},\mathcal{Q}](0)\neq 0.
\end{equation*}
For the two component Hamiltonian $\underline{\mathcal{H}}$, that is,
$\mathcal{H}^{(+)}$ combined with $\mathcal{H}^{(-)}$, the self-adjointness is
recovered for the choice
\begin{equation}
  A^{(+)}=A^{(-)}=1,
  \label{A+=A-=1}
\end{equation}
as the criterion \eqref{fHgN-HfgN2} is satisfied
\begin{equation*}
  \lim_{N\to\infty}\bigl({\bin{\bm{f}}{\underline{\mathcal{H}}\,\bm{g}}}_N
  -{\bin{\underline{\mathcal{H}}\,\bm{f}}{\bm{g}}}_N\bigr)
  =\bigl(A^{(+)}-A^{(-)}\bigr)(-acq)(1-q)
  \text{W}[\mathcal{P},\mathcal{Q}](0)=0.
\end{equation*}
With the above choice of $A^{(\pm)}$, the ground state vectors
\eqref{phi02} read
\begin{equation*}
  \phi^{(\pm)}_0(x;\bm{\lambda})^2
  =\pm\eta^{(\pm)}(x;\bm{\lambda})
  \frac{\bigl(a^{-1}\eta^{(\pm)}(x;\bm{\lambda}),
  c^{-1}\eta^{(\pm)}(x;\bm{\lambda});q\bigr)_{\infty}}
  {\bigl(\eta^{(\pm)}(x;\bm{\lambda}),
  bc^{-1}\eta^{(\pm)}(x;\bm{\lambda});q\bigr)_{\infty}},
\end{equation*}
implying that the inner product of $\bm{f}$ and $\bm{g}$ \eqref{binpro}
is expressed by the Jackson integral,
\begin{equation*}
  \bin{\bm{f}}{\bm{g}}= 
  \frac{1}{1-q}\int_{cq}^{aq}\!\!d_qy\,
  \frac{(a^{-1}y,c^{-1}y;q)_{\infty}}{(y,bc^{-1}y;q)_{\infty}}
  \mathcal{P}(y)\mathcal{Q}(y).
\end{equation*}
Thus we arrive at the orthogonality relation of the big $q$-Jacobi polynomials
\begin{align}
  \bin{\bm{\phi}_n}{\bm{\phi}_m}
  &=\bigl(\phi^{(+)}_n,\phi^{(+)}_m\bigr)
  +\bigl(\phi^{(-)}_n,\phi^{(-)}_m\bigr)
  =\frac{\delta_{nm}}{d_n(\bm{\lambda})^2}
  \ \ (n,m=0,1,\ldots)
  \label{bqjortho}\\
  &=\sum_{\epsilon=\pm}\sum_{x=0}^{\infty}
  \phi^{(\epsilon)}_0(x;\bm{\lambda})^2
  P_n\bigl(\eta^{(\epsilon)}(x;\bm{\lambda});\bm{\lambda}\bigr)
  P_m\bigl(\eta^{(\epsilon)}(x;\bm{\lambda});\bm{\lambda}\bigr)\n
  &=aq\sum_{k=0}^{\infty}\frac{(q^{k+1},ac^{-1}q^{k+1};q)_{\infty}}
  {(aq^{k+1},abc^{-1}q^{k+1};q)_{\infty}}
  P_n(aq^{k+1};\bm{\lambda})P_m(aq^{k+1};\bm{\lambda})q^k\n
  &\quad-cq\sum_{k=0}^{\infty}\frac{(q^{k+1},a^{-1}cq^{k+1};q)_{\infty}}
  {(cq^{k+1},bq^{k+1};q)_{\infty}}
  P_n(cq^{k+1};\bm{\lambda})P_m(cq^{k+1};\bm{\lambda})q^k\n
  &=\frac{1}{1-q}\int_{cq}^{aq}\!\!d_qy\,
  \frac{(a^{-1}y,c^{-1}y;q)_{\infty}}{(y,bc^{-1}y;q)_{\infty}}
  P_n(y;\bm{\lambda})P_m(y;\bm{\lambda}),
\end{align}
where the normalisation constant $d_n(\bm{\lambda})^2$ is
\begin{align}
  d_n(\bm{\lambda})^2&\eqdef
  (-acq^2)^{-n}q^{-\frac12n(n-1)}\frac{1-abq^{2n+1}}{1-abq^{n+1}}
  \frac{(aq,abq^2,cq;q)_n}{(q,bq,abc^{-1}q;q)_n}
  \times d_0(\bm{\lambda})^2,\n
  d_0(\bm{\lambda})^2&\eqdef
  \frac{(aq,bq,cq,abc^{-1}q;q)_{\infty}}
  {aq(q,abq^2,ac^{-1}q,a^{-1}c;q)_{\infty}}.
\end{align}
The normalised eigenvectors are
\begin{equation*}
  \hat{\!\bm{\phi}}_n(x;\bm{\lambda})
  =\biggl(\begin{matrix}
  \hat{\phi}^{(+)}_n(x;\bm{\lambda})\\
  \hat{\phi}^{(-)}_n(x;\bm{\lambda})
  \end{matrix}\biggr),\quad
  \hat{\phi}^{(\pm)}_n(x;\bm{\lambda})
  =d_n(\bm{\lambda}){\phi}^{(\pm)}_0(x;\bm{\lambda})
  \check{P}^{(\pm)}_n(x;\bm{\lambda}),
\end{equation*}
and they satisfy the orthogonality relation 
\begin{equation*}
  \bin{\,\hat{\!\bm{\phi}}_n}{\,\hat{\!\bm{\phi}}_m}
  =\bigl(\hat{\phi}^{(+)}_n,\hat{\phi}^{(+)}_n\bigr)
  +\bigl(\hat{\phi}^{(-)}_n,\hat{\phi}^{(-)}_n\bigr)
  =\delta_{nm}.
\end{equation*}

The forward/backward shift relations for the polynomial
$P_n(\eta;\bm{\lambda})$ are
\begin{align*}
  &\frac{(1-aq)(1-cq)}{q\eta}
  \bigl(P_n(\eta;\bm{\lambda})-P_n(q\eta;\bm{\lambda})\bigr)
  =\mathcal{E}(n;\bm{\lambda})
  P_{n-1}(q\eta;\bm{\lambda}+\bm{\delta}),\\
  &
  \bigl(B^{\text{J}}(\eta;\bm{\lambda})P_{n-1}(q\eta;\bm{\lambda}+\bm{\delta})
  -q^{-1}D^{\text{J}}(\eta;\bm{\lambda})P_{n-1}(\eta;\bm{\lambda}+\bm{\delta})
  \bigr)\frac{q\eta}{(1-aq)(1-cq)}
  =P_n(\eta;\bm{\lambda}).\!\!
\end{align*}
These relations can be rewritten as the forward/backward shift relations
for the polynomials $\check{P}^{(\pm)}_n(x;\bm{\lambda})$ on
the $x$-lattice, \eqref{FB}--\eqref{FP=}
\begin{align}
  &\mathcal{F}^{(\pm)}(\bm{\lambda})=
  \varphi^{(\pm)}(x;\bm{\lambda})^{-1}(1-e^{\partial}),
  \ \ \mathcal{B}^{(\pm)}(\bm{\lambda})=
  \bigl(B^{(\pm)}(x;\bm{\lambda})
  -D^{(\pm)}(x;\bm{\lambda})e^{-\partial}\bigr)
  \varphi^{(\pm)}(x;\bm{\lambda}),\!\!
  \label{FBe}\\
  &\mathcal{F}^{(\pm)}(\bm{\lambda})
  \check{P}^{(\pm)}_n(x;\bm{\lambda})
  =\mathcal{E}(n;\bm{\lambda})
  \check{P}^{(\pm)}_{n-1}(x;\bm{\lambda}+\bm{\delta}),
  \ \ \mathcal{B}^{(\pm)}(\bm{\lambda})
  \check{P}^{(\pm)}_{n-1}(x;\bm{\lambda}+\bm{\delta})
  =\check{P}^{(\pm)}_n(x;\bm{\lambda}),
  \label{FBPne}
\end{align}
where the auxiliary functions $\varphi^{(\pm)}(x;\bm{\lambda})$ are
\begin{equation}
  \varphi^{(\pm)}(x;\bm{\lambda})\eqdef
  \frac{q\eta^{(\pm)}(x;\bm{\lambda})}{(1-aq)(1-cq)}.
  \label{bqJ:varphipm}
\end{equation}

It is interesting to note that the following parameter substitution
(involution)
\begin{equation}
  (a,b,c)\to(c,abc^{-1},a),
  \label{bqjsym}
\end{equation}
gives rise to the interchange of the $(+)$ and $(-)$ systems
\begin{align}
  &\begin{array}{ccc}
  B^{(+)}(x;\bm{\lambda})\ \leftrightarrow\ B^{(-)}(x;\bm{\lambda}),&
  &\eta^{(+)}(x;\bm{\lambda})\ \leftrightarrow\ \eta^{(-)}(x;\bm{\lambda}),
  \\[2pt]
  D^{(+)}(x;\bm{\lambda})\ \leftrightarrow\ D^{(-)}(x;\bm{\lambda}),&
  &\check{P}^{(+)}_n(x;\bm{\lambda})\ \leftrightarrow
  \ \check{P}^{(-)}_n(x;\bm{\lambda}),
  \end{array}\n
  &B^{\text{J}}(\eta;\bm{\lambda}),\ D^{\text{J}}(\eta;\bm{\lambda}),
  \ \mathcal{E}(n;\bm{\lambda}),\ P_n(\eta;\bm{\lambda})
  :\ \text{invariant},
  \label{bqJ:sym}
\end{align}
but the parameter range \eqref{bqjpararange} is not preserved.

In \cite{os14} we have presented a unified prescription to obtain
(quasi-)exactly solvable difference Schr\"odinger equations
in which the potential functions are constructed in terms of the
sinusoidal coordinates.
Here we show how the b$q$J system fits in that scheme, although the
solvability is obvious from the start in the two component Hamiltonian
formalism.
By taking $v_{k,l}$ as
\begin{equation}
  \begin{array}{ll}
  v_{0,0}=-(1-q)(1-q^2)ac,\\[2pt]
  v_{1,0}=(1-q)\bigl(a(b+c)-q(a+c)\bigr),&
  v_{0,1}=(1-q)\bigl(a+c-qa(b+c)\bigr),\\[2pt]
  v_{2,0}=(1-q)(1-ab)+v_{0,2},&
  v_{1,1}=(1-q^{-1})(1-abq^2)-(q+q^{-1})v_{0,2},
  \end{array}
  \label{vkl}
\end{equation}
($v_{0,2}$ is arbitrary), the potential functions are expressed as follows:
\begin{align}
  B^{(\pm)}(x;\bm{\lambda})&=\frac{
  \sum_{\genfrac{}{}{0pt}{}{k,l\geq 0}{k+l\leq 2}}
  v_{k,l}\eta^{(\pm)}(x;\bm{\lambda})^k\eta^{(\pm)}(x+1;\bm{\lambda})^l}
  {\bigl(\eta^{(\pm)}(x+1;\bm{\lambda})-\eta^{(\pm)}(x;\bm{\lambda})\bigr)
  \bigl(\eta^{(\pm)}(x+1;\bm{\lambda})-\eta^{(\pm)}(x-1;\bm{\lambda})\bigr)},
  \label{os14B}\\
  D^{(\pm)}(x;\bm{\lambda})&=\frac{
  \sum_{\genfrac{}{}{0pt}{}{k,l\geq 0}{k+l\leq 2}}
  v_{k,l}\eta^{(\pm)}(x;\bm{\lambda})^k\eta^{(\pm)}(x-1;\bm{\lambda})^l}
  {\bigl(\eta^{(\pm)}(x-1;\bm{\lambda})-\eta^{(\pm)}(x;\bm{\lambda})\bigr)
  \bigl(\eta^{(\pm)}(x-1;\bm{\lambda})-\eta^{(\pm)}(x+1;\bm{\lambda})\bigr)}.
  \label{os14D}
\end{align}
All the above parameters $v_{k,l}$ \eqref{vkl} are invariant under the
parameter substitution \eqref{bqjsym}.
In the logic of \cite{os14}, the upper triangularity of the similarity
transformed Hamiltonians $\widetilde{\mathcal{H}}^{(\pm)}$  with
respect to the bases
$\{\eta^{(\pm)}(x;\bm{\lambda})^k\}_{k=0,1,\ldots}$ is shown
through the expansion:
\begin{equation*}
  \frac{\eta^{(\pm)}(x+1;\bm{\lambda})^{n+1}
  -\eta^{(\pm)}(x-1;\bm{\lambda})^{n+1}}
  {\eta^{(\pm)}(x+1;\bm{\lambda})-\eta^{(\pm)}(x-1;\bm{\lambda})}
  =\sum_{k=0}^ng_n^{(k)\,(\pm)}\eta^{(\pm)}(x;\bm{\lambda})^{n-k},
\end{equation*}
which is the consequence of the properties \eqref{etaprop1}--\eqref{etaprop2}.
In fact, the above formula is rather trivial in the present case.

\subsubsection{dual big $q$-Jacobi polynomials}
\label{sec:dbqJ}

As explained in \S\,\ref{sec:rdQM}, dual polynomials are usually defined
from the original polynomials by interchanging the roles of $x$ and $n$
\eqref{duality1}--\eqref{duality2}.
Fixing $x\in\mathbb{Z}_{\ge0}$ in $\check{P}^{(\pm)}_n(x;\bm{\lambda})$
\eqref{Pn+}--\eqref{Pn-}, however, does not provide degree $x$ polynomials
in $\mathcal{E}(n;\bm{\lambda})$ \eqref{bqjeta}.
The reason is that the original eigenpolynomials
$\check{P}^{(\pm)}_n(x;\bm{\lambda})$ \eqref{cPnpm} fail to satisfy
the universal normalisation condition
$\check{P}^{(\pm)}_n(0;\bm{\lambda})\neq1$ \eqref{Pzero}.

The dual big $q$-Jacobi polynomials (db$q$J) are defined from the original
polynomials $\check{P}^{(\pm)}_n(x;\bm{\lambda})$ \eqref{cPnpm} by
removing the $n$ dependent factor $\alpha^{(\pm)}_n(\bm{\lambda})$,
which can be identified by the transformation formula of the truncated
${}_3\phi_2$ function listed as eq.(III.12) in the Appendix III of
\cite{gasper}:
\begin{equation*}
  {}_3\phi_2\Bigl(\genfrac{}{}{0pt}{}
  {q^{-n},\,b,\,c}{d,\,e}\!\Bigm|\!q\,;q\Bigr)
  =c^n\frac{(c^{-1}e;q)_n}{(e;q)_n}{}_3\phi_2\Bigl(\genfrac{}{}{0pt}{}
  {q^{-n},\,c,\,b^{-1}d}{d,\,ce^{-1}q^{1-n}}\!\Bigm|\!q\,;\frac{bq}{e}\Bigr).
\end{equation*}
Let us introduce $\check{Q}^{(\pm)}_x(n;\bm{\lambda})$ by
\begin{align}
  &\check{P}^{(\pm)}_n(x;\bm{\lambda})
  =\alpha^{(\pm)}_n(\bm{\lambda})
  \check{Q}^{(\pm)}_x(n;\bm{\lambda}),\quad 
  \check{Q}^{(\pm)}_x(0;\bm{\lambda})=1,\quad
  (\pm1)^n\alpha^{(\pm)}_n(\bm{\lambda})>0,
  \label{bqjdualpoly}\\
  &\check{Q}^{(+)}_x(n;\bm{\lambda})
  \eqdef Q^{(+)}_x\bigl(\mathcal{E}(n;\bm{\lambda});\bm{\lambda}\bigr)=
  {}_3\phi_2\Bigl(\genfrac{}{}{0pt}{}
  {q^{-n},\,abq^{n+1},\,q^{-x}}{aq,\,abc^{-1}q}\!\Bigm|\!q\,;
  ac^{-1}q^{x+1}\Bigr),
  \label{Qx+}\\
  &\check{Q}^{(-)}_x(n;\bm{\lambda})
  \eqdef Q^{(-)}_x\bigl(\mathcal{E}(n;\bm{\lambda});\bm{\lambda}\bigr)=
  {}_3\phi_2\Bigl(\genfrac{}{}{0pt}{}
  {q^{-n},\,abq^{n+1},\,q^{-x}}{bq,\,cq}\!\Bigm|\!q\,;
  a^{-1}cq^{x+1}\Bigr),
  \label{Qx-}\\
  &\alpha^{(+)}_n(\bm{\lambda})\eqdef
  (-c)^nq^{\frac12n(n+1)}\frac{(abc^{-1}q;q)_n}{(cq;q)_n},\quad
  \alpha^{(-)}_n(\bm{\lambda})\eqdef
  (-a)^nq^{\frac12n(n+1)}\frac{(bq;q)_n}{(aq;q)_n}.
  \label{alphan-}
\end{align}
For fixed $x\in\mathbb{Z}_{\ge0}$,
$Q^{(\pm)}_x(\mathcal{E}(n;\bm{\lambda});\bm{\lambda})$ are
degree $x$ polynomials in $\mathcal{E}(n;\bm{\lambda})$ \eqref{bqjeta}.
The difference equation for the b$q$J \eqref{diffeq} can be rewritten as the
three term recurrence relation for
$\check{Q}^{(\pm)}_x(n;\bm{\lambda})$
\begin{align}
  \mathcal{E}(n;\bm{\lambda})\check{Q}^{(\pm)}_x(n;\bm{\lambda})
  &=-B^{(\pm)}(x;\bm{\lambda})\check{Q}^{(\pm)}_{x+1}(n;\bm{\lambda})
  -D^{(\pm)}(x;\bm{\lambda})\check{Q}^{(\pm)}_{x-1}(n;\bm{\lambda})\n
  &\quad
  +\bigl(B^{(\pm)}(x;\bm{\lambda})+D^{(\pm)}(x;\bm{\lambda})\bigr)
  \check{Q}^{(\pm)}_x(n;\bm{\lambda}).
  \label{3termQxe}
\end{align}
In order to write down the difference equations for the dual polynomials
$\check{Q}^{(\pm)}_x(n;\bm{\lambda})$,
let us introduce $A^{(\pm)}_n(\bm{\lambda})$,
$C^{(\pm)}_n(\bm{\lambda})$ by rescaling the coefficients
$A_n(\bm{\lambda})$ and $C_n(\bm{\lambda})$ \eqref{bqjAn}--\eqref{bqjCn}
of the three term recurrence relation of the big $q$-Jacobi polynomial
\eqref{bqj3} by $\alpha^{(\pm)}_n(\bm{\lambda})$:
\begin{equation}
  A^{(\pm)}_n(\bm{\lambda})\eqdef
  \pm \frac{\alpha^{(\pm)}_{n+1}(\bm{\lambda})}{\alpha^{(\pm)}_n(\bm{\lambda})}
  A_n(\bm{\lambda}),\quad
  C^{(\pm)}_n(\bm{\lambda})\eqdef
  \pm \frac{\alpha^{(\pm)}_{n-1}(\bm{\lambda})}{\alpha^{(\pm)}_n(\bm{\lambda})}
  C_n(\bm{\lambda}).
  \label{AneCne}
\end{equation}
To be more explicit, we have
\begin{align*}
  &A^{(+)}_n(\bm{\lambda})
  =cq^{n+1}\frac{(1-aq^{n+1})(1-abq^{n+1})(1-abc^{-1}q^{n+1})}
  {(1-abq^{2n+1})(1-abq^{2n+2})},\\
  &A^{(-)}_n(\bm{\lambda})
  =-aq^{n+1}\frac{(1-bq^{n+1})(1-abq^{n+1})(1-cq^{n+1})}
  {(1-abq^{2n+1})(1-abq^{2n+2})},\\
  &C^{(+)}_n(\bm{\lambda})
  =-aq\frac{(1-q^n)(1-bq^n)(1-cq^n)}{(1-abq^{2n})(1-abq^{2n+1})},\\
  &C^{(-)}_n(\bm{\lambda})
  =cq\frac{(1-q^n)(1-aq^n)(1-abc^{-1}q^n)}{(1-abq^{2n})(1-abq^{2n+1})}.
\end{align*}
The sign change factor $\pm$ in \eqref{AneCne} is to compensate the
sign change of $\alpha^{(-)}_n(\bm{\lambda})$ \eqref{alphan-}.
The potential functions of the dual Hamiltonians are {\em bounded} and
positive semi-definite 
\begin{equation*}
  -A^{(\pm)}_n(\bm{\lambda})>0\ \ (n\geq 0),\quad
  -C^{(\pm)}_n(\bm{\lambda})>0\ \ (n\geq 1),
  \ \ C^{(\pm)}_0(\bm{\lambda})=0,
\end{equation*}
and they satisfy the relations
\begin{align}
  A_n(\bm{\lambda})C_{n+1}(\bm{\lambda}) 
  &=A^{(\pm)}_n(\bm{\lambda})C^{(\pm)}_{n+1}(\bm{\lambda}),
  \label{bqJ:AnCn+1}\\
  A_n(\bm{\lambda})+C_n(\bm{\lambda})
  &=\pm\bigl(
  A^{(\pm)}_n(\bm{\lambda})+C^{(\pm)}_n(\bm{\lambda})\bigr)
  -\bigl(1-\eta^{(\pm)}(0;\bm{\lambda})\bigr).
  \label{bqJ:An+Cn}
\end{align}
Under the parameter substitution \eqref{bqjsym}, the $(+)$ and $(-)$
quantities are interchanged
\begin{equation*}
  \check{Q}^{(+)}_x(n;\bm{\lambda})\leftrightarrow
  \check{Q}^{(-)}_x(n;\bm{\lambda}),\quad
  A^{(+)}_n(\bm{\lambda})\leftrightarrow-A^{(-)}_n(\bm{\lambda}),\quad
  C^{(+)}_n(\bm{\lambda})\leftrightarrow-C^{(-)}_n(\bm{\lambda}).
\end{equation*}
Thanks to \eqref{bqJ:An+Cn} and
$\eta^{(\pm)}(x;\bm{\lambda})=\eta^{(\pm)}(0;\bm{\lambda})q^x$,
the three term recurrence relation of the b$q$J \eqref{bqj3}
is rewritten as the difference equations for the dual b$q$J
\eqref{Qx+}--\eqref{Qx-}:
\begin{align}
  &\quad-A^{(\pm)}_n(\bm{\lambda})
  \bigl(\check{Q}^{(\pm)}_x(n;\bm{\lambda})
  -\check{Q}^{(\pm)}_x(n+1;\bm{\lambda})\bigr)
  -C^{(\pm)}_n(\bm{\lambda})
  \bigl(\check{Q}^{(\pm)}_x(n;\bm{\lambda})
  -\check{Q}^{(\pm)}_x(n-1;\bm{\lambda})\bigr)\n
  &=\pm\eta^{(\pm)}(0;\bm{\lambda})\eta(x)
  \check{Q}^{(\pm)}_x(n;\bm{\lambda}),\quad
  \eta(x)\eqdef 1-q^x,
  \label{diffeqQxe}
\end{align}
with positive semi-definite and {\em bounded} eigenvalues
\begin{equation}
  \mathcal{E}^{\text{d}\,(\pm)}(x;\bm{\lambda})
  \eqdef\pm\eta^{(\pm)}(0;\bm{\lambda})\eta(x)=\left\{
  \begin{array}{rl}
  aq(1-q^x)&:(+)\\
  -cq(1-q^x)&:(-)
  \end{array}\right..
  \label{bqJ:Edpm}
\end{equation}
In terms of the similarity transformed dual Hamiltonians
\begin{equation}
  \widetilde{\mathcal{H}}^{\text{d}\,(\pm)}(\bm{\lambda})
  =-A^{(\pm)}_n(\bm{\lambda})(1- e^{\partial_n})
  -C^{(\pm)}_n(\bm{\lambda})(1- e^{-\partial_n}),
 \label{dbqJtrham}
\end{equation}
the above equations \eqref{diffeqQxe} reads succinctly
\begin{equation*}
  \widetilde{\mathcal{H}}^{\text{d}\,(\pm)}(\bm{\lambda})
  \check{Q}^{(\pm)}_x(n;\bm{\lambda})
  =\pm\eta^{(\pm)}(0;\bm{\lambda})\eta(x)\check{Q}^{(\pm)}_x(n;\bm{\lambda}).
\end{equation*}
It is straightforward to verify that the similarity transformed dual
Hamiltonians $\widetilde{\mathcal{H}}^{\text{d}\,(\pm)}(\bm{\lambda})$
\eqref{dbqJtrham} are triangular in the basis
$\{1,\mathcal{E}(n;\bm{\lambda}),\ldots, \mathcal{E}(n;\bm{\lambda})^k\}$
with the above positive semi-definite eigenvalues.
Their eigenpolynomials
$\check{Q}^{(\pm)}_x(n;\bm{\lambda})\sim\mathcal{E}(n;\bm{\lambda})^x
\simeq q^{-nx}$ grow rapidly as the coordinate $n\to\infty$.

Let us define two dual Hamiltonians
$\mathcal{H}^{\text{d}\,(\pm)}(\bm{\lambda})$ as follows:
\begin{align}
  \mathcal{H}^{\text{d}\,(\pm)}(\bm{\lambda})&\eqdef
  \mp\sqrt{A^{(\pm)}_n(\bm{\lambda})C^{(\pm)}_{n+1}(\bm{\lambda})}
  \,e^{\partial_n}
  \mp\sqrt{A^{(\pm)}_{n-1}(\bm{\lambda})C^{(\pm)}_n(\bm{\lambda})}
  \,e^{-\partial_n}
  -\bigl(A^{(\pm)}_n(\bm{\lambda})+C^{(\pm)}_n(\bm{\lambda})\bigr)\n
  &=\mathcal{A}^{\text{d}\,(\pm)}(\bm{\lambda})^\dagger
  \mathcal{A}^{\text{d}\,(\pm)}(\bm{\lambda}),
  \label{dbqJham}\\
  \mathcal{A}^{\text{d}\,(\pm)}(\bm{\lambda})&\eqdef
  \sqrt{-A^{(\pm)}_n(\bm{\lambda})}
  \mp e^{\partial_n}\sqrt{-C^{(\pm)}_n(\bm{\lambda})}.\nonumber
\end{align}
Note that $\mathcal{H}^{\text{d}\,(-)}(\bm{\lambda})$ and
$\mathcal{A}^{\text{d}\,(-)}(\bm{\lambda})$ are related to the standard
forms \eqref{Hdef2}--\eqref{A,Ad} by the similarity transformation in
terms of the diagonal matrix $(-1)^n=\text{diag}(1,-1,1,-1,\ldots)$,
\begin{align*}
  (-1)^n\circ\mathcal{H}^{\text{d}\,(-)}\circ(-1)^n
  &=-\sqrt{A^{(-)}_nC^{(-)}_{n+1}}\,e^{\partial_n}
  -\sqrt{A^{(-)}_{n-1}C^{(-)}_n}\,e^{-\partial_n}
  -\bigl(A^{(-)}_n+C^{(-)}_n\bigr),\n
  (-1)^n\circ\mathcal{A}^{\text{d}\,(-)}\circ(-1)^n
  &=\sqrt{-A^{(-)}_n}-e^{\partial_n}\sqrt{-C^{(\pm)}_n},
\end{align*}
where we have suppressed $\bm{\lambda}$.
The ground state eigenvectors $\phi^{\text{d}\,(\pm)}_0(n;\bm{\lambda})$
characterised by
$\mathcal{A}^{\text{d}\,(\pm)}(\bm{\lambda})
\phi^{\text{d}\,(\pm)}_0(n;\bm{\lambda})=0$ and
$\phi^{\text{d}\,(\pm)}_0(0;\bm{\lambda})=1$ are
\begin{equation}
  \phi^{\text{d}\,(\pm)}_0(n;\bm{\lambda})
  =(\pm1)^n\prod_{m=0}^{n-1}\sqrt{\frac{A^{(\pm)}_m(\bm{\lambda})}
  {C^{(\pm)}_{m+1}(\bm{\lambda})}}
  =\alpha^{(\pm)}_n(\bm{\lambda})\frac{d_n(\bm{\lambda})}{d_0(\bm{\lambda})},
  \quad(\pm1)^n\phi^{\text{d}\,(\pm)}_0(n;\bm{\lambda})>0.
  \label{dbqJpsi0}
\end{equation}
The dual polynomials $\check{Q}^{(\pm)}_x(n;\bm{\lambda})$ give eigenvectors
\begin{align}
  &\phi^{\text{d}\,(\pm)}_x(n;\bm{\lambda})\eqdef
  \phi^{\text{d}\,(\pm)}_0(n;\bm{\lambda})
  \check{Q}^{(\pm)}_x(n;\bm{\lambda}),
  \label{bqJ:phidpm}\\
  &\mathcal{H}^{\text{d}\,(\pm)}(\bm{\lambda})
  \phi^{\text{d}\,(\pm)}_x(n;\bm{\lambda})
  =\mathcal{E}^{\text{d}\,(\pm)}(x;\bm{\lambda})
  \phi^{\text{d}\,(\pm)}_x(n;\bm{\lambda}).
  \label{bqJ:Ed}
\end{align}
By using \eqref{bqJ:AnCn+1}--\eqref{bqJ:An+Cn}, we can show that the two
dual Hamiltonians $\mathcal{H}^{\text{d}\,(\pm)}(\bm{\lambda})$ are
related as follows:
\begin{equation}
  \mathcal{H}^{\text{d}\,(+)}(\bm{\lambda})-\eta^{(+)}(0;\bm{\lambda})
  =-\mathcal{H}^{\text{d}\,(-)}(\bm{\lambda})-\eta^{(-)}(0;\bm{\lambda}).
  \label{HmHprel}
\end{equation}
Therefore $\mathcal{H}^{\text{d}\,(\pm)}(\bm{\lambda})$ have another set
of eigenvectors $\phi^{\text{d}\,(\mp)}_x(n;\bm{\lambda})$,
\begin{equation}
  \mathcal{H}^{\text{d}\,(\pm)}(\bm{\lambda})
  \phi^{\text{d}\,(\mp)}_x(n;\bm{\lambda})
  =\mathcal{E}^{\prime\,\text{d}\,(\pm)}(x;\bm{\lambda})
  \phi^{\text{d}\,(\mp)}_x(n;\bm{\lambda}),\quad
  \mathcal{E}^{\prime\,\text{d}\,(\pm)}(x;\bm{\lambda})\eqdef\left\{
  \begin{array}{ll}
  q(a-cq^x)&:(+)\\
  q(-c+aq^x)&:(-)
  \end{array}\right..
  \label{bqJ:E'd}
\end{equation}
The `dual Hamiltonian' $\mathcal{H}^{\text{d}}(\bm{\lambda})$ \eqref{dbqJHam}
may be constructed from the original coefficients $A_n(\bm{\lambda})$ and
$C_n(\bm{\lambda})$ \eqref{bqjAn}--\eqref{bqjCn} of the three term
recurrence relation of the b$q$J:
\begin{align}
  \mathcal{H}^{\text{d}}(\bm{\lambda})&\eqdef
  -\sqrt{A_n(\bm{\lambda})C_{n+1}(\bm{\lambda})}\,e^{\partial_n}-
  \sqrt{A_{n-1}(\bm{\lambda})C_{n}(\bm{\lambda})}\,e^{-\partial_n}
  -\bigl(A_n(\bm{\lambda})+C_{n}(\bm{\lambda})\bigr)
  \label{dbqJHam}\\
  &=\pm\mathcal{H}^{\text{d}\,(\pm)}(\bm{\lambda})
  +1-\eta^{(\pm)}(0;\bm{\lambda}).\nonumber
\end{align}
However, the zero mode of this `dual Hamiltonian'
$\mathcal{H}^{\text{d}}(\bm{\lambda})$ is not square summable
\begin{equation*}
  \phi^{\text{d}}_0(n;\bm{\lambda})^2\propto
  \prod_{m=0}^{n-1}\frac{A_m(\bm{\lambda})}{C_{m+1}(\bm{\lambda})}
  =\frac{d_n(\bm{\lambda})^2}{d_0(\bm{\lambda})^2}
  \simeq(-ac)^{-n}q^{-\frac12n(n+3)}\bigl(1+O(q^n)\bigr)\times\text{const}
  \ \ (n\to\infty),
\end{equation*}
and it fails to deliver the dual big $q$-Jacobi polynomials.
This is why we have to introduce two dual Hamiltonians
$\mathcal{H}^{\text{d}\,(\pm)}(\bm{\lambda})$ \eqref{dbqJham}.

We arrive at the dual orthogonality relations ($x,y=0,1,\ldots$)
\begin{align}
  (\phi^{\text{d}\,(\pm)}_x,\phi^{\text{d}\,(\pm)}_y)
  &=\sum_{n=0}^{\infty}\phi^{\text{d}\,(\pm)}_0(n;\bm{\lambda})^2
  \check{Q}^{(\pm)}_x(n;\bm{\lambda})
  \check{Q}^{(\pm)}_y(n;\bm{\lambda})
  =\frac{\delta_{xy}}{\phi^{(\pm)}_0(x;\bm{\lambda})^2d_0(\bm{\lambda})^2},
  \label{Q-ortho1}\\
  (\phi^{\text{d}\,(+)}_x,\phi^{\text{d}\,(-)}_y)
  &=\sum_{n=0}^{\infty}\phi^{\text{d}\,(+)}_0(n;\bm{\lambda})
  \phi^{\text{d}\,(-)}_0(n;\bm{\lambda})
  \check{Q}^{(+)}_x(n;\bm{\lambda})
  \check{Q}^{(-)}_y(n;\bm{\lambda})=0,
  \label{Q-ortho3}\\
  (\hat{\phi}^{\text{d}\,(\epsilon)}_x,\hat{\phi}^{\text{d}\,(\epsilon')}_y)
  &=\delta_{\epsilon\,\epsilon'}\delta_{xy},\quad
  \hat{\phi}^{\text{d}\,(\pm)}_x(n;\bm{\lambda})\eqdef
  \phi^{(\pm)}_0(x;\bm{\lambda})d_0(\bm{\lambda})
  \phi^{\text{d}\,(\pm)}_x(n;\bm{\lambda}).
  \label{Q-ortho4}
\end{align}
These are the consequences of the self-adjointness of the Hamiltonian
$\mathcal{H}^{\text{d}\,(+)}(\bm{\lambda})$ \eqref{dbqJham}, for which
$\{\phi^{\text{d}\,(+)}_x(n;\bm{\lambda})\}$ and
$\{\phi^{\text{d}\,(-)}_x(n;\bm{\lambda})\}$ are eigenvectors.
The spectrum of the latter
$\{\mathcal{E}^{\prime\,\text{d}\,(+)}(x;\bm{\lambda})=q(a-cq^x)\}$
\eqref{bqJ:E'd}, which is monotonously decreasing with $x$,
lies above that of the former,
$\{\mathcal{E}^{\text{d}\,(+)}(x)=aq(1-q^x)\}$ \eqref{bqJ:Ed},
which is monotonously increasing with $x$.
They share the same accumulation point $aq$.
The ground state vector $\phi^{\text{d}\,(+)}_0(n;\bm{\lambda})$ is positive
and the excited state vector $\phi^{\text{d}\,(+)}_x(n;\bm{\lambda})$ has
$x$ `zeros' due to $x$ zeros of the polynomial
$Q^{(+)}_x(\mathcal{E};\bm{\lambda})$. Here the number of `zeros' means the
number of sign changing in $n\in[0,\infty)$.
All the eigenvectors $\{\phi^{\text{d}\,(-)}_x(n;\bm{\lambda})\}$ have
infinitely many `zeros' due to the infinitely many eigenlevels of
$\{\phi^{\text{d}\,(+)}_x(n;\bm{\lambda})\}$ lying below.
These infinite `zeros' are partly taken care of by the alternating sign
factor $(-1)^n$ or $\alpha^{(-)}_n(\bm{\lambda})$ in
$\phi^{\text{d}\,(-)}_0(n;\bm{\lambda})$ \eqref{dbqJpsi0}.
Since zeros of $Q^{(-)}_x(\mathcal{E};\bm{\lambda})$ cancels the sign
change of $(-1)^n$, it is easy to ``understand'' that the number of
`zeros' in $\phi^{\text{d}\,(-)}_x(n;\bm{\lambda})$
would increase as the degree $x$ decreases.
This is consistent with the oscillation theorem.
The orthogonality relation \eqref{Q-ortho3} simply means
$\sum\limits_{n=0}^{\infty}\phi^{\text{d}\,(+)}_0(n;\bm{\lambda})
\phi^{\text{d}\,(-)}_0(n;\bm{\lambda})\mathcal{E}(n;\bm{\lambda})^k=0$,
namely
\begin{equation}
  \sum_{n=0}^{\infty}(-1)^n\frac{q^{\frac12n(n-1)}}{(q;q)_n}
  \frac{1-abq^{2n+1}}{1-abq^{n+1}}(abq^2;q)_n
  \,\mathcal{E}(n;\bm{\lambda})^k
  =0\ \ (k\in\mathbb{Z}_{\geq 0}),
  \label{bqJ:orthopm}
\end{equation}
as mentioned in \cite{atakishi4}.

The matrix elements of $\mathcal{H}^{\text{d}\,(-)}(\bm{\lambda})$
\eqref{dbqJham} are all non-negative.
The eigenvectors and eigenvalues are
$\{\phi^{\text{d}\,(-)}_x(n;\bm{\lambda})\}$ with
$\{\mathcal{E}^{\text{d}\,(-)}(x)=-cq(1-q^x)\}$ \eqref{bqJ:Ed} and
$\{\phi^{\text{d}\,(+)}_x(n;\bm{\lambda})\}$ with
$\{\mathcal{E}^{\prime\,\text{d}\,(-)}(x;\bm{\lambda})=q(-c+aq^x)\}$
\eqref{bqJ:E'd}.
The maximum eigenvalue is
$\mathcal{E}^{\prime\,\text{d}\,(-)}(0;\bm{\lambda})=q(a-c)$ and
the corresponding eigenvector
$\phi^{\text{d}\,(+)}_0(n;\bm{\lambda})$ is positive everywhere.
This is consistent with the Perron-Frobenius theorem.
The ground state vector
$\phi^{\text{d}\,(-)}_0(n;\bm{\lambda})$ has infinitely many `zeros'.

The alternating sign factor has no effect for the orthogonality among
$\{\check{Q}^{(-)}_x(n;\bm{\lambda})\}$ \eqref{Q-ortho1}, which is just
an ordinary orthogonality relation among polynomials.
It should be stressed that the weight function for the orthogonality of
the two different kinds of polynomials,
$\{\check{Q}^{(+)}_x(n;\bm{\lambda})\}$ and
$\{\check{Q}^{(-)}_x(n;\bm{\lambda})\}$ \eqref{Q-ortho3} is
{\em not positive definite} due to the alternating sign factor in
$\phi^{\text{d}\,(-)}_0(n;\bm{\lambda})$ \eqref{dbqJpsi0}.
If the weight function were positive definite, the orthogonality of
$\check{Q}^{(+)}_0(n;\bm{\lambda})=1=\check{Q}^{(-)}_0(n;\bm{\lambda})$
could not be attained.
As is well known, the overall scale of the Hamiltonian is immaterial.
In order to get the common eigenvalues $\{1-q^x\}$ for both Hamiltonians
$\mathcal{H}^{\text{d}\,(\pm)}(\bm{\lambda})$ \eqref{dbqJham},
one should divide $A^{(+)}_n$, $C^{(+)}_n$ by $aq$ in
$(+)$ sector and $A^{(-)}_n$, $C^{(-)}_n$ by $-cq$ in $(-)$ sector.

As discussed in \S\,\ref{sec:rdQM}, dual polynomials are defined through the
completeness relation \eqref{dualortho} of the original orthogonal polynomials.
In the case of big $q$-Jacobi polynomials, one multiplies normalised
eigenfunction $\hat{\phi}^{(\epsilon)}_n(y)\eqdef d_n{\phi}^{(\epsilon)}_n(y)$
to the original normalised orthogonality relation \eqref{bqjortho}
$(\hat{\phi}^{(+)}_n,\hat{\phi}^{(+)}_m)
+(\hat{\phi}^{(-)}_n,\hat{\phi}^{(-)}_m)=\delta_{nm}$
and sum over $n$. Since the original orthogonality relation is uniformly
convergent, one obtains the completeness relation or the dual orthogonality
relation
\begin{equation}
  \sum_{n=0}^{\infty}\hat{\phi}^{(+)}_n(x)\hat{\phi}^{(+)}_n(y)
  =\delta_{xy},\quad
  \sum_{n=0}^{\infty}\hat{\phi}^{(-)}_n(x)\hat{\phi}^{(-)}_n(y)
  =\delta_{xy},\quad
  \sum_{n=0}^{\infty}\hat{\phi}^{(+)}_n(x)\hat{\phi}^{(-)}_n(y)=0.
  \label{naivdual}
\end{equation}
Rewriting the b$q$J $\check{P}^{(\pm)}_n(x;\bm{\lambda})$ in the above
orthogonality relations \eqref{naivdual} in terms of the dual polynomials
$\check{Q}^{(\pm)}_x(n;\bm{\lambda})$ \eqref{bqjdualpoly} provides the
orthogonality relations \eqref{Q-ortho1}--\eqref{Q-ortho3}.

The dual big $q$-Jacobi polynomials have been introduced and discussed in
\cite{atakishi2,atakishi4}.

\subsection{Limits of big $q$-Jacobi}
\label{sec:lim_bqJ}

By restricting the parameters or by taking certain limits of the b$q$J
system presented in the previous subsection, we obtain three systems
described by the big $q$-Laguerre (b$q$L), Al-Salam-Carlitz $\I$ (ASC\,$\I$)
and discrete $q$-Hermite $\I$ (d$q$H\,$\I$) polynomials. The orthogonality of
these polynomials requires the two component Hamiltonian formalism for the
self-adjointness and their orthogonality measures are of Jackson integral
type. The duals of these polynomials have also been discussed in
\cite{atakishi4} in some detail. We will provide main results skipping
most of derivation.
The coordinate $x$ takes non-negative integer values,
$x\in\mathbb{Z}_{\geq 0}$. 

\subsubsection{big $q$-Laguerre}
\label{sec:bqL}

The big $q$-Laguerre polynomial is obtained from the big $q$-Jacobi (b$q$J)
polynomial by setting
$(a^{\text{b$q$J}},b^{\text{b$q$J}},c^{\text{b$q$J}})=(a,0,b)$.
There is no difficulty in this limit.
The basic data are:
\begin{align}
  &  q^{\bm{\lambda}}=(a,b),\quad
  \bm{\delta}=(1,1),\quad\kappa=q^{-1},\quad
  0<a<q^{-1},\quad b<0,\n
  &\mathcal{E}(n)\eqdef q^{-n}-1,\quad
  \eta^{(+)}(x;\bm{\lambda})\eqdef aq^{x+1},\quad
  \eta^{(-)}(x;\bm{\lambda})\eqdef bq^{x+1},  
  \label{bqleta}\\
  &B^{\text{J}}(\eta;\bm{\lambda})\eqdef\eta^{-2}(-abq)(1-\eta),\quad
  D^{\text{J}}(\eta;\bm{\lambda})\eqdef\eta^{-2}(aq-\eta)(\eta-bq).
  \label{BtDtEn_bqL}
\end{align}
The potential functions \eqref{tilddef} are
\begin{align}
  &B^{(+)}(x;\bm{\lambda})=-a^{-1}bq^{-2x-1}(1-aq^{x+1}),\quad
  D^{(+)}(x;\bm{\lambda})=q^{-2x}(1-q^x)(q^x-a^{-1}b),
  \label{bqlB+}\\
  &B^{(-)}(x;\bm{\lambda})=-ab^{-1}q^{-2x-1}(1-bq^{x+1}),\quad
  D^{(-)}(x;\bm{\lambda})=q^{-2x}(1-q^x)(q^x-ab^{-1}).
  \label{bqlB-}
\end{align}
It is easy to verify the triangularity \eqref{gentri} and the
eigenvalues $\mathcal{E}(n)$ \eqref{bqleta}.
The sinusoidal coordinate $\eta^{(\pm)}(x;\bm{\lambda})$ satisfies
\eqref{etaprop}--\eqref{etaprop2}.

The big $q$-Laguerre polynomial $P_n(\eta;a,b;q)$ ($n\in\mathbb{Z}_{\geq 0}$),
\begin{align}
  P_n(\eta;\bm{\lambda})&\eqdef P_n(\eta;a,b;q)
  \eqdef{}_3\phi_2\Bigl(\genfrac{}{}{0pt}{}
  {q^{-n},\,0,\,\eta}{aq,\,bq}\!\Bigm|\!q\,;q\Bigr)\n
  &=\frac{(-b)^nq^{\frac12n(n+1)}}{(bq;q)_n}
  \,{}_2\phi_1\Bigl(\genfrac{}{}{0pt}{}
  {q^{-n},\,aq\eta^{-1}}{aq}\!\Bigm|\!q\,;b^{-1}\eta\Bigr),\quad
  P_n(1;\bm{\lambda})=1,
  \label{bqlpol2}
\end{align}
is the polynomial solution of degree $n$ in $\eta$ of the second order
difference equation \eqref{diffeq} with $B^{\text{J}}(\eta)$ and
$D^{\text{J}}(\eta)$ in \eqref{BtDtEn_bqL}.
It also satisfies the same type of three term recurrence relation as that
of the b$q$J \eqref{bqj3} with the coefficients
\begin{equation}
  A_n(\bm{\lambda})\eqdef-(1-aq^{n+1})(1-bq^{n+1}),\quad
  C_n(\bm{\lambda})\eqdef abq^{n+1}(1-q^{n}).
  \label{bqlACn}
\end{equation}
By setting $\eta\to aq^{x+1}$ in the second expression \eqref{bqlpol2},
we see that the big $q$-Laguerre polynomial is proportional to the dual
$q$-Meixner polynomial \eqref{dualqmeixner0} introduced in \S\,\ref{sec:naive}.

The eigenvectors of the two Hamiltonians $\mathcal{H}^{(\pm)}(\bm{\lambda})$
have the forms \eqref{phine} with \eqref{cPnpm},
in which $\phi^{(\pm)}_0(x;\bm{\lambda})$ are calculated to be 
\begin{align}
  &\phi^{(\pm)}_0(x;\bm{\lambda})^2
  =\pm A^{(\pm)}\eta^{(\pm)}(x;\bm{\lambda})
  \frac{\bigl(a^{-1}\eta^{(\pm)}(x;\bm{\lambda}),
  b^{-1}\eta^{(\pm)}(x;\bm{\lambda});q\bigr)_{\infty}}
  {\bigl(\eta^{(\pm)}(x;\bm{\lambda});q\bigr)_{\infty}},
  \label{bqL:phipm0}\\
  &A^{(+)}=\frac{\phi^{(+)}_0(0;\bm{\lambda})^2(aq;q)_{\infty}}
  {aq(q,ab^{-1}q;q)_{\infty}},\quad
  A^{(-)}=\frac{\phi^{(-)}_0(0;\bm{\lambda})^2(bq;q)_{\infty}}
  {-bq(q,a^{-1}bq;q)_{\infty}}.\nonumber
\end{align}
As shown in \S\,\ref{sec:bqJ}, the choice $A^{(+)}=A^{(-)}=1$ \eqref{A+=A-=1}
leads to the self-adjoint two component Hamiltonian system
$\mathcal{H}^{(\pm)}(\bm{\lambda})$ with the
potential functions $B^{(\pm)}(x;\bm{\lambda})$ and
$D^{(\pm)}(x;\bm{\lambda})$ \eqref{bqlB+}--\eqref{bqlB-}.
The ground state vectors $\phi^{(\pm)}_0(x;\bm{\lambda})>0$
\eqref{bqL:phipm0} read explicitly
\begin{equation}
  \phi^{(+)}_0(x;\bm{\lambda})^2=aq^{x+1}
  \frac{(q^{x+1},ab^{-1}q^{x+1};q)_{\infty}}{(aq^{x+1};q)_{\infty}},
  \ \ \phi^{(-)}_0(x;\bm{\lambda})^2=-bq^{x+1}
  \frac{(q^{x+1},a^{-1}bq^{x+1};q)_{\infty}}{(bq^{x+1};q)_{\infty}}.
  \label{bqL:phi0}
\end{equation}
We arrive at the orthogonality relation with the Jackson integral measure
\begin{align}
  \bin{\bm{\phi}_n}{\bm{\phi}_m}
  &=\bigl(\phi^{(+)}_n,\phi^{(+)}_m\bigr)
  +\bigl(\phi^{(-)}_n,\phi^{(-)}_m\bigr)
  =\frac{\delta_{nm}}{d_n(\bm{\lambda})^2}
  \label{bqlortho}\\
  &=aq\sum_{k=0}^{\infty}\frac{(q^{k+1},ab^{-1}q^{k+1};q)_{\infty}}
  {(aq^{k+1};q)_{\infty}}
  P_n(aq^{k+1};\bm{\lambda})P_m(aq^{k+1};\bm{\lambda})q^k\n
  &\quad-bq\sum_{k=0}^{\infty}\frac{(q^{k+1},a^{-1}bq^{k+1};q)_{\infty}}
  {(bq^{k+1};q)_{\infty}}
  P_n(bq^{k+1};\bm{\lambda})P_m(bq^{k+1};\bm{\lambda})q^k\n
  &=\frac{1}{1-q}\int_{bq}^{aq}\!\!d_qy\,
  \frac{(a^{-1}y,b^{-1}y;q)_{\infty}}{(y;q)_{\infty}}
  P_n(y;\bm{\lambda})P_m(y;\bm{\lambda}),
\end{align}
where the normalisation constant $d_n(\bm{\lambda})^2$ is
\begin{equation}
  d_n(\bm{\lambda})^2=
  (-abq^2)^{-n}q^{-\frac12n(n-1)}\frac{(aq,bq;q)_n}{(q;q)_n}
  \times d_0(\bm{\lambda})^2,\quad
  d_0(\bm{\lambda})^2=\frac{(aq,bq;q)_{\infty}}
  {aq(q,ab^{-1}q,a^{-1}b;q)_{\infty}}.
  \label{bqL:dn}
\end{equation}
The system is shape invariant, since the potential functions
$B^{(\pm)}(x;\bm{\lambda})$ and $D^{(\pm)}(x;\bm{\lambda})$
satisfy \eqref{shapeinv1cond1}--\eqref{shapeinv1cond2}.
The forward/backward shift relations for the b$q$L have the form
\eqref{FBe}--\eqref{FBPne} with the auxiliary functions
$\varphi^{(\pm)}(x;\bm{\lambda})$,
\begin{equation}
  \varphi^{(\pm)}(x;\bm{\lambda})\eqdef
  \frac{q\eta^{(\pm)}(x;\bm{\lambda})}{(1-aq)(1-bq)}.
  \label{bqL:varphipm}
\end{equation}
Under the parameter substitution (involution) $(a,b)\to(b,a)$,
the $(+)$ and $(-)$ systems are interchanged as \eqref{bqJ:sym}.

\subsubsection{dual big $q$-Laguerre}
\label{sec:dbqL}

The two types of the dual big $q$-Laguerre (db$q$L) polynomials are
\begin{align}
  \check{Q}^{(+)}_x(n;\bm{\lambda})&\eqdef
  Q^{(+)}_x\bigl(\mathcal{E}(n);\bm{\lambda}\bigr)\eqdef
  {}_2\phi_1\Bigl(\genfrac{}{}{0pt}{}{q^{-n},\,q^{-x}}{aq}\!\Bigm|\!q\,;
  ab^{-1}q^{x+1}\Bigr),
  \label{bqlQx+}\\
  \check{Q}^{(-)}_x(n;\bm{\lambda})&\eqdef
  Q^{(-)}_x\bigl(\mathcal{E}(n);\bm{\lambda}\bigr)\eqdef
  {}_2\phi_1\Bigl(\genfrac{}{}{0pt}{}{q^{-n},\,q^{-x}}{bq}\!\Bigm|\!q\,;
  a^{-1}bq^{x+1}\Bigr).
  \label{bqlQx-}
\end{align}
They are degree $x$ polynomials in $\mathcal{E}(n)$ \eqref{bqleta} and
related to $\check{P}^{(\pm)}_n(x;\bm{\lambda})$ as \eqref{bqjdualpoly} with
\begin{equation}
  \alpha^{(+)}_n(\bm{\lambda})\eqdef
  \frac{(-b)^nq^{\frac12n(n+1)}}{(bq;q)_n},\quad
  \alpha^{(-)}_n(\bm{\lambda})\eqdef
  \frac{(-a)^nq^{\frac12n(n+1)}}{(aq;q)_n}.
  \label{bqL:alphapmn}
\end{equation}
Following the argument in \S\,\ref{sec:dbqJ}, $A^{(\pm)}_n(\bm{\lambda})$
and $C^{(\pm)}_n(\bm{\lambda})$ \eqref{AneCne}
and $\mathcal{H}^{\text{d}\,(\pm)}(\bm{\lambda})$ \eqref{dbqJham} are
introduced:
\begin{alignat*}{2}
  A^{(+)}_n(\bm{\lambda})&=bq^{n+1}(1-aq^{n+1}),
  &\ \ C^{(+)}_n(\bm{\lambda})&=-aq(1-q^n)(1-bq^n),\\
  A^{(-)}_n(\bm{\lambda})&=-aq^{n+1}(1-bq^{n+1}),
  &\ \ C^{(-)}_n(\bm{\lambda})&=bq(1-q^n)(1-aq^n).
\end{alignat*}
They satisfy \eqref{bqJ:AnCn+1}--\eqref{bqJ:An+Cn}, which implies
\eqref{HmHprel}.
The eigenvectors of $\mathcal{H}^{\text{d}\,(\pm)}(\bm{\lambda})$ are
\begin{align}
  &\phi^{\text{d}\,(\pm)}_x(n;\bm{\lambda})\eqdef
  \phi^{\text{d}\,(\pm)}_0(n;\bm{\lambda})
  \check{Q}^{(\pm)}_x(n;\bm{\lambda}),
  \label{bqL:phidpm}\\
  &\phi^{\text{d}\,(\pm)}_0(n;\bm{\lambda})=
  (\pm1)^n\prod_{m=0}^{n-1}\sqrt{\frac{A^{(\pm)}_m(\bm{\lambda})}
  {C^{(\pm)}_{m+1}(\bm{\lambda})}}
  =\alpha^{(\pm)}_n(\bm{\lambda})
  \frac{d_n(\bm{\lambda})}{d_0(\bm{\lambda})},
  \ \ (\pm1)^n\phi^{\text{d}\,(\pm)}_0(n;\bm{\lambda})>0,\n
  &\mathcal{H}^{\text{d}\,(\pm)}(\bm{\lambda})
  \phi^{\text{d}\,(\pm)}_x(n;\bm{\lambda})
  =\mathcal{E}^{\text{d}\,(\pm)}(x;\bm{\lambda})
  \phi^{\text{d}\,(\pm)}_x(n;\bm{\lambda}),\n
  &\mathcal{H}^{\text{d}\,(\pm)}(\bm{\lambda})
  \phi^{\text{d}\,(\mp)}_x(n;\bm{\lambda})
  =\mathcal{E}^{\prime\,\text{d}\,(\pm)}(x;\bm{\lambda})
  \phi^{\text{d}\,(\mp)}_x(n;\bm{\lambda}),\n
  &\mathcal{E}^{\text{d}\,(\pm)}(x;\bm{\lambda})\eqdef\left\{
  \begin{array}{rl}
  aq(1-q^x)&:(+) \\
  -bq(1-q^x)&:(-)
  \end{array}
  \right.,\quad
  \mathcal{E}^{\prime\,\text{d}\,(\pm)}(x;\bm{\lambda})\eqdef\left\{
  \begin{array}{ll}
  q(a-bq^x)&:(+)\\
  q(-b+aq^x)&:(-)
  \end{array}\right..
  \label{bqL:E'd}
\end{align}
They satisfy the dual orthogonality relations in the same form as
\eqref{Q-ortho1}--\eqref{Q-ortho3},
in which $\phi^{(\pm)}_0(x;\bm{\lambda})^2$ and $d_0(\bm{\lambda})^2$ are
defined in \eqref{bqL:phi0} and \eqref{bqL:dn}.
The orthogonality relation \eqref{Q-ortho3} means simply
\begin{equation}
  \sum_{n=0}^{\infty}(-1)^n\frac{q^{\frac12n(n-1)}}{(q;q)_n}q^{-kn}=0
  \ \ (k\in\mathbb{Z}_{\geq 0}),\quad
  \phi^{\text{d}\,(+)}_0(n;\bm{\lambda})\phi^{\text{d}\,(-)}_0(n;\bm{\lambda})
  =(-1)^n\frac{q^{\frac12n(n-1)}}{(q;q)_n},
  \label{bqlQ-ortho4}
\end{equation}
which can be shown by \cite{koeswart}
\begin{equation*}
  \sum_{n=0}^{\infty}\frac{q^{\frac12n(n-1)}}{(q;q)_n}z^n=(-z;q)_{\infty}.
\end{equation*}

As expected from the big $q$-Laguerre and the dual $q$-Meixner
correspondence, the dual big $q$-Laguerre polynomials
\eqref{bqlQx+}--\eqref{bqlQx-} are equal to the $q$-Meixner polynomial
and its completeness partner \eqref{qM:Pnm} under certain parameter
identification.

\subsubsection{Al-Salam-Carlitz $\I$}
\label{sec:ASCI}

The Al-Salam-Carlitz $\I$ polynomial is obtained from the big $q$-Laguerre
(b$q$L) polynomial by setting $(a^{\text{b$q$L}},b^{\text{b$q$L}})=(t,ta)$,
$\eta^{\text{b$q$L}}=tq\eta$ and taking $t\to 0$ limit.
In this limit we need appropriate overall rescalings for various quantities.
The basic data are:
\begin{align}
  &q^{\bm{\lambda}}=a,\quad
  \bm{\delta}=0,\quad\kappa=q^{-1},\quad a<0,\n
  &\mathcal{E}(n)=q^{-n}-1,\quad
  \eta^{(+)}(x;\bm{\lambda})\eqdef q^x,\quad
  \eta^{(-)}(x;\bm{\lambda})\eqdef aq^x,\n
  &B^{\text{J}}(\eta;\bm{\lambda})\eqdef\eta^{-2}(-aq^{-1}),\quad
  D^{\text{J}}(\eta;\bm{\lambda})\eqdef\eta^{-2}(1-\eta)(\eta-a).
  \label{BtDtEn_ASCI}
\end{align}
The potential functions \eqref{tilddef} are
\begin{alignat}{2}
  B^{(+)}(x;\bm{\lambda})&=-aq^{-2x-1},&
  D^{(+)}(x;\bm{\lambda})&=q^{-2x}(1-q^x)(q^x-a),
  \label{ASCIB+}\\
  B^{(-)}(x;\bm{\lambda})&=-a^{-1}q^{-2x-1},&\quad
  D^{(-)}(x;\bm{\lambda})&=q^{-2x}(1-q^x)(q^x-a^{-1}).
  \label{ASCIB-}
\end{alignat}
The two sinusoidal coordinates satisfy \eqref{etaprop1}--\eqref{etaprop2}.
Throughout this paper, our definition of Al-Salam-Carlitz $\I$ (ASC\,$\I$)
polynomial is a simple $q$-hypergeometric function
\begin{equation}
  P_n(\eta;\bm{\lambda})\eqdef{}_2\phi_1\Bigl(\genfrac{}{}{0pt}{}
  {q^{-n},\,\eta^{-1}}{0}\!\Bigm|\!q\,;a^{-1}q\eta\Bigr)
  =(-a)^{-n}q^{-\frac12n(n-1)}U^{(a)}_n(\eta;q),\quad
  P_n(1;\bm{\lambda})=1,
  \label{ASCIpol}
\end{equation}
which is obtained by rescaling the conventional Al-Salam-Carlitz $\I$
polynomial $U^{(a)}_n(\eta;q)$ ($n\in\mathbb{Z}_{\geq 0}$) \cite{koeswart}.
It is the polynomial solution of degree $n$ in $\eta$ of the difference
equation \eqref{diffeq} with $B^{\text{J}}(\eta)$ and
$D^{\text{J}}(\eta)$ in \eqref{BtDtEn_ASCI}.
With the rescaling our ASC\,$\I$ satisfy the same form of the three term
recurrence relation as that of b$q$J \eqref{bqj3} and its coefficients are
\begin{equation}
  A_n(\bm{\lambda})\eqdef aq^n,\quad C_n\eqdef -(1-q^n).
  \label{ASCI:An,Cn}
\end{equation}
The forward/backward shift relations are
\begin{align*}
  &\frac{-a}{q\eta}
  \bigl(P_n(\eta;\bm{\lambda})-P_n(q\eta;\bm{\lambda})\bigr)
  =\mathcal{E}(n)P_{n-1}(\eta;\bm{\lambda}+\bm{\delta}),\\
  &\bigl({B}^{\text{J}}(\eta;\bm{\lambda})
  P_{n-1}(\eta;\bm{\lambda}+\bm{\delta})
  -q^{-1}{D}^{\text{J}}(\eta;\bm{\lambda})
  P_{n-1}(q^{-1}\eta;\bm{\lambda}+\bm{\delta})\bigr)\frac{q\eta}{-a}
  =P_n(\eta;\bm{\lambda}).
\end{align*}
We strongly advocate this natural definition of Al-Salam-Carlitz $\I$
polynomial \eqref{ASCIpol}, since many important formulas take the same
form as those for the b$q$J family.

By setting $\eta\to q^{x}$ ($-aq^x$) in \eqref{ASCIpol}, we see that the
ASC\,$\I$ polynomial $P_n(\eta;\bm{\lambda})$ is equal to the dual $q$-Charlier
polynomial \eqref{dualqchar-1} (\eqref{dualqchar-}) introduced in
\S\,\ref{sec:dqChar} with parameter identification $a\to-a$.

The eigenvectors of the two Hamiltonians $\mathcal{H}^{(\pm)}(\bm{\lambda})$
are \eqref{phine} with \eqref{cPnpm},
in which $\phi^{(\pm)}_0(x;\bm{\lambda})$ are calculated as
\begin{align}
  &\phi^{(\pm)}_0(x;\bm{\lambda})^2
  =\pm A^{(\pm)}\eta^{(\pm)}(x;\bm{\lambda})
  (q\eta^{(\pm)}(x;\bm{\lambda}),
  a^{-1}q\eta^{(\pm)}(x;\bm{\lambda});q\bigr)_{\infty},
  \label{ASCI:phipm0}\\
  &A^{(+)}=\frac{\phi^{(+)}_0(0;\bm{\lambda})^2}
  {(q,a^{-1}q;q)_{\infty}},\quad
  A^{(-)}=\frac{\phi^{(-)}_0(0;\bm{\lambda})^2}
  {-a(q,aq;q)_{\infty}}.\nonumber
\end{align}
As shown in \S\,\ref{sec:bqJ}, the choice $A^{(+)}=A^{(-)}=1$ \eqref{A+=A-=1}
leads to the self-adjoint two component Hamiltonian system
$\mathcal{H}^{(\pm)}(\bm{\lambda})$ with the
potential functions $B^{(\pm)}(x;\bm{\lambda})$ and
$D^{(\pm)}(x;\bm{\lambda})$ \eqref{ASCIB+}--\eqref{ASCIB-}.
The ground state vectors $\phi^{(\pm)}_0(x;\bm{\lambda})>0$
\eqref{ASCI:phipm0} read explicitly
\begin{equation}
  \phi^{(+)}_0(x;\bm{\lambda})^2
  =q^x(q^{x+1},a^{-1}q^{x+1};q)_{\infty},\quad
  \phi^{(-)}_0(x;\bm{\lambda})^2
  =-aq^x(q^{x+1},aq^{x+1};q)_{\infty}.
  \label{ASCI:phi0}
\end{equation}
We arrive at the orthogonality relation with the Jackson integral measure
\begin{align}
  \bin{\bm{\phi}_n}{\bm{\phi}_m}
  &=\bigl(\phi^{(+)}_n,\phi^{(+)}_m\bigr)
  +\bigl(\phi^{(-)}_n,\phi^{(-)}_m\bigr)
  =\frac{\delta_{nm}}{d_n(\bm{\lambda})^2}
  \label{ASCI:ortho}\\
  &=\sum_{k=0}^{\infty}(q^{k+1},a^{-1}q^{k+1};q)_{\infty}
  P_n(q^k;\bm{\lambda})P_m(q^k;\bm{\lambda})q^k\n
  &\quad-a\sum_{k=0}^{\infty}(q^{k+1},aq^{k+1};q)_{\infty}
  P_n(aq^k;\bm{\lambda})P_m(aq^k;\bm{\lambda})q^k\n
  &=\frac{1}{1-q}\int_a^1\!\!d_qy\,
  (qy,a^{-1}qy;q)_{\infty}
  P_n(y;\bm{\lambda})P_m(y;\bm{\lambda}),
\end{align}
where the normalisation constant $d_n(\bm{\lambda})^2$ is
\begin{equation}
  d_n(\bm{\lambda})^2=
  \frac{(-a)^nq^{\frac12n(n-1)}}{(q;q)_n}
  \times d_0(\bm{\lambda})^2,\quad
  d_0(\bm{\lambda})^2=\frac{1}{(q,a,a^{-1}q;q)_{\infty}}.
  \label{ASCI:dn}
\end{equation}
The system is shape invariant, since the potential functions
$B^{(\pm)}(x;\bm{\lambda})$ and $D^{(\pm)}(x;\bm{\lambda})$
satisfy \eqref{shapeinv1cond1}--\eqref{shapeinv1cond2}.
The forward/backward shift relations for the ASC\,$\I$ have the form
\eqref{FBe}--\eqref{FBPne} with the auxiliary functions
$\varphi^{(\pm)}(x;\bm{\lambda})$,
\begin{equation}
  \varphi^{(\pm)}(x;\bm{\lambda})\eqdef
  \frac{q\eta^{(\pm)}(x;\bm{\lambda})}{-a}.
  \label{ASCI:varphipm}
\end{equation}

\subsubsection{dual Al-Salam-Carlitz $\I$}
\label{sec:dASCI}

The two types of the dual Al-Salam-Carlitz $\I$ (dASC\,$\I$) polynomials are
\begin{align}
  \check{Q}^{(+)}_x(n;\bm{\lambda})&\eqdef
  Q^{(+)}_x\bigl(\mathcal{E}(n);\bm{\lambda}\bigr)\eqdef
  {}_2\phi_1\Bigl(\genfrac{}{}{0pt}{}{q^{-n},\,q^{-x}}{0}\!\Bigm|\!q\,;
  a^{-1}q^{x+1}\Bigr),
  \label{ASCIQx+}\\
  \check{Q}^{(-)}_x(n;\bm{\lambda})&\eqdef
  Q^{(-)}_x\bigl(\mathcal{E}(n);\bm{\lambda}\bigr)\eqdef 
  {}_2\phi_1\Bigl(\genfrac{}{}{0pt}{}{q^{-n},\,q^{-x}}{0}\!\Bigm|\!q\,;
  aq^{x+1}\Bigr).
  \label{ASCIQx-}
\end{align}
They are degree $x$ polynomials in $\mathcal{E}(n)$ and related to
$\check{P}^{(\pm)}_n(x;\bm{\lambda})$ as \eqref{bqjdualpoly} with
\begin{equation}
  \alpha^{(+)}_n(\bm{\lambda})\eqdef 1,\quad
  \alpha^{(-)}_n(\bm{\lambda})\eqdef a^{-n}.
  \label{ASCI:alphapmn}
\end{equation}
Following the argument in \S\,\ref{sec:dbqJ}, $A^{(\pm)}_n(\bm{\lambda})$
and $C^{(\pm)}_n(\bm{\lambda})$ \eqref{AneCne}
and $\mathcal{H}^{\text{d}\,(\pm)}(\bm{\lambda})$ \eqref{dbqJham} are
introduced:
\begin{equation*}
  A^{(+)}_n(\bm{\lambda})=aq^n,\quad
  A^{(-)}_n(\bm{\lambda})=-q^n,\quad
  C^{(+)}_n(\bm{\lambda})=-(1-q^n),\quad C^{(-)}_n(\bm{\lambda})=a(1-q^n).
\end{equation*}
They satisfy \eqref{bqJ:AnCn+1}--\eqref{bqJ:An+Cn}, which implies
\eqref{HmHprel}.
The eigenvectors of $\mathcal{H}^{\text{d}\,(\pm)}(\bm{\lambda})$ are
\begin{align}
  &\phi^{\text{d}\,(\pm)}_x(n;\bm{\lambda})\eqdef
  \phi^{\text{d}\,(\pm)}_0(n;\bm{\lambda})
  \check{Q}^{(\pm)}_x(n;\bm{\lambda}),
  \label{ASCI:phidpm}\\
  &\phi^{\text{d}\,(\pm)}_0(n;\bm{\lambda})=
  (\pm1)^n\prod_{m=0}^{n-1}\sqrt{\frac{A^{(\pm)}_m(\bm{\lambda})}
  {C^{(\pm)}_{m+1}(\bm{\lambda})}}
  =\alpha^{(\pm)}_n(\bm{\lambda})
  \frac{d_n(\bm{\lambda})}{d_0(\bm{\lambda})},
  \ \ (\pm1)^n\phi^{\text{d}\,(\pm)}_0(n;\bm{\lambda})>0,\\
  &\mathcal{H}^{\text{d}\,(\pm)}(\bm{\lambda})
  \phi^{\text{d}\,(\pm)}_x(n;\bm{\lambda})
  =\mathcal{E}^{\text{d}\,(\pm)}(x;\bm{\lambda})
  \phi^{\text{d}\,(\pm)}_x(n;\bm{\lambda}),\n
  &\mathcal{H}^{\text{d}\,(\pm)}(\bm{\lambda})
  \phi^{\text{d}\,(\mp)}_x(n;\bm{\lambda})
  =\mathcal{E}^{\prime\,\text{d}\,(\pm)}(x;\bm{\lambda})
  \phi^{\text{d}\,(\mp)}_x(n;\bm{\lambda}),\n
  &\mathcal{E}^{\text{d}\,(\pm)}(x;\bm{\lambda})\eqdef\left\{
  \begin{array}{rl}
  1-q^x\phantom{)}&:(+)\\
  -a(1-q^x)&:(-)
  \end{array}\right.,\quad
  \mathcal{E}^{\prime\,\text{d}\,(\pm)}(x;\bm{\lambda})\eqdef\left\{
  \begin{array}{ll}
   1-aq^x&:(+)\\
   -a+q^x&:(-)
  \end{array}\right..
  \label{ASCI:E'd}
\end{align}
They satisfy the dual orthogonality relations in the same form as
\eqref{Q-ortho1}--\eqref{Q-ortho3},
in which $\phi^{(\pm)}_0(x;\bm{\lambda})^2$ and $d_0(\bm{\lambda})^2$ are
given in \eqref{ASCI:phi0} and \eqref{ASCI:dn}.
As in the big $q$-Laguerre case, the orthogonality relation \eqref{Q-ortho3}
means the same relation as \eqref{bqlQ-ortho4}.
In line with the Al-Salam-Carlitz $\I$ and the dual $q$-Charlier
correspondence, the dual Al-Salam-Carlitz $\I$ polynomials
\eqref{ASCIQx+}--\eqref{ASCIQx-} are equal to the $q$-Charlier polynomial
and its completeness supplement \eqref{qCsup0}--\eqref{qCsup}.

\subsubsection{discrete $q$-Hermite $\I$}
\label{sec:dqHI}

The discrete $q$-Hermite polynomial $h_n(\eta;q)$ is obtained from the
Al-Salam-Carlitz $\I$ (ASC\,$\I$) polynomial by setting
$a^{\text{ASC\,$\I$}}=-1$: $h_n(\eta;q)=U^{(-1)}_n(\eta;q)$.
Various formulas are easily obtained from those of ASC\,$\I$ by setting
$a^{\text{ASC\,$\I$}}=-1$.

\subsection{Complete set involving $q$-Meixner polynomials}
\label{sec:compqm}

By constructing the duals of the dual $q$-Meixner polynomials in a similar
way as the dual big $q$-Jacobi and big $q$-Laguerre cases developed in
\S\,\ref{sec:dbqJ}, \S\,\ref{sec:dbqL}, we obtain the complete set of
orthogonal vectors involving the $q$-Meixner polynomials.
Since the dual of $\check{P}^{(+)}_n(x)$ \eqref{dqMpm4} is the original
$q$-Meixner \eqref{qMpol}, that of $\check{P}^{(-)}_n(x)$
\begin{equation*}
  {}_2\phi_1\Bigl(
  \genfrac{}{}{0pt}{}{q^{-n},\,-c^{-1}q^{-x}}{bq}\!\Bigm|\!q\,;q^{x+1}\Bigr)
  =\frac{(-bcq;q)_n}{(-c)^n(bq;q)_n}{}_2\phi_1\Bigl(
  \genfrac{}{}{0pt}{}{q^{-n},\,q^{-x}}{-bcq}\!\Bigm|\!q\,;-cq^{x+1}\Bigr),
\end{equation*}
provides the supplementary orthogonal vectors. By interchanging
$x\leftrightarrow n$ of $\check{P}^{(\pm)}_n(x)$, we obtain two polynomials
in $\eta(x)=q^{-x}-1$,
\begin{align}
  \check{P}_n(x;\bm{\lambda})
  \eqdef P_n\bigl(\eta(x);\bm{\lambda}\bigr)
  &\eqdef{}_2\phi_1\Bigl(\genfrac{}{}{0pt}{}
  {q^{-n},q^{-x}}{bq}\!\!\Bigm|\!q\,;-c^{-1}q^{n+1}\Bigr)
  =M_n(q^{-x};b,c;q),
  \label{qM:Pnp}\\
  \check{P}^{(-)}_n(x;\bm{\lambda})
  \eqdef P^{(-)}_n\bigl(\eta(x);\bm{\lambda}\bigr)
  &\eqdef{}_2\phi_1\Bigl(\genfrac{}{}{0pt}{}
  {q^{-n},q^{-x}}{-bcq}\!\!\Bigm|\!q\,;-cq^{n+1}\Bigr)
  =M_n(q^{-x};-bc,c^{-1};q),
  \label{qM:Pnm}
\end{align}
{\em i.e.\/} the original $q$-Meixner \eqref{qMpol} and its supplement.
The parameter change (involution)
\begin{equation*}
  (b,c)\to(-bc,c^{-1}),
\end{equation*}
exchanges these two and it also provides $B^{(-)}(x)$, $D^{(-)}(x)$,
$A^{(-)}_n$ and $C^{(-)}_n$ governing the latter
\begin{align}
  &B^{(-)}(x;\bm{\lambda})\eqdef c^{-1}q^x(1+bcq^{x+1}),\quad
  D^{(-)}(x;\bm{\lambda})\eqdef (1-q^x)(1-bq^x),\n
  &A^{(-)}_n(\bm{\lambda})\eqdef -c^{-1}q^{-2n-1}(1+bcq^{n+1}),\quad
  C^{(-)}_n(\bm{\lambda})\eqdef -q^{-2n}(1-q^n)(q^n+c^{-1}),
  \label{qM:Amn,Cmn}
\end{align}
from $B(x)$, $D(x)$, $A_n$ and $C_n$ for the original $q$-Meixner
\eqref{qMB&D}, \eqref{3coef}. By noting
\begin{align*}
  B(x;\bm{\lambda})D(x+1;\bm{\lambda})
  &=c^2B^{(-)}(x;\bm{\lambda})D^{(-)}(x+1;\bm{\lambda}),\\
  B(x;\bm{\lambda})+D(x;\bm{\lambda})-1
  &=-c\bigl(B^{(-)}(x;\bm{\lambda})+D^{(-)}(x;\bm{\lambda})-1\bigr),
\end{align*}
we find that the two Hamiltonians $\mathcal{H}(\bm{\lambda})$ and
$\mathcal{H}^{(-)}(\bm{\lambda})$
\begin{align*}
  \mathcal{H}(\bm{\lambda})
  &=-\sqrt{B(x;\bm{\lambda})D(x+1;\bm{\lambda})}\,e^{\partial}
  -\sqrt{B(x-1;\bm{\lambda})D(x;\bm{\lambda})}\,e^{-\partial}
  +B(x;\bm{\lambda})+D(x;\bm{\lambda}),\n
  \mathcal{H}^{(-)}(\bm{\lambda})&=
  \sqrt{B^{(-)}(x;\bm{\lambda})D^{(-)}(x+1;\bm{\lambda})}\,e^{\partial}
  +\sqrt{B^{(-)}(x-1;\bm{\lambda})D^{(-)}(x;\bm{\lambda})}\,e^{-\partial}\n
  &\quad+B^{(-)}(x;\bm{\lambda})+D^{(-)}(x;\bm{\lambda}),
\end{align*}
are linearly related
\begin{equation*}
  c\mathcal{H}^{(-)}=-\mathcal{H}+1+c.
\end{equation*}
The eigenvectors are
\begin{align} 
  &\phi_n(x;\bm{\lambda})=\phi_0(x;\bm{\lambda})
  \check{P}_n(x;\bm{\lambda}),
  \ \ \phi_0(x;\bm{\lambda})^2
  =c^xq^{\frac12x(x-1)}\frac{(bq;q)_x}{(q,-bcq;q)_x},
  \ \ \phi_0(x;\bm{\lambda})>0,
  \label{qMori}\\
  &\phi^{(-)}_n(x;\bm{\lambda})=\phi^{(-)}_0(x;\bm{\lambda})
  \check{P}^{(-)}_n(x;\bm{\lambda}),\quad
  \phi^{(-)}_0(x;\bm{\lambda})=(-1)^x\prod_{y=0}^{x-1}
  \sqrt{\frac{B^{(-)}(y;\bm{\lambda})}{D^{(-)}(y+1;\bm{\lambda})}},\n
  &\qquad\phi^{(-)}_0(x;\bm{\lambda})^2
  =c^{-x}q^{\frac12x(x-1)}\frac{(-bcq;q)_x}{(q,bq;q)_x},\quad
  (-1)^x\phi^{(-)}_0(x;\bm{\lambda})>0,
  \label{qMsup}\\
  &\mathcal{H}(\bm{\lambda})\phi_n(x;\bm{\lambda})
  =\mathcal{E}(n)\phi_n(x;\bm{\lambda}),
  \ \mathcal{H}^{(-)}(\bm{\lambda})\phi^{(-)}_n(x;\bm{\lambda})
  =\mathcal{E}(n)\phi^{(-)}_n(x;\bm{\lambda}),
  \ \mathcal{E}(n)=1-q^n.\nonumber
\end{align}
The original Hamiltonian $\mathcal{H}(\bm{\lambda})$ governing the
$q$-Meixner polynomials has another infinite set of eigenvectors
$\{\phi^{(-)}_n(x;\bm{\lambda})\}$ \eqref{qMsup},
\begin{equation}
  \mathcal{H}(\bm{\lambda})\phi^{(-)}_n(x;\bm{\lambda})
  =\mathcal{E}^{\prime}(n;\bm{\lambda})\phi^{(-)}_n(x;\bm{\lambda}),\quad
  \mathcal{E}^{\prime}(n;\bm{\lambda})\eqdef 1+cq^n, 
  \label{qMsupeig}
\end{equation}
which lies above the original eigenvectors.
Similarly we have
\begin{equation*}
  \mathcal{H}^{(-)}(\bm{\lambda})\phi^{(+)}_n(x;\bm{\lambda})
  =\mathcal{E}^{\prime\,(-)}(n;\bm{\lambda})
  \phi^{(+)}_n(x;\bm{\lambda}),\quad
  \mathcal{E}^{\prime\,(-)}(n;\bm{\lambda})\eqdef 1+c^{-1}q^n.
\end{equation*}
The orthogonality relations are ($n,m=0,1,\ldots$)
\begin{align}
  (\phi_n,\phi_m)&=\sum_{x=0}^{\infty}
  \phi_n(x;\bm{\lambda})\phi_m(x;\bm{\lambda})
  =\frac{\delta_{nm}}{d_n(\bm{\lambda})^2},
  \label{qM:orthorel0}\\
  \bigl(\phi^{(-)}_n,\phi^{(-)}_m\bigr)&=\sum_{x=0}^{\infty}
  \phi^{(-)}_n(x;\bm{\lambda})\phi^{(-)}_m(x;\bm{\lambda})
  =\frac{\delta_{nm}}{d^{(-)}_n(\bm{\lambda})^2},
  \ \ d^{(-)}_n(\bm{\lambda})\eqdef
  d_n(\bm{\lambda})\bigl|_{(b,c)\to(-bc,c^{-1})},\n
  \bigl(\phi_n,\phi^{(-)}_m\bigr)&=\sum_{x=0}^{\infty}
  \phi_n(x;\bm{\lambda})\phi^{(-)}_m(x;\bm{\lambda})=0.
  \label{qM:orthorel2}
\end{align}
The last orthogonality relation simply means
\begin{equation*}
  \sum_{x=0}^{\infty}(-1)^x\frac{q^{\frac12x(x-1)}}{(q;q)_x}q^{-kx}=0
  \ \ (k\in\mathbb{Z}_{\geq 0}),\quad
  \phi_0(x;\bm{\lambda})\phi^{(-)}_0(x;\bm{\lambda})
  =(-1)^x\frac{q^{\frac12x(x-1)}}{(q;q)_x},
\end{equation*}
which is the same as \eqref{bqlQ-ortho4}.

\subsubsection{complete set involving $q$-Charlier polynomials}
\label{sec:compqc}

The $q$-Charlier polynomials are obtained from those of $q$-Meixner by
setting $(b^{\text{$q$M}},c^{\text{$q$M}})=(0,a)$.
The complete set of orthogonal vectors involving the $q$-Charlier
polynomials are
\begin{align}
  &\phi_n(x;\bm{\lambda})=\phi_0(x;\bm{\lambda})
  {}_2\phi_1\Bigl(\genfrac{}{}{0pt}{}{q^{-n},\,q^{-x}}{0}\!\Bigm|\!q\,;
  -a^{-1}q^{n+1}\Bigr),
  \ \ \phi_0(x;\bm{\lambda})=\sqrt{\frac{a^xq^{\frac12x(x-1)}}{(q;q)_x}}\,,
  \label{qCsup0}\\
  &\phi^{(-)}_n(x;\bm{\lambda})=\phi^{(-)}_0(x;\bm{\lambda})
  {}_2\phi_1\Bigl(\genfrac{}{}{0pt}{}{q^{-n},\,q^{-x}}{0}\!\Bigm|\!q\,;
  -aq^{n+1}\Bigr),
  \ \ \phi^{(-)}_0(x;\bm{\lambda})
  =(-1)^x\sqrt{\frac{a^{-x}q^{\frac12x(x-1)}}{(q;q)_x}}\,,
  \label{qCsup}\\
  &\mathcal{H}(\bm{\lambda})\phi_n(x;\bm{\lambda})
  =\mathcal{E}(n)\phi_n(x;\bm{\lambda}),\quad
  \mathcal{E}(n)=1-q^n,\\
  &\mathcal{H}(\bm{\lambda})\phi^{(-)}_n(x;\bm{\lambda})
  =\mathcal{E}^{\prime}(n;\bm{\lambda})\phi^{(-)}_n(x;\bm{\lambda}),\quad
  \mathcal{E}^{\prime}(n;\bm{\lambda})\eqdef 1+aq^n.
  \label{qCsupeig}
\end{align}
The orthogonality relations have the same form as
\eqref{qM:orthorel0}--\eqref{qM:orthorel2} with $d_n(\bm{\lambda})$ given
in \eqref{qCdn} and
$d^{(-)}_n(\bm{\lambda})\eqdef d_n(\bm{\lambda})|_{a\to a^{-1}}$.
The third one $\bigl(\phi_n,\phi^{(-)}_m\bigr)=0$ is the same as
\eqref{bqlQ-ortho4}.

\section{Discrete $q$-Hermite $\II$}
\label{sec:dqHII}

In this section we discuss the discrete $q$-Hermite polynomial $\II$,
which is defined on the integer lattice of the full line $x\in\mathbb{Z}$. 
In order to place it under proper perspectives with respect to the other
orthogonal polynomials having Jackson integral measures, let us go back to
the general structure of the difference equations and the corresponding
potential functions \eqref{p1p2} giving rise to orthogonal polynomials
with Jackson integral type measures.
The essential point for the Jackson integral measure is that $p_2(\eta)$
($\eta\propto q^x$) has two distinct roots other than zero. If it has one
zero root, $\eta=0\leftrightarrow x=+\infty$, then $p_1(\eta)$ must
have one zero root, too. Otherwise the zero mode vector $\phi_0(x)$ would not
be square summable. This is the case for the little $q$-Jacobi (Laguerre)
polynomials. One could say that their orthogonality measures are of Jackson
integral type of the form $\int_0^1\!\!d_qx$.
If $p_2(\eta)$ is a linear function of $\eta$, the second Hamiltonian
cannot be introduced. The difference operator (Hamiltonian) is not
self-adjoint due to the rapid growth at $x=\infty$.
If $p_2$ is a constant, however, new possibility of achieving
self-adjointness emerges when $p_1$ is even. In this case the boundary
point of the variable $x$ disappears and it can take full integer values,
$-\infty<x<\infty$ and the scale of the variable $q^x$ can become arbitrary.
Due to the evenness of the potential functions, a polynomial solution
$P(cq^x)$ ($c\in\mathbb{R}_{>0}$) is always accompanied by another $P(-cq^x)$.
Thus two component Hamiltonian method can be employed to take care of the
rapid growth of the potential functions at $x\to+\infty$ and the exponential
increase of the eigenpolynomials $P(\pm cq^x)$ at $x\to-\infty$.
The fact that $x$ can take negative values makes it impossible to introduce
dual polynomials by the interchange $n\leftrightarrow x$
\eqref{duality1}--\eqref{duality2}, since $n$ is always a non-negative integer.

Let us now consider the general difference equation
\eqref{Hdef}--\eqref{A,Ad} on the full integer lattice.
Now the potential functions $B(x)$ and $D(x)$ are everywhere positive
\begin{equation}
  B(x)>0,\quad D(x)>0\ \ (x\in\mathbb{Z}).
  \label{BDcondinf}
\end{equation}
The zero mode equation $\mathcal{A}\phi_0(x)=0$ \eqref{phi0eq}, being
a two term recurrence relation, can be solved easily by starting from
$x=0$ into the positive and negative directions:
\begin{equation}
  \phi_0(x)=\phi_0(0)\times\left\{
  \begin{array}{ll}
  {\displaystyle\prod_{y=0}^{x-1}\sqrt{\frac{B(y)}{D(y+1)}}}
  &(x\in\mathbb{Z}_{\geq 0})\\
  {\displaystyle\prod_{y=0}^{-x-1}\sqrt{\frac{D(-y)}{B(-y-1)}}}
  &(x\in\mathbb{Z}_{<0})
  \end{array}\right..
  \label{phi0inf}
\end{equation}
The similarity transformed Hamiltonian $\widetilde{\mathcal{H}}$ has the
same expression as in the half line case \eqref{Ht}
\begin{equation*}
  \widetilde{\mathcal{H}}
  \eqdef\phi_0(x)^{-1}\circ\mathcal{H}\circ\phi_0(x)
  =B(x)(1-e^{\partial})+D(x)(1-e^{-\partial}).
\end{equation*}
The potential functions for the discrete $q$-Hermite $\II$ are
\begin{equation}
  B^{\text{J}}(\eta)\eqdef\eta^{-2}(1+\eta^2),\quad
  D^{\text{J}}(\eta)\eqdef\eta^{-2}q,\quad
  \mathcal{E}(n)=1-q^n,\quad\eta(x)\propto q^x.
  \label{BtDtEn_dqHII}
\end{equation}
It is obvious that the above $\widetilde{\mathcal{H}}$ is triangular with
respect to the basis $\{1,\eta, \eta^2,\ldots,\eta^n\}$ and its eigenvalues
are $\mathcal{E}(n)=1-q^n$.
This is a quite interesting situation that the similarity transformed
Hamiltonian $\widetilde{\mathcal{H}}$ is unbounded but the spectrum is
bounded; markedly different from those having the Jackson integral measures,
the big $q$-Jacobi family.
The corresponding polynomial solution of degree $n$ in $\eta$ is the
discrete $q$-Hermite polynomial $\II$ \cite{koeswart},
\begin{align}
  P_n(\eta)&\eqdef q^{\frac12n(n-1)}\tilde{h}_n(\eta;q)
  \eqdef i^{-n}{}_2\phi_0\Bigl(\genfrac{}{}{0pt}{}
  {q^{-n},\,i\eta}{-}\!\!\Bigm|\!q\,;-q^n\Bigr)\n
  &=q^{\frac12n(n-1)}\eta^n{}_2\phi_1\Bigl(\genfrac{}{}{0pt}{}
  {q^{-n},\,q^{1-n}}{0}\!\!\Bigm|\!q^2\,;-q^2\eta^{-2}\Bigr),\qquad
  P_n(-i)=(-i)^n,
  \label{dqHII:Pn}
\end{align}
which has definite parity, $P_n(-\eta)=(-1)^nP_n(\eta)$ due to the
evenness of the potential functions.
It is obtained from the Al-Salam-Carlitz $\II$ polynomial $V^{(a)}_n(\eta;q)$
by setting $a=-1$ and $\eta\to i\eta$,
$\tilde{h}_n(\eta;q)=i^{-n}V^{(-1)}_n(i\eta;q)$.
The difference equation, forward/backward shift relations and
three term recurrence relation are
\begin{align}
  &B^{\text{J}}(\eta)\bigl(P_n(\eta)-P_n(q\eta)\bigr)
  +D^{\text{J}}(\eta)\bigl(P_n(\eta)-P_n(q^{-1}\eta)\bigr)
  =\mathcal{E}(n)P_n(\eta),
  \label{dqHII:diffeq}\\
  &\eta^{-1}\bigl(P_n(\eta)-P_n(q\eta)\bigr)=\mathcal{E}(n)P_{n-1}(q\eta),
  \label{dqHII:forward}\\
  &\bigl(B^{\text{J}}(\eta)P_{n-1}(q\eta)
  -q^{-1}D^{\text{J}}(\eta)P_{n-1}(\eta)\bigr)\eta=P_n(\eta),
  \label{dqHII:backward}\\
  &\eta P_n(\eta)=A_nP_{n+1}(\eta)+C_nP_{n-1}(\eta),\quad
  A_n\eqdef q^{-n},\ \ C_n\eqdef q^{-n}-1.
  \label{dqHII:3term}
\end{align}

In the big $q$-Jacobi family, the two zeros of $D(x)$ determine the two
scales of $q^x$ appearing in the Jackson integral.
In the present case, the scale is arbitrary $\pm cq^x$, $c>0$, since
$D^{\text{J}}(\eta)$ has no zero and the potentials are even.
From now on we consider $c$ as the system parameter:
\begin{equation}
  q^{\bm{\lambda}}=c,\quad\bm{\delta}=1,\quad\kappa=q,\quad c>0.
  \label{dqHII:pararange}
\end{equation}
As shown for the big $q$-Jacobi family, in the single component formulation,
the self-adjointness is broken by the unboundedness of the potentials at
$x\to+\infty$ and the behaviour of the polynomial solutions $P_n(cq^x)$
at $x\to+\infty$. These will lead to the two component Hamiltonian
formulation with the sinusoidal coordinates
$\eta^{(\pm)}(x;\bm{\lambda})$:
\begin{equation}
  \eta^{(+)}(x;\bm{\lambda})\eqdef cq^x,\quad
  \eta^{(-)}(x;\bm{\lambda})\eqdef-cq^x,
  \label{eta_dqHII}
\end{equation}
and the corresponding potential functions \eqref{tilddef}
\begin{equation}
  B^{(\pm)}(x;\bm{\lambda})=1+c^{-2}q^{-2x}\eqdef B(x;\bm{\lambda})>0,\quad
  D^{(\pm)}(x;\bm{\lambda})=c^{-2}q^{1-2x}\eqdef D(x;\bm{\lambda})>0,
  \label{BD_dqHII}
\end{equation}
which are equal for $(\pm)$.
The eigenpolynomials are different for $(\pm)$:
\begin{equation}
  \check{P}^{(\pm)}_n(x;\bm{\lambda})\eqdef
  P_n\bigl(\eta^{(\pm)}(x;\bm{\lambda})\bigr)=
  i^{-n}{}_2\phi_0\Bigl(\genfrac{}{}{0pt}{}
  {q^{-n},\,\pm icq^x}{-}\!\Bigm|\!q\,;-q^n\Bigr).
\end{equation}
As in the b$q$J family, the appropriate ratio of the two zero mode vectors
$\phi^{(\pm)}_0(0;\bm{\lambda})$ achieve the self-adjointness of the
two component Hamiltonians, as demonstrated shortly.
Eq.\,\eqref{phi0inf} gives
\begin{align}
  \frac{\phi_0(x;\bm{\lambda})^2}{\phi_0(0;\bm{\lambda})^2}
  &=\left\{
  \begin{array}{lll}
  q^x(-c^2;q^2)_x&={\displaystyle
  \frac{q^x(-c^2;q^2)_{\infty}}{(-c^2q^{2x};q^2)_{\infty}}}
  &(x\in\mathbb{Z}_{\geq 0})\\[8pt]
  {\displaystyle
  \frac{c^{2x}q^{x^2}}{(-c^{-2}q^2;q^2)_{-x}}}&={\displaystyle
  \frac{c^{2x}q^{x^2}(-c^{-2}q^{2-2x};q^2)_{\infty}}
  {(-c^{-2}q^2;q^2)_{\infty}}}
  &(x\in\mathbb{Z}_{<0})
  \end{array}\right.,
  \label{phi02_dqHII_0}\\[6pt]
  &=\frac{q^x(-c^2;q^2)_{\infty}}{(-c^2q^{2x};q^2)_{\infty}}
  \ \ (x\in\mathbb{Z}).
  \label{phi02_dqHII_1}
\end{align}
In the last equality we have used an identity
\begin{equation*}
  1=\frac{(-c^2,-c^{-2}q^2;q^2)_{\infty}}
  {c^{2x}q^{x(x-1)}(-c^2q^{2x},-c^{-2}q^{2-2x};q^2)_{\infty}}
  \ \ (x\in\mathbb{Z}).
\end{equation*}

Now let us concentrate on the issue of self-adjointness.
Since it could be broken either at $x=+\infty$ or $x=-\infty$ or both,
the inner product of two single component vectors $f=(f(x))$ and $g=(g(x))$
is defined by
\begin{equation}
  (f,g)\eqdef\lim_{N,N'\to\infty}(f,g)_{N,N'},\quad
  (f,g)_{N,N'}\eqdef\sum_{x=-N'}^Nf(x)g(x).
 \label{inproNN'}
\end{equation}
As shown in \S\,\ref{formulation} \eqref{fHgN}--\eqref{fHgNlast}, for two
vectors $f(x)=\phi_0(x)\check{\mathcal{P}}(x)$ and
$g(x)=\phi_0(x)\check{\mathcal{Q}}(x)$, we have
\begin{align*}
  (f,\mathcal{H}g)_{N,N'}
  &=\sum_{x=-N'}^N\phi_0(x)^2\bigl(B(x)+D(x)\bigr)\check{\mathcal{P}}(x)
  \check{\mathcal{Q}}(x)
  -\sum_{x=-N'}^N\phi_0(x)^2B(x)\check{\mathcal{P}}(x)
  \check{\mathcal{Q}}(x+1)\n
  &\qquad
  -\!\!\sum_{x=-N'-1}^{N-1}\!\!\phi_0(x)^2B(x)\check{\mathcal{P}}(x+1)
  \check{\mathcal{Q}}(x),
\end{align*}
and obtain
\begin{align}
  &\quad(f,\mathcal{H}g)_{N,N'}-(\mathcal{H}f,g)_{N,N'}\n
  &=\phi_0(N)^2B(N)\bigl(\check{\mathcal{P}}(N+1)\check{\mathcal{Q}}(N)
  -\check{\mathcal{P}}(N)\check{\mathcal{Q}}(N+1)\bigr)\n
  &\quad+\phi_0(-N'-1)^2B(-N'-1)\bigl(\check{\mathcal{P}}(-N'-1)
  \check{\mathcal{Q}}(-N')
  -\check{\mathcal{P}}(-N')\check{\mathcal{Q}}(-N'-1)\bigr).
  \label{fHgNN'-HfgNN'}
\end{align}
At $x=-N'$, a degree $n$ polynomial grows like $\sim q^{-nN'}$, whereas the
potential damps as $B(-N')\sim q^{2N'}$ and the zero mode vector
$\phi_0(-N')^2$ provides much stronger damping as given in
\eqref{phi02_dqHII_0},
$\phi_0(-N';\bm{\lambda})^2\sim c^{-2N'}q^{N'^2}$.
Thus the second term in \eqref{fHgNN'-HfgNN'} vanishes as $N'\to\infty$ and 
the behaviour at $x=-\infty$ causes no problem for self-adjointness.
The behaviour at $x=+\infty$ is the same as that for the b$q$J family case
and we have to rely on the two component formulation. In this case the
two Hamiltonians $\mathcal{H}^{(\pm)}$ are exactly the same.
The eigenvectors of $\mathcal{H}^{(\pm)}$ have the same structures as those
of the big $q$ Jacobi \eqref{phine}.

The two component system has been introduced in \eqref{vec2}--\eqref{Ht2}.
The inner product of two vectors
$\bm{f}(x)=\genfrac{(}{)}{0pt}{}{f^{(+)}(x)}{f^{(-)}(x)}$ and
$\bm{g}(x)=\genfrac{(}{)}{0pt}{}{g^{(+)}(x)}{g^{(-)}(x)}$ is defined by
\begin{equation}
  \bin{\bm{f}}{\bm{g}}\eqdef\lim_{N\to\infty}
  {\bin{\bm{f}}{\bm{g}}}_N,\quad
  {\bin{\bm{f}}{\bm{g}}}_N\eqdef\sum_{x=-\infty}^N\bigl(
  f^{(+)}(x)g^{(+)}(x)+f^{(-)}(x)g^{(-)}(x)\bigr).
 \label{inproNN'2}
\end{equation}
For two vectors
$\bm{f}$ with $f^{(\pm)}(x)=\phi^{(\pm)}_0(x)\check{\mathcal{P}}^{(\pm)}(x)$
and
$\bm{g}$ with $g^{(\pm)}(x)=\phi^{(\pm)}_0(x)\check{\mathcal{Q}}^{(\pm)}(x)$,
eq.\,\eqref{fHgNN'-HfgNN'} gives
\begin{align}
  &\quad{\bin{\bm{f}}{\underline{\mathcal{H}}\,\bm{g}}}_N
  -{\bin{\underline{\mathcal{H}}\,\bm{f}}{\bm{g}}}_N\n
  &=\phi^{(+)}_0(N)^2B(N)\bigl(\check{\mathcal{P}}^{(+)}(N+1)
  \check{\mathcal{Q}}^{(+)}(N)
  -\check{\mathcal{P}}^{(+)}(N)\check{\mathcal{Q}}^{(+)}(N+1)\bigr)\n
  &\quad
  +\phi^{(-)}_0(N)^2B(N)\bigl(\check{\mathcal{P}}^{(-)}(N+1)
  \check{\mathcal{Q}}^{(-)}(N)
  -\check{\mathcal{P}}^{(-)}(N)\check{\mathcal{Q}}^{(-)}(N+1)\bigr).
  \label{fHgNN'-HfgNN'2}
\end{align}
The situation is simpler than that in the big $q$-Jacobi family case,
since $\phi^{(\pm)}_0(N)^2$ are only different by the unspecified
overall factor \eqref{phi02_dqHII_1}
\begin{equation*}
  \phi^{(\pm)}_0(N;\bm{\lambda})^2=\phi^{(\pm)}_0(0;\bm{\lambda})^2
  \frac{q^N(-c^2;q^2)_{\infty}}{(-c^2q^{2N};q^2)_{\infty}}.
\end{equation*}
It is very easy to see that for the equal choice
\begin{equation}
  \phi^{(\pm)}_0(0;\bm{\lambda})=\frac1{\sqrt{(-c^2;q^2)_\infty}}\
  \Rightarrow
  \phi^{(\pm)}_0(x;\bm{\lambda})^2=\phi_0(x;\bm{\lambda})^2
  =\frac{q^x}{(-c^2q^{2x};q^2)_\infty},
\end{equation}
the Hamiltonian is self-adjoint and the inner product of the two component
eigenvectors
\begin{equation}
  \bm{\phi}_n(x;\bm{\lambda})
  =\biggl(\begin{matrix}
  \phi^{(+)}_n(x;\bm{\lambda})\\\phi^{(-)}_n(x;\bm{\lambda})
  \end{matrix}\biggr)
  =\phi_0(x;\bm{\lambda})
  \begin{pmatrix}P_n(cq^x)\\P_n(-cq^x)\end{pmatrix},
  \label{dqHII:phipmn}
\end{equation}
is given by
\begin{align}
  \bin{\bm{\phi}_n}{\bm{\phi}_m}
  &=\bigl(\phi^{(+)}_n,\phi^{(+)}_m\bigr)
  +\bigl(\phi^{(-)}_n,\phi^{(-)}_m\bigr)
  =\frac{2\delta_{nm}}{d_n(\bm{\lambda})^2}\ \ (n,m=0,1,\ldots)
  \label{dqh2ortho}\\
  &=\sum_{x=-\infty}^{\infty}\bigl(
  \phi^{(+)}_n(x;\bm{\lambda})\phi^{(+)}_m(x;\bm{\lambda})
  +\phi^{(-)}_n(x;\bm{\lambda})\phi^{(-)}_m(x;\bm{\lambda})\bigr)\n
  &=\sum_{x=-\infty}^{\infty}\frac{q^x}{(-c^2q^{2x};q^2)_{\infty}}
  \bigl(P_n(cq^x)P_m(cq^x)+P_n(-cq^x)P_m(-cq^x)\bigr).
  \nonumber
\end{align}
Here the normalisation constant $d_n(\bm{\lambda})>0$ is
\begin{equation}
  d_n(\bm{\lambda})^2\eqdef
  \frac{q^n}{(q;q)_n}\times d_0(\bm{\lambda})^2,\quad
  d_0(\bm{\lambda})^2\eqdef
  \frac{(q,-c^2,-c^{-2}q^2;q^2)_{\infty}}{(q^2,-c^2q,-c^{-2}q;q^2)_{\infty}}.
  \label{dqhIIdn}
\end{equation}
Due to the definite parity $P_n(-\eta)=(-1)^nP_n(\eta)$, this orthogonal
relation is rewritten as
\begin{equation}
  \sum_{x=-\infty}^{\infty}\frac{1+(-1)^{n+m}}{2}
  \phi^{(+)}_n(x;\bm{\lambda})\phi^{(+)}_m(x;\bm{\lambda})
  =\frac{\delta_{nm}}{d_n(\bm{\lambda})^2}.
  \label{dqHII:(phipn,phipm)}
\end{equation}

The system is shape invariant, since the potential functions
$B^{(\pm)}(x;\bm{\lambda})$ and $D^{(\pm)}(x;\bm{\lambda})$
satisfy \eqref{shapeinv1cond1}--\eqref{shapeinv1cond2}.
The forward/backward shift relations for the d$q$H\,$\II$ have the form
\eqref{FBe}--\eqref{FBPne} with the auxiliary functions
$\varphi^{(\pm)}(x;\bm{\lambda})$,
\begin{equation}
  \varphi^{(\pm)}(x;\bm{\lambda})\eqdef
  \eta^{(\pm)}(x;\bm{\lambda}).
  \label{dqHII:varphipm}
\end{equation}
The triangularity of the potentials can be explained along the line of
\cite{os14}.
By taking $v_{k,l}$ as
\begin{equation*}
  \begin{array}{ll}
  v_{0,0}=q^{-1}(1-q)(1-q^2),&
  v_{1,0}=v_{0,1}=0,\\[2pt]
  v_{2,0}=q^{-1}-1+v_{0,2},&
  v_{1,1}=q-1-(q+q^{-1})v_{0,2},
  \end{array}
\end{equation*}
($v_{0,2}$ is arbitrary), the potential functions are expressed as
\eqref{os14B}--\eqref{os14D}.

\subsection{$q$-Laguerre}
\label{sec:qLag}

As is well known, the Laguerre polynomials of $\alpha=\mp\tfrac12$ are
related to the Hermite polynomials:
\begin{equation*}
  H_{2n}(x)=(-1)^n\,2^{2n}n!\,L^{(-\frac12)}_n(x^2),\quad
  H_{2n+1}(x)=(-1)^n\,2^{2n}n!\,2xL^{(\frac12)}_n(x^2).
\end{equation*}
The $q$-Laguerre ($q$L) polynomial $P_n(\eta;\bm{\lambda})$
($n\in\mathbb{Z}_{\geq0}$) \cite{koeswart},
\begin{align}
  P_n(\eta;\bm{\lambda})&\eqdef P_n(\eta;a;q),\quad P_n(-1;\bm{\lambda})=1,
  \label{qLPn}\\
  &\eqdef{}_2\phi_1\Bigl(\genfrac{}{}{0pt}{}
  {q^{-n},\,-\eta}{0}\!\!\Bigm|\!q\,;aq^{n+1}\Bigr)
  =(q;q)_nL^{(\alpha)}_n(\eta;q),\quad a=q^{\alpha},\n
  &=(-a)^nq^{n^2}\eta^n
  {}_2\phi_1\Bigl(\genfrac{}{}{0pt}{}
  {q^{-n},\,a^{-1}q^{-n}}{0}\!\!\Bigm|\!q\,;-q\eta^{-1}\Bigr),\n
  &=(aq;q)_n\,{}_1\phi_1\Bigl(\genfrac{}{}{0pt}{}
  {q^{-n}}{aq}\!\!\Bigm|\!q\,;-aq^{n+1}\eta\Bigr),\nonumber
\end{align}
is the counterpart of the discrete $q$-Hermite $\II$ polynomial defined on
the full integer lattice. For $\alpha=\mp\tfrac12$, they are related in a
similar way as the Hermite to the Laguerre:
\begin{equation}
  P^{\text{d$q$H\,$\II$}}_{2n}(\eta^{\frac12};q^{\frac12})
  =(-1)^nP_n(\eta;q^{-\frac12};q),\quad
  P^{\text{d$q$H\,$\II$}}_{2n+1}(\eta^{\frac12};q^{\frac12})
  =(-1)^n\eta^{\frac12}P_n(\eta;q^{\frac12};q).
\end{equation}
The $q$L is related to the Laguerre polynomial $L^{(\alpha)}_n(\eta)$
by the following limit,
\begin{equation}
  \lim_{q\to 1}(1-q)^{-n}P_n\bigl((1-q)\eta;q^{\alpha};q\bigr)
  =n!\,L^{(\alpha)}_n(\eta).
\end{equation}

The basic data of $q$L are
\begin{align}
  &q^{\bm{\lambda}}=(a,c),\ \ \bm{\delta}=(1,1),\ \ \kappa=q,
  \ \ 0<a<q^{-1},\ \ c>0,\\
  &\mathcal{E}(n)\eqdef 1-q^n,\\
  &B^{\text{J}}(\eta)\eqdef\eta^{-1}(\eta+1),\quad
  D^{\text{J}}(\eta;\bm{\lambda})\eqdef a^{-1}\eta^{-1},\\
  &A_n(\bm{\lambda})\eqdef-a^{-1}q^{-2n-1},\quad
  C_n(\bm{\lambda})\eqdef-a^{-1}q^{-2n}(1-q^n)(1-aq^n).
\end{align}
The $q$L polynomial $P_n(\eta;\bm{\lambda})$ is the degree
$n$ polynomial solution in $\eta$ of the second order difference equation
\begin{equation}
  B^{\text{J}}(\eta)\bigl(P_n(\eta;\bm{\lambda})-P_n(q\eta;\bm{\lambda})\bigr)
  +D^{\text{J}}(\eta;\bm{\lambda})\bigl(P_n(\eta;\bm{\lambda})
  -P_n(q^{-1}\eta;\bm{\lambda})\bigr)
  =\mathcal{E}(n)P_n(\eta;\bm{\lambda}).
  \label{qL:diffeq}
\end{equation}
Its recurrence relation and the forward/backward shift relations are
\begin{align} 
  &(1+\eta)P_n(\eta;\bm{\lambda})=A_n(\bm{\lambda})P_{n+1}(\eta;\bm{\lambda})
  -\bigl(A_n(\bm{\lambda})+C_n(\bm{\lambda})\bigr)P_n(\eta;\bm{\lambda})
  +C_n(\bm{\lambda})P_{n-1}(\eta;\bm{\lambda}),
  \label{qL3term}\\
  &(-aq\eta)^{-1}\bigl(P_n(\eta;\bm{\lambda})-P_n(q\eta;\bm{\lambda})\bigr)
  =\mathcal{E}(n)P_{n-1}(q\eta;\bm{\lambda}+\bm{\delta}),\\
  &\bigl(B^{\text{J}}(\eta)P_{n-1}(q\eta;\bm{\lambda}+\bm{\delta})
  -q^{-1}D^{\text{J}}(\eta;\bm{\lambda})
  P_{n-1}(\eta;\bm{\lambda}+\bm{\delta})\bigr)(-aq\eta)
  =P_n(\eta;\bm{\lambda}).
\end{align}

In terms of the sinusoidal coordinate
\begin{equation}
  \eta(x;\bm{\lambda})\eqdef cq^x,\quad
  \check{P}_n(x;\bm{\lambda})\eqdef
  P_n\bigl(\eta(x;\bm{\lambda});\bm{\lambda}\bigr),
\end{equation}
and the potential functions
\begin{equation}
  B(x;\bm{\lambda})\eqdef
  B^{\text{J}}\bigl(\eta(x;\bm{\lambda})\bigr)
  =1+c^{-1}q^{-x},\quad
  D(x;\bm{\lambda})\eqdef
  D^{\text{J}}\bigl(\eta(x;\bm{\lambda});\bm{\lambda}\bigr)
  =a^{-1}c^{-1}q^{-x},
  \label{qL:Dpm}
\end{equation}
the corresponding Hamiltonian $\mathcal{H}(\bm{\lambda})$ is defined on the
full integer lattice, $x\in\mathbb{Z}$.
Its ground state vector $\phi_0(x;\bm{\lambda})$ is determined in the same
way as for the d$q$H\,$\II$ \eqref{phi0inf},
\begin{equation}
  \phi_0(x;\bm{\lambda})^2=\frac{(aq)^x}{(-cq^x;q)_{\infty}},
  \label{qLphi0}
\end{equation}
in which we have adopted as the boundary condition
$\phi_0(0;\bm{\lambda})^2=(-c;q)_{\infty}^{-1}$.
Like the d$q$H\,$\II$ case, this measure function decreases very rapidly
at $x\to-\infty$ but the unboundedness of the potentials at large $N>0$ is
milder $\sim q^{-N}$ than that of the d$q$H\,$\II$. Thus the single
component Hamiltonian $\mathcal{H}(\bm{\lambda})$ is self-adjoint and its
eigenvectors are 
\begin{equation}
  \mathcal{H}(\bm{\lambda})\phi_n(x;\bm{\lambda})
  =\mathcal{E}(n)\phi_n(x;\bm{\lambda}),\quad 
  \phi_n(x;\bm{\lambda})
  =\phi_0(x;\bm{\lambda})\check{P}_n(x;\bm{\lambda})\ \ (n=0,1,\ldots).
  \label{qLeig}
\end{equation}
The orthogonality relation is the standard one
\begin{equation}
  (\phi_n,\phi_m)=\sum_{x=-\infty}^{\infty}
  \phi_n(x;\bm{\lambda})\phi_m(x;\bm{\lambda})
  =\frac{\delta_{nm}}{d_n(\bm{\lambda})^2}\ \ (n,m\in\mathbb{Z}_{\geq 0}),
\end{equation}
in which $d_n(\bm{\lambda})>0$ is given by
\begin{equation}
  d_n(\bm{\lambda})^2=\frac{q^n}{(q,aq;q)_n}\times d_0(\bm{\lambda})^2,\quad
  d_0(\bm{\lambda})^2=\frac{(aq,-c,-c^{-1}q;q)_{\infty}}
  {(q,-acq,-a^{-1}c^{-1};q)_{\infty}}.
\end{equation}
The system is shape invariant.
The forward/backward shift relations for the $q$L have the form
\eqref{FBe}--\eqref{FBPne} (without superscript $(\pm)$)
with the auxiliary functions $\varphi(x;\bm{\lambda})$,
\begin{equation}
  \varphi(x;\bm{\lambda})\eqdef-aq\eta(x;\bm{\lambda}).
  \label{qL:varphi}
\end{equation}

\subsection{Universal Rodrigues formula for polynomials with Jackson integrals}
\label{sec:uniRod}

The universal Rodrigues formula for the polynomials with Jackson integral
type measures has the same structure \eqref{univrod} as that for the rest
of the classical orthogonal polynomials of a discrete variable
(without superscript $(\pm)$ for $q$L):
\begin{align}
  P^{(\pm)}_n\bigl(\eta^{(\pm)}(x;\bm{\lambda});\bm{\lambda}\bigr)
  &=\phi^{(\pm)}_0(x;\bm{\lambda})^{-2}\,
  \bar{\mathcal{D}}^{(\pm)}(\bm{\lambda})
  \bar{\mathcal{D}}^{(\pm)}(\bm{\lambda}+\bm{\delta})\cdots
  \bar{\mathcal{D}}^{(\pm)}\bigl(\bm{\lambda}+(n-1)\bm{\delta}\bigr)\n
  &\quad\times
  \phi^{(\pm)}_0(x;\bm{\lambda}+n\bm{\delta})^2,
  \label{Junivrod}\\
  \bar{\mathcal{D}}^{(\pm)}(\bm{\lambda})&\eqdef
  (1-e^{-\partial})\bar{\varphi}^{(\pm)}(x;\bm{\lambda})^{-1},\quad
  \bar{\varphi}^{(\pm)}(x;\bm{\lambda})\propto 
  {\varphi}^{(\pm)}(x;\bm{\lambda}).
  \label{JDlamdef}
\end{align}
However, reflecting the different normalisation (boundary) conditions for
the ground state wavefunctions $\phi^{(\pm)}_0(x;\bm{\lambda})$ and for the
polynomials;
$P_n(1;\bm{\lambda})=1$ for b$q$J \eqref{Pn}, b$q$L \eqref{bqlpol2},
ASC\,$\I$ \eqref{ASCIpol}, $P_n(-i;\bm{\lambda})=(-i)^n$ for d$q$H\,$\II$
\eqref{dqHII:Pn}, $P_n(-1;\bm{\lambda})=1$ for $q$L \eqref{qLPn},
the function $\varphi(x;\bm{\lambda})$ in \eqref{Dlamdef} gets an extra
constant factor.
The modified functions $\bar{\varphi}^{(\pm)}(x;\bm{\lambda})$ are
\begin{align}
  \bar{\varphi}^{(\pm)}(x;\bm{\lambda})
  &=\eta^{(\pm)}(x;\bm{\lambda})\times\left\{
  \begin{array}{ll}
  (1-aq)(1-cq)(-acq)^{-1}&:\text{b$q$J}\\
  (1-aq)(1-bq)(-abq)^{-1} &:\text{b$q$L}\\
  1&:\text{ASC\,$\I$,\,d$q$H\,$\II$}
  \end{array}\right.,\n
  \bar{\varphi}(x;\bm{\lambda})
  &=\eta(x;\bm{\lambda})(-acq)^{-1}\ :\text{$q$L}.
\end{align}

\section{Other Topics}
\label{sec:other}

\subsection{Birth and Death processes related with Jackson integral measures}
\label{sec:bdp}

As an application of orthogonal polynomials of a discrete variable,
here we comment on birth and death processes \cite{ismail}.
A Birth and Death (BD) process is a typical stationary Markov process
with a one-dimensional discrete state space and the transitions occur
only between nearest neighbours.
The transition probability per unit time from $x$ to $x+1$ is $B(x)$
(birth rate) and $x$ to $x-1$ is $D(x)$ (death rate).
It was shown in \cite{bdproc} that classical orthogonal polynomials
obtained from hermitian matrices in I, {\em e.g.\/} the ($q$-)Racah and
(dual) ($q$-)Hahn polynomials, provide {\em exactly solvable birth and death
processes}.
Here we show that the orthogonal polynomials discussed in this paper,
{\em i.e.\/} those having Jackson integral measures, also supply examples
of exactly solvable birth and death processes with two component discrete
state spaces.
The dual polynomials of the big $q$-Jacobi family also give rise to
interesting examples of exactly solvable BD processes.
 
Let us summarise the essence of BD process within the framework of
`orthogonal polynomials from hermitian matrices' as explained in
\S\,\ref{sec:rdQM}.
The birth and death equation reads
\begin{equation}
  \frac{\partial}{\partial t}\mathcal{P}(x;t)
  =(L_{\text{BD}}\mathcal{P})(x;t),\quad
  \mathcal{P}(x;t)\ge0,\quad\sum_x\mathcal{P}(x;t)=1,
  \label{bdeqformal}
\end{equation}
in which $\mathcal{P}(x;t)$ is the probability distribution over a certain
discrete set of the parameter $x$.
The {\em birth and death operator\/} $L_{\text{BD}}$ is obtained by
{\em inverse similarity transformation\/} of the Hamiltonian $\mathcal{H}$:
\begin{align}
  &L_{\text{BD}}\eqdef-\phi_0\circ\mathcal{H}\circ\phi_0^{-1}
  =(e^{-\partial}-1)B(x)+(e^{\partial}-1)D(x),
  \label{LBDdef}\\
  &L_{\text{BD}}\phi_0(x)\phi_n(x)=-\mathcal{E}(n)\phi_0(x)\phi_n(x)
  \ \ (n=0,1,\ldots).
  \nonumber
\end{align}
Given an arbitrary initial probability distribution $\mathcal{P}(x;0)\ge0$
(with $\sum_x\mathcal{P}(x;0)=1$), the probability distribution at a
later time $t$ is
\begin{equation}
  \mathcal{P}(x;t)=\hat{\phi}_0(x)\sum_{n=0}^{\infty}
  c_n\,e^{-\mathcal{E}(n)t}\hat{\phi}_n(x)\ \ (t>0),
  \label{simpleBDsol}
\end{equation}
in which the constants $\{c_n\}$ are the expansion coefficients in terms
of the {\em complete set of normalised eigenfunctions\/}
$\{\hat{\phi}_n(x)\}$ of the Hamiltonian $\mathcal{H}$:
\begin{equation}
  \mathcal{P}(x;0)=\hat{\phi}_0(x)
  \sum_{n=0}^{\infty}c_n\hat{\phi}_n(x),
  \ \ c_n=\bigl(\hat{\phi}_n(x),\hat{\phi}_0(x)^{-1}
  \mathcal{P}(x;0)\bigr)
  \ \ (\Rightarrow c_0=1).
\end{equation}
For a concentrated initial distribution at $y$ ({\em e.g.\/}
$\mathcal{P}(x;0)=\delta_{xy}$), the transition probability from $y$ to
$x$ is given by
\begin{equation}
  \mathcal{P}(x,y;t)=\hat{\phi}_0(x)
  \Bigl(\sum_{n=0}^{\infty}e^{-\mathcal{E}(n)t}
  \hat{\phi}_n(x)\hat{\phi}_n(y)\Bigr)\hat{\phi}_0(y)^{-1}\ \ (t>0),
  \label{tranprobdisc1}
\end{equation}
which satisfies the so-called Chapman-Kolmogorov equation \cite{ismail}
\begin{equation}
  \mathcal{P}(x,y;t)=\sum_{z=0}^{\infty}\mathcal{P}(x,z;t-t')
  \mathcal{P}(z,y;t')\ \ (0<t'<t),
  \label{ChKol}
\end{equation}
as the consequence of the  normalised eigenfunctions
$\sum\limits_{z=0}^{\infty}\hat{\phi}_n(z)\hat{\phi}_m(z)=\delta_{nm}$.
The notation for the transition probability here is slightly changed
from \cite{bdproc}.

\subsubsection{two component BD}

Now let us formulate the two component BD processes matching with the big
$q$-Jacobi family \S\,\ref{sec:bqJ}.
Corresponding to the two component Hamiltonian \eqref{vec2}--\eqref{Ht2},
the state space has also two components,
\begin{equation}
  \bm{\mathcal{P}}(x;t)=
  \begin{pmatrix}
  \mathcal{P}^{(+)}(x;t)\\\mathcal{P}^{(-)}(x;t)
  \end{pmatrix},\quad
  \mathcal{P}^{(\pm)}(x;t)\ge0,\quad
  \sum_{x=0}^{\infty}\bigl(\mathcal{P}^{(+)}(x;t)
  +\mathcal{P}^{(-)}(x;t)\bigr)=1.
\end{equation}
The two component birth and death equation reads
\begin{align}
  &\frac{\partial}{\partial t}\bm{\mathcal{P}}(x;t)
  =\underline{L}_{\text{BD}}\bm{\mathcal{P}}(x;t),\\
  &\underline{L}_{\text{BD}}
  =-\underline{\phi_0}(x)\circ\underline{\mathcal{H}}\circ
  \underline{\phi_0}(x)^{-1}
  =(\underline{e^{-\partial}\!\!}\,-1)\underline{B}(x)
  +(\underline{e^{\partial}}-1)\underline{D}(x),
  \label{BD2}\\[2pt]
  &\underline{L}_{\text{BD}}
  \begin{pmatrix}
  \hat{\phi}^{(+)}_0(x)\hat{\phi}^{(+)}_n(x)\\[2pt]
  \hat{\phi}^{(-)}_0(x)\hat{\phi}^{(-)}_n(x)
  \end{pmatrix}=-\mathcal{E}(n)
  \begin{pmatrix}
  \hat{\phi}^{(+)}_0(x)\hat{\phi}^{(+)}_n(x)\\[2pt]
  \hat{\phi}^{(-)}_0(x)\hat{\phi}^{(-)}_n(x)
  \end{pmatrix}
  \ \ (n=0,1,\ldots).
\end{align}
For each type $(+)$ and $(-)$, the equation reads
\begin{equation}
  \frac{\partial}{\partial t}\mathcal{P}^{(\pm)}(x;t)=
  \bigl((e^{-\partial}-1)B^{(\pm)}(x)+(e^{\partial}-1)D^{(\pm)}(x)\bigr)
  \mathcal{P}^{(\pm)}(x;t),
  \label{BD3}
\end{equation}
in which $B^{(\pm)}(x)$ and $D^{(\pm)}(x)$ are given in
\eqref{bqJB+}--\eqref{bqJB-} for the b$q$J, \eqref{bqlB+}--\eqref{bqlB-}
for the b$q$L, and \eqref{ASCIB+}--\eqref{ASCIB-} for the ASC\,$\I$.
As seen from the equation, the direct interactions are between the nearest
neighbours of the same type only. The $(+)$ and $(-)$ sites are correlated
through the asymptotic points in such a way that the self-adjointness of
the corresponding Hamiltonian is realised.
For an arbitrary initial probability distribution
$\mathcal{P}^{(\pm)}(x;0)\ge0$
(with $\sum_{\epsilon}\sum_{x} \mathcal{P}^{(\epsilon)}(x;0)=1$),
the probability distribution at a later time $t$ is
\begin{equation}
  \bm{\mathcal{P}}(x;t)=\sum_{n=0}^{\infty}c_n\,e^{-\mathcal{E}(n)t}
  \begin{pmatrix}
  \hat{\phi}^{(+)}_0(x)\hat{\phi}^{(+)}_n(x)\\[2pt]
  \hat{\phi}^{(-)}_0(x)\hat{\phi}^{(-)}_n(x)
  \end{pmatrix}\ \ (t>0),
  \label{2comppt}
\end{equation}
in which the constants $\{c_n\}$ are the expansion coefficients in terms
of the complete set of normalised eigenvectors $\{\,\hat{\!\bm{\phi}}_n(x)\}$
of the Hamiltonian $\underline{\mathcal{H}}$:
\begin{align}
  &\bbin{\,\hat{\!\bm{\phi}}_n}{\,\hat{\!\bm{\phi}}_m}
  =\sum_{x=0}^{\infty}\bigl(\hat{\phi}^{(+)}_n(x)\hat{\phi}^{(+)}_m(x)+
  \hat{\phi}^{(-)}_n(x)\hat{\phi}^{(-)}_m(x)\bigr)=\delta_{nm},
  \label{2typortho}\\
  &c_n=\bbin{\,\hat{\!\bm{\phi}}_n(x)}
  {\underline{\hat{\phi}}_0(x)^{-1}\bm{\mathcal{P}}(x;0)}
  \ \ (\Rightarrow c_0=1)\n
  &\phantom{c_n}
  =\sum_{x=0}^{\infty}\bigl(\hat{\phi}^{(+)}_n(x)\hat{\phi}^{(+)}_0(x)^{-1}
  \mathcal{P}^{(+)}(x;0)+
  \hat{\phi}^{(-)}_n(x)\hat{\phi}^{(-)}_0(x)^{-1}\mathcal{P}^{(-)}(x;0)\bigr).
\end{align}
Asymptotically the above probability distribution \eqref{2comppt} approaches
to the stationary distribution
\begin{equation*}
  \lim_{t\to\infty}\bm{\mathcal{P}}(x;t)=
  \biggl(\begin{matrix}
  \hat{\phi}^{(+)}_0(x)^2\\[2pt]\hat{\phi}^{(-)}_0(x)^2
  \end{matrix}\biggr).
\end{equation*}

There are four types of transition probabilities, starting from $(+)$ or
$(-)$ component at site $y$ at $t=0$ and arriving at $(+)$ or $(-)$
component at site $x$ at a later time $t$. They form a $2\times2$ matrix
generalisation of the one component formula \eqref{tranprobdisc1},
\begin{equation}
  \bm{\mathcal{P}}(x,y;t)=\sum_{n=0}^{\infty}\,e^{-\mathcal{E}(n)t}
  \begin{pmatrix}
  \hat{\phi}^{(+)}_0(x)\hat{\phi}^{(+)}_n(x)\\[2pt]
  \hat{\phi}^{(-)}_0(x)\hat{\phi}^{(-)}_n(x)
  \end{pmatrix}
  \begin{pmatrix}
  \hat{\phi}^{(+)}_n(y)\hat{\phi}^{(+)}_0(y)^{-1}
  &\hspace*{-2pt}\hat{\phi}^{(-)}_n(y)\hat{\phi}^{(-)}_0(y)^{-1}
  \end{pmatrix},
\end{equation}
which also satisfies the matrix form Chapman-Kolmogorov equation
\begin{equation}
  \bm{\mathcal{P}}(x,y;t)=\sum_{z=0}^{\infty}\bm{\mathcal{P}}(x,z;t-t')
  \bm{\mathcal{P}}(z,y;t')\ \ (0<t'<t),
\end{equation}
as a consequence of the above orthogonality relation \eqref{2typortho}.

\subsubsection{discrete $q$-Hermite $\II$}

For the BD process corresponding to the discrete $q$-Hermite $\II$
\S\,\ref{sec:dqHII}, the above formulas require slight modifications.
The site $x$ is now on the full integer lattice and
$\hat{\phi}^{(+)}_0(x)=\hat{\phi}^{(-)}_0(x)=\hat{\phi}_0(x)$ and the
potentials $B^{(\pm)}(x)=B(x)$, $D^{(\pm)}(x)=D(x)$ are given in
\eqref{BD_dqHII}.
The transition probability is simply given by
\begin{equation}
  \bm{\mathcal{P}}(x,y;t)
  =\hat{\phi}_0(x)\sum_{n=0}^{\infty}\,e^{-\mathcal{E}(n)t}
  \begin{pmatrix}
  \hat{\phi}^{(+)}_n(x)\\[2pt]\hat{\phi}^{(-)}_n(x)
  \end{pmatrix}
  \begin{pmatrix}
  \hat{\phi}^{(+)}_n(y)&\hspace*{-2pt}\hat{\phi}^{(-)}_n(y)
  \end{pmatrix}\hat{\phi}_0(y)^{-1}.
\end{equation}

\subsubsection{$q$-Laguerre}

Except for the demographic interpretation, which is not applicable because of
the full integer lattice sites, all the formulas of the standard BD process
\eqref{bdeqformal}--\eqref{ChKol} are valid with the birth/death rates given
by \eqref{qL:Dpm} and the eigenvectors by \eqref{qLphi0}--\eqref{qLeig}.

\subsubsection{dual big $q$-Jacobi}

To write down various formulas for the BD processes corresponding to the
dual big $q$-Jacobi families \S\,\ref{sec:dbqJ} is an interesting exercise.
We fix one Hamiltonian, say $\mathcal{H}^{\text{d}\,(+)}$, with $A^{(+)}_n$
and $C^{(+)}_n$. The sites are now parametrised by $n\in\mathbb{Z}_{\ge0}$
and
\begin{equation}
  L_{\text{BD}}^{\text{d}}\eqdef-\phi^{\text{d}\,(+)}_0\circ
  \mathcal{H}^{\text{d}\,(+)}\circ{\phi^{\text{d}\,(+)}_0}^{-1}
  =-(e^{-\partial_n}-1)A^{(+)}_n-(e^{\partial_n}-1)C^{(+)}_n.
  \label{dualLBDdef}
\end{equation}
There are two series of eigenvectors
$\{\hat{\phi}^{\text{d}\,(+)}_x(n),\mathcal{E}^{\text{d}\,(+)}(x)\}$
\eqref{bqJ:Ed} and
$\{\hat{\phi}^{\text{d}\,(-)}_x(n),\mathcal{E}^{\prime\,\text{d}\,(+)}(x)\}$
\eqref{bqJ:E'd}.
The transition probability from site $m$ at time $t=0$ to site $n$ at a
later time $t$ is
\begin{align}
  \mathcal{P}^{\text{d}}(n,m;t)
  &=\hat{\phi}^{\text{d}\,(+)}_0(n)\sum_{x=0}^{\infty}\bigl(
  e^{-\mathcal{E}^{\text{d}\,(+)}(x)t}
  \hat{\phi}^{\text{d}\,(+)}_x(n)\hat{\phi}^{\text{d}\,(+)}_x(m)
  +e^{-\mathcal{E}^{\prime\,\text{d}\,(+)}(x)t}
  \hat{\phi}^{\text{d}\,(-)}_x(n)\hat{\phi}^{\text{d}\,(-)}_x(m)\bigr)\n
  &\quad\times\hat{\phi}^{\text{d}\,(+)}_0(m)^{-1}\ \ (t>0).
  \label{dualtranprobdisc1}
\end{align}
It also satisfies the Chapman-Kolmogorov equation \eqref{ChKol} thanks to
the orthogonality relation \eqref{Q-ortho4},
\begin{equation*}
  (\hat{\phi}^{\text{d}\,(\epsilon)}_x,\hat{\phi}^{\text{d}\,(\epsilon')}_y)
  =\sum_{n=0}^{\infty}\hat{\phi}^{\text{d}\,(\epsilon)}_x(n)
  \hat{\phi}^{\text{d}\,(\epsilon')}_y(n)
  =\delta_{\epsilon\,\epsilon'}\delta_{xy}.
\end{equation*}

\subsubsection{complete $q$-Meixner}

For the exact solvability of BD processes, the completeness of the
corresponding orthogonal polynomials is essential. Thus the simple form
of solutions \eqref{simpleBDsol}, \eqref{tranprobdisc1}, given in
\cite{bdproc} for the $q$-Meixner and $q$-Charlier is flawed.
The form of the equation is the same as given in \cite{bdproc},
{\em i.e.\/} \eqref{bdeqformal} with the birth/death rates given by
\eqref{qMB&D} or \eqref{qCB&D}. The correct transition probability from
site $y$ at time $t=0$ to site $x$ at a later time $t$ has a similar form
to that of the dual big $q$-Jacobi \eqref{dualtranprobdisc1},
\begin{equation}
  \mathcal{P}(x,y;t)=\hat{\phi}_0(x)
  \sum_{n=0}^{\infty}\Bigl(
  e^{-\mathcal{E}(n)t}\hat{\phi}_n(x)\hat{\phi}_n(y)
  +e^{-\mathcal{E}^{\prime}(n)t}\hat{\phi}^{(-)}_n(x)\hat{\phi}^{(-)}_n(y)
  \Bigr)\hat{\phi}_0(y)^{-1}\ \ (t>0),
  \label{qMdualtranprob}
\end{equation}
in which the supplementary vectors $\{\hat{\phi}^{(-)}_n(x)\}$ contribute.
Here $\mathcal{E}(n)=1-q^n$ and $\mathcal{E}^{\prime}(n)=1+cq^n$
\eqref{qMsupeig} for $q$M and $\mathcal{E}^{\prime}(n)=1+aq^n$
\eqref{qCsupeig} for $q$C. The eigenvectors are
\eqref{qMori}--\eqref{qMsup} for $q$M and \eqref{qCsup0}--\eqref{qCsup}
for $q$C. The transition probability \eqref{qMdualtranprob} satisfies the
Chapman-Kolmogorov equation \eqref{ChKol}.

\subsection{Proposal for new normalisation for some orthogonal polynomials}
\label{sec:nnorm}

In general, the normalisation of orthogonal polynomials is a matter of
historical conventions. However, as shown in \cite{os12,os14} and
recapitulated in \S\,\ref{polys}, there is an unambiguous and universal rule
for normalisation for orthogonal polynomials of a discrete variable defined
on a finite integer lattice or on a semi-infinite integer lattice
\eqref{Pzero},
\begin{equation*}
  \check{P}_n(0;\bm{\lambda})
  =P_n\bigl(\eta(0;\bm{\lambda});\bm{\lambda}\bigr)=P_n(0;\bm{\lambda})=1
  \ \ (n=0,1,2,\ldots).
\end{equation*}
Most of the conventional normalisation follow the above rule,
{\em except for \/} the little $q$-Jacobi, little $q$-Laguerre,
Al-Salam-Carlitz $\II$ and the alternative $q$-Charlier ($q$-Bessel)
polynomials. Below we present the explicit hypergeometric expressions for
the above four polynomials satisfying the universal normalisation rule
\eqref{Pzero}. They all have $q^{-n}$ and $q^{-x}$ among the upper indices
$\{a_1,\ldots,a_r\}$ of the $q$-hypergeometric function
${}_r\phi_s\Bigl(\genfrac{}{}{0pt}{}{a_1,\,\ldots,\,a_r}
{b_1,\,\ldots,\,b_s}\!\Bigm|\!q\,;z\Bigr)$.
In I, the universal normalisation rules were applied to these four polynomials
by multiplying appropriate rescaling factors to the conventional definitions,
but the explicit hypergeometric expressions were not reported.

\subsubsection{little $q$-Jacobi}

The conventional form \cite{koeswart} is
\begin{equation*}
  p^{\text{conv}}_n(q^x;a,b|q)={}_2\phi_1\Bigl(
  \genfrac{}{}{0pt}{}{q^{-n},\,abq^{n+1}}{aq}\!\Bigm|\!q\,;q^{x+1}\Bigr).
\end{equation*}
The proposed new form is
\begin{align}
  P_n\bigl(\eta(x);\bm{\lambda}\bigr)&={}_3\phi_1\Bigl(
  \genfrac{}{}{0pt}{}{q^{-n},\,abq^{n+1},\,q^{-x}}{bq}\!\Bigm|\!q\,;
  \frac{q^{x}}{a}\Bigr)
  =(-a)^{-n}q^{-\frac12n(n+1)}\frac{(aq;q)_n}{(bq;q)_n}
  p^{\text{conv}}_n(q^x;a,b|q),\n[2pt]
  \eta(x)&=1-q^x,\quad\mathcal{E}(n;\bm{\lambda})=(q^{-n}-1)(1-abq^{n+1}).
\end{align}
The new form is obtained by using transformation formula (III.8) in
\cite{gasper}.

\subsubsection{little $q$-Laguerre}

The conventional form \cite{koeswart} is
\begin{equation*}
  p^{\text{conv}}_n(q^x;a|q)={}_2\phi_1\Bigl(
  \genfrac{}{}{0pt}{}{q^{-n},\,0}{aq}\!\Bigm|\!q\,;q^{x+1}\Bigr).
\end{equation*}
The universal form \cite{os12} is
\begin{align}
  P_n\bigl(\eta(x);\bm{\lambda}\bigr)&={}_2\phi_0\Bigl(
  \genfrac{}{}{0pt}{}{q^{-n},\,q^{-x}}{-}\!\Bigm|\!q\,;\frac{q^x}{a}\Bigr)
  =(a^{-1}q^{-n};q)_n\,p^{\text{conv}}_n(q^x;a|q),
  \label{lqL}\n[2pt]
  \eta(x)&=1-q^x,\quad\mathcal{E}(n)=q^{-n}-1.
\end{align}

\subsubsection{Al-Salam-Carlitz $\II$}

The conventional form \cite{koeswart} is
\begin{equation*}
  V^{(a)}_n(q^{-x};q)=(-a)^nq^{-\frac12n(n-1)}{}_2\phi_0\Bigl(
  \genfrac{}{}{0pt}{}{q^{-n},\,q^{-x}}{-}\!\Bigm|\!q\,;\frac{q^{n}}{a}\Bigr).
\end{equation*}
The universal form \cite{os12} is
\begin{align}
  P_n\bigl(\eta(x);\bm{\lambda}\bigr)&={}_2\phi_0\Bigl(
  \genfrac{}{}{0pt}{}{q^{-n},\,q^{-x}}{-}\!\Bigm|\!q\,;\frac{q^n}{a}\Bigr)
  =(-a)^{-n}q^{\frac12n(n-1)}\,V^{(a)}_n(q^{-x};q),
  \label{ASC2}\n[2pt]
  \eta(x)&=q^{-x}-1,\quad\mathcal{E}(n)=1-q^n.
\end{align}
It is obvious that the little $q$-Laguerre polynomial \eqref{lqL} and the
Al-Salam-Carlitz $\II$ polynomial \eqref{ASC2} are dual to each other,
that is, they are interchanged by $x\leftrightarrow n$.

\subsubsection{alternative $q$-Charlier ($q$-Bessel)}

The conventional form \cite{koeswart} is
\begin{equation*}
  K_n(q^x;a;q)={}_2\phi_1\Bigl(
  \genfrac{}{}{0pt}{}{q^{-n},\,-aq^n}{0}\!\Bigm|\!q\,;q^{x+1}\Bigr).
\end{equation*}
The proposed new form is
\begin{align}
  P_n\bigl(\eta(x);\bm{\lambda}\bigr)&={}_3\phi_0\Bigl(
  \genfrac{}{}{0pt}{}{q^{-n},\,-aq^n,\,q^{-x}}{-}\!\Bigm|\!q\,;
  -\frac{q^x}{a}\Bigr)
  =q^{nx}{}_2\phi_1\Bigl(
  \genfrac{}{}{0pt}{}{q^{-n},\,q^{-x}}{0}\!\Bigm|\!q\,;
  -\frac{q^{1-n}}{a}\Bigr)\n
  &=(-aq^n)^{-n}\,K_n(q^x;a;q),\n[2pt]
  \eta(x)&=1-q^x,\quad\mathcal{E}(n;\bm{\lambda})=(q^{-n}-1)(1+aq^n).
\end{align}
The new form is obtained by using transformation formula (III.8)
(with $c\to 0$) in \cite{gasper}.

Let us stress, however, the universal normalisation cannot be applied to
those having Jackson integral measures. They have two component
eigenpolynomials $P^{(\pm)}_n$ and their fixed ratios due to the
self-adjointness do not allow $n$-dependent rescaling for each of
$P^{(\pm)}_n$.

\section{Summary and Comments}
\label{sec:summary}

A brief summary of the project ``orthogonal polynomials from hermitian
matrices'' is given in section two. In essence the project offers unified
understanding of various properties of orthogonal polynomials with purely
discrete orthogonality measures through the framework of eigenvalue problems
of certain Jacobi matrices, to be called the Hamiltonians.
It is pointed out that for the orthogonal polynomials defined on the
semi-infinite integer lattice, the completeness relation or the dual
orthogonality relation does not necessarily hold.
In section three the failure of the orthogonality relation for the proposed
dual $q$-Meixner and dual $q$-Charlier polynomials in \cite{os12} is
attributed to the breakdown of the self-adjointness due to the asymptotic
behaviours of the potential functions and the polynomials.
This also means that the $q$-Meixner ($q$-Charlier) polynomials do not
form a complete set of orthogonal vectors in the $\ell^2$ Hilbert space.
The self-adjointness is recovered by the introduction of another Hamiltonian
system which is obtained by rescaling the sinusoidal coordinate.
In section four, the two component Hamiltonian formalism with two types
of polynomials and the corresponding potentials is applied to the
polynomials belonging to the big $q$-Jacobi polynomials family.
It is shown that the recovery of the self-adjointness leads naturally to
the Jackson integral measures for these polynomials.
The striking features of the difference equations governing the families
of dual big $q$-Jacobi polynomials are also explored in section four.
The infinite dimensional Jacobi matrix (Hamiltonian) corresponding to the
dual difference equation, or the three term recurrence relation of the
rescaled original polynomial, has a spectrum consisting of two infinite
series sharing one accumulation point.
In \S\,\ref{sec:compqm} the complete set of orthogonal vectors involving
the $q$-Meixner ($q$-Charlier) polynomials is presented.
In section five, the orthogonality of the discrete $q$-Hermite $\II$
polynomial is explained within the framework of the extended two component
Hamiltonian formalism.
Since the system is defined on the full integer lattice $x\in\mathbb{Z}$,
the corresponding dual polynomials cannot be defined by the interchange
$\mathbb{Z}_{\ge0}\ni n\leftrightarrow x\in\mathbb{Z}$.
As another example of orthogonal polynomials defined on the full integer
lattice, the $q$-Laguerre polynomial and its properties are explored in
\S\ref{sec:qLag}.
The universal Rodrigues formulas for those having Jackson integral measures
are presented in \S\ref{sec:uniRod}.
The birth and death processes corresponding to the polynomials having
Jackson integral measures, {\em i.e.\/} the big $q$-Jacobi family and
the discrete $q$-Hermite $\II$ are explored in \S\,\ref{sec:bdp}.
The corresponding stationary stochastic process has two states $(\pm)$
at each site interacting with the nearest neighbours of the same type only.
In \S\,\ref{sec:nnorm} we propose that the universal normalisation condition
should be applied to four classical orthogonal polynomials, the little
$q$-Jacobi, little $q$-Laguerre, Al-Salam-Carlitz $\II$ and the alternative
$q$-Charlier ($q$-Bessel) polynomials. Their explicit forms are presented.
By adopting the universal normalisation condition \eqref{Pzero} with the
corresponding simple $q$-hypergeometric expressions, the identification of
the dual polynomial is automatic by the interchange $x\leftrightarrow n$.
In Appendix, an alternative solution method based on closure relations is
applied to those polynomials discussed in the main text. The Heisenberg
operator solutions for the sinusoidal coordinates, the creation/annihilation
operators, the self-consistent equations determining the eigenvalues, etc.
are derived.

The self-adjointness of the systems having Jackson integral measures has
been discussed in \cite{koeswart2} from a different point of view.

In these days infinitely many examples of orthogonal polynomials satisfying
second order differential/difference equations are known
\cite{krein}--\cite{os25}.
Since they have `holes' in the degrees or the lowest degree is greater
than 1, they do not satisfy the three term recurrence relations.
They are not classical.
It is an interesting challenge to try and construct the multi-indexed
deformations of various classical orthogonal polynomials defined on the
semi-infinite and full infinite integer lattices.
As shown for the other classical orthogonal polynomials
\cite{os25,os23,os26,os27}, the method (multiple Darboux transformations)
is well established. One only needs to find appropriate seed solutions.

\section*{Acknowledgements}

S.\,O. is supported in part by Grant-in-Aid for Scientific Research
from the Ministry of Education, Culture, Sports, Science and Technology
(MEXT), No.25400395.

\addcontentsline{toc}{section}{Appendix: Closure Relation}
\section*{Appendix: Closure Relation}
\label{append}
\setcounter{equation}{0}
\setcounter{section}{1}
\renewcommand{\theequation}{A.\arabic{equation}}
\renewcommand{\thesection}{\Alph{section}}

Here we discuss another symmetry of classical orthogonal polynomials
called {\em closure relation\/}.
While shape invariance leads to exact solvability in the Schr\"odinger
picture, closure relation means exact solvability in the Heisenberg picture.
That is, the Heisenberg operator solution of the sinusoidal coordinate
$e^{i\mathcal{H}t}\eta(x)e^{-i\mathcal{H}t}$ can be obtained explicitly
and its positive/negative frequency parts provide the annihilation/creation
operators, which enable algebraic determination of all eigenvalues and
eigenvectors.
It has been fully discussed for the classical orthogonal polynomials of
a discrete variable in I, except for
\romannumeral1) those having Jackson integral measures, {\em i.e.\/}
the big $q$-Jacobi family,
\romannumeral2) duals of the big $q$-Jacobi family,
\romannumeral2') the $q$-Meixner and $q$-Charlier,
\romannumeral3) the discrete $q$-Hermite $\II$,
and \romannumeral3') the $q$-Laguerre.

The closure relation is the commutation relation between the Hamiltonian
$\mathcal{H}$ and the sinusoidal coordinate $\eta=\eta(x)$,
\begin{align}
  &\bigl[\mathcal{H},[\mathcal{H},\eta]\bigr]
  =\eta R_0(\mathcal{H})
  +[\mathcal{H},\eta]R_1(\mathcal{H})+R_{-1}(\mathcal{H}),
  \label{closurerel1}\\
  \text{or}\quad&
  \bigl[\widetilde{\mathcal{H}},[\widetilde{\mathcal{H}},\eta]\bigr]
  =\eta R_0(\widetilde{\mathcal{H}})
  +[\widetilde{\mathcal{H}},\eta]R_1(\widetilde{\mathcal{H}})
  +R_{-1}(\widetilde{\mathcal{H}}),
  \label{closurerel1t}
\end{align}
in which $R_i(z)$ are polynomials in $z$,
\begin{equation}
  R_1(z)=r_1^{(1)}z+r_1^{(0)},
  \ \ R_0(z)=r_0^{(2)}z^2+r_0^{(1)}z+r_0^{(0)},
  \ \ R_{-1}(z)=r_{-1}^{(2)}z^2+r_{-1}^{(1)}z+r_{-1}^{(0)},
  \label{Ricoeff}
\end{equation}
and $r^{(j)}_k$ are real constants.
The structure of the closure relation is unchanged by 
an affine transformation of $\eta$,
$\eta^{\text{new}}(x)=a\eta(x)+b$ ($a,b$: constant, $a\neq 0$), 
\begin{align*}
  &\bigl[\mathcal{H},[\mathcal{H},\eta^{\text{new}}]\bigr]
  =\eta^{\text{new}}R^{\text{new}}_0(\mathcal{H})
  +[\mathcal{H},\eta^{\text{new}}]R^{\text{new}}_1(\mathcal{H})
  +R^{\text{new}}_{-1}(\mathcal{H}),\\
  &R^{\text{new}}_1(z)=R_1(z),\quad
  R^{\text{new}}_0(z)=R_0(z),\quad
  R^{\text{new}}_{-1}(z)=aR_{-1}(z)-bR_0(z).
\end{align*}
The Heisenberg operator solution for $\eta(x)$ and the creation/annihilation
operator $a^{(\pm)}$ are given by \cite{os12,os7},
\begin{align}
  &e^{it\mathcal{H}}\eta(x)e^{-it\mathcal{H}}
  =a^{(+)}e^{i\alpha_+(\mathcal{H})t}+a^{(-)}e^{i\alpha_-(\mathcal{H})t}
  -R_{-1}(\mathcal{H})R_0(\mathcal{H})^{-1},
  \label{Heisensol}\\
  &\alpha_{\pm}(z)\eqdef\tfrac12\bigl(R_1(z)
  \pm\sqrt{R_1(z)^2+4R_0(z)}\,\bigr),
  \label{alphapmexp}\\
  &\qquad\qquad
  R_1(z)=\alpha_+(z)+\alpha_-(z),\quad
  R_0(z)=-\alpha_+(z)\alpha_-(z),
  \label{R1R0}\\
  &a^{(\pm)}\eqdef\pm\Bigl([\mathcal{H},\eta(x)]
  -\bigl(\eta(x)+R_{-1}(\mathcal{H})R_0(\mathcal{H})^{-1}\bigr)
  \alpha_{\mp}(\mathcal{H})\Bigr)
  \bigl(\alpha_+(\mathcal{H})-\alpha_-(\mathcal{H})\bigr)^{-1}\n
  &\phantom{a^{(\pm)}}=
  \pm\bigl(\alpha_+(\mathcal{H})-\alpha_-(\mathcal{H})\bigr)^{-1}
  \Bigl([\mathcal{H},\eta(x)]
  +\alpha_{\pm}(\mathcal{H})\bigl(\eta(x)
  +R_{-1}(\mathcal{H})R_0(\mathcal{H})^{-1}\bigr)\Bigr).
  \label{a^{(pm)}}
\end{align}
The necessary and sufficient condition for \eqref{closurerel1t} is
(I.4.62)--(I.4.66), in which $\eta(0)=0$ is not assumed. From these conditions
the following relations are derived (I.4.69), (I.4.70),
\begin{align}
  &r_0^{(2)}=r_1^{(1)},\quad r_0^{(1)}=2r_1^{(0)},
  \label{cransr}\\
  &\eta(x+2)-(2+r_1^{(1)})\eta(x+1)+\eta(x)=r_{-1}^{(2)}.
  \label{crcond1pp}
\end{align}
The closure relation is intimately related to the three term recurrence
relation. For the orthogonal polynomials $P_n(\eta)$ in I ($P_n(0)=1$,
$\eta(0)=0$), the three term recurrence relation has the following form:
\begin{equation}
  \eta P_n(\eta)=A_nP_{n+1}(\eta)-(A_n+C_n)P_n(\eta)+C_nP_{n-1}(\eta).
  \label{3term_P}
\end{equation}
Since eigenvectors have the form $\phi_n(x)=\phi_0(x)P_n(\eta(x))$,
this gives
\begin{align}
  &\alpha_{\pm}\bigl(\mathcal{E}(n)\bigr)
  =\mathcal{E}(n\pm 1)-\mathcal{E}(n),
  \label{alphapmE}\\
  &R_{-1}\bigl(\mathcal{E}(n)\bigr)R_0\bigl(\mathcal{E}(n)\bigr)^{-1}
  =A_n+C_n,
  \label{Rm1E}\\
  &a^{(+)}\phi_n(x)=A_n\phi_{n+1}(x),\quad
  a^{(-)}\phi_n(x)=C_n\phi_{n-1}(x).
  \label{apmphin}
\end{align}

\subsection{Big $q$-Jacobi family}

The polynomials $P_n(\eta)$,
\eqref{Pn} for b$q$J, \eqref{bqlpol2} for b$q$L and
\eqref{ASCIpol} for ASC\,$\I$ ($a=-1$: d$q$H\,$\I$),
satisfy $P_n(1)=1$ and the three term recurrence relation has the
following form \eqref{bqj3}
\begin{equation}
  (1-\eta)P_n(\eta)=A_nP_{n+1}(\eta)-(A_n+C_n)P_n(\eta)+C_nP_{n-1}(\eta),
  \label{3term_P_i)}
\end{equation}
where $A_n$ and $C_n$ are given in
\eqref{bqjAn}--\eqref{bqjCn}, \eqref{bqlACn} and \eqref{ASCI:An,Cn}.
{}From this form, we take the sinusoidal coordinate for the closure relation
\eqref{closurerel1} as $\eta(x)=1-\eta^{(\pm)}(x)$.
Since $B^{(\pm)}(x)$ and $D^{(\pm)}(x)$ are obtained from $B^{\text{J}}(x)$
and $D^{\text{J}}(x)$ as \eqref{tilddef}, the two Hamiltonians
$\mathcal{H}^{(\pm)}$ satisfy the closure relation \eqref{closurerel1}
with the same $R_i(z)$.
We can verify \eqref{alphapmE}, \eqref{Rm1E} and
\begin{equation}
  a^{(+)}\phi^{(\pm)}_n(x)=A_n\phi^{(\pm)}_{n+1}(x),\quad
  a^{(-)}\phi^{(\pm)}_n(x)=C_n\phi^{(\pm)}_{n-1}(x).
  \label{apmphin_i)}
\end{equation}
It should be emphasised that when $z$ is replaced by the actual spectrum
$\mathcal{E}(n)$ ($\mathcal{E}(n)=(q^{-n}-1)(1-abq^{n+1})$ for b$q$J,
$\mathcal{E}(n)=q^{-n}-1$ for the rest), $R_1(z)^2+4R_0(z)$ in
\eqref{alphapmexp} becomes a complete square, as in all the other cases
discussed in I.

The data of $R_i(x)$ are\\
\underline{big $q$-Jacobi} :
\begin{align}
  R_1(z)&=(q^{\frac12}-q^{-\frac12})^2(z+1+abq),\n
  R_0(z)&=(q^{\frac12}-q^{-\frac12})^2
  \bigl((z+1+abq)^2-(1+q)^2ab\bigr),\\
  R_{-1}(z)&=-(q^{\frac12}-q^{-\frac12})^2
  \Bigl(z^2+\bigl(2-(a+c-ab+ac)q\bigr)z+(1-ab)(1-aq)(1-cq)\Bigr),\nonumber
\end{align}
\underline{big $q$-Laguerre} :
\begin{align}
  R_1(z)&=(q^{\frac12}-q^{-\frac12})^2(z+1),\n
  R_0(z)&=(q^{\frac12}-q^{-\frac12})^2(z+1)^2
  \ \ \Bigl(\Rightarrow\sqrt{R_1(z)^2+4R_0(z)}=(q^{-1}-q)|1+z|\Bigr),\\
  R_{-1}(z)&=-(q^{\frac12}-q^{-\frac12})^2
  \Bigl(z^2+\bigl(2-(a+b+ab)q\bigr)z+(1-aq)(1-bq)\Bigr),\nonumber
\end{align}
\underline{Al-Salam-Carlitz I} :
\begin{align}
  R_1(z)&=(q^{\frac12}-q^{-\frac12})^2(z+1),\n
  R_0(z)&=(q^{\frac12}-q^{-\frac12})^2(z+1)^2
  \ \ \Bigl(\Rightarrow\sqrt{R_1(z)^2+4R_0(z)}=(q^{-1}-q)|1+z|\Bigr),\\
  R_{-1}(z)&=-(q^{\frac12}-q^{-\frac12})^2
  (z+1)(z-a).\nonumber
\end{align}

\subsection{Dual big $q$-Jacobi family}

The polynomials $Q^{(\pm)}_x(\mathcal{E})$
\eqref{Qx+}--\eqref{Qx-} for b$q$J, \eqref{bqlQx+}--\eqref{bqlQx-} for b$q$L
and \eqref{ASCIQx+}--\eqref{ASCIQx-} for ASC\,$\I$ ($a=-1$: d$q$H\,$\I$),
satisfy $Q^{(\pm)}_x(0)=1$ and the three term recurrence relation has the
following form \eqref{3termQxe}
\begin{equation}
  \mathcal{E}Q^{(\pm)}_x(\mathcal{E})
  =-B^{(\pm)}(x)Q^{(\pm)}_{x+1}(\mathcal{E})
  +\bigl(B^{(\pm)}(x)+D^{(\pm)}(x)\bigr)Q^{(\pm)}_x(\mathcal{E})
  -D^{(\pm)}(x)Q^{(\pm)}_{x-1}(\mathcal{E}),
  \label{3term_P_ii)}
\end{equation}
where $B^{(\pm)}(x)$ and $D^{(\pm)}(x)$ are given in
\eqref{bqJB+}--\eqref{bqJB-}, \eqref{bqlB+}--\eqref{bqlB-} and
\eqref{ASCIB+}--\eqref{ASCIB-}.
{}From this form, we take the sinusoidal coordinate as
$\mathcal{E}(n)$ satisfying $\mathcal{E}(0)=0$,
$\mathcal{E}(n)=(q^{-n}-1)(1-abq^{n+1})$ for b$q$J and
$\mathcal{E}(n)=q^{-n}-1$ for the rest.
The closure relation for $\mathcal{H}^{\text{d}\,(+)}$ reads
\begin{equation}
  \bigl[\mathcal{H}^{\text{d}\,(+)},
  [\mathcal{H}^{\text{d}\,(+)},\mathcal{E}(n)]\bigr]
  =\mathcal{E}(n)R_0\bigl(\mathcal{H}^{\text{d}\,(+)}\bigr)
  +[\mathcal{H}^{\text{d}\,(+)},\mathcal{E}(n)]
  R_1\bigl(\mathcal{H}^{\text{d}\,(+)}\bigr)
  +R_{-1}\bigl(\mathcal{H}^{\text{d}\,(+)}\bigr).
  \label{dualbqJ:closurerel1}
\end{equation}

The system has two types of the eigenvectors, the ordinary
$\{\phi^{\text{d}\,(+)}_x(n)\}$ and
the supplementary $\{\phi^{\text{d}\,(-)}_x(n)\}$,
\eqref{bqJ:phidpm}, \eqref{bqL:phidpm} and \eqref{ASCI:phidpm},
with the corresponding eigenvalues
$\{\mathcal{E}^{\text{d}\,(+)}(x)\}$ and
$\{\mathcal{E}^{\prime\,\text{d}\,(+)}(x)\}$ 
\eqref{bqJ:Edpm}, \eqref{bqJ:E'd}, \eqref{bqL:E'd} and \eqref{ASCI:E'd},
respectively.
The latter is a monotonously decreasing function of $x$ and it lies above the
former $\mathcal{E}^{\prime\,\text{d}\,(+)}(x')>\mathcal{E}^{\text{d}\,(+)}(x)$.
We can verify that corresponding to \eqref{alphapmE}--\eqref{apmphin},
the ordinary sector is controlled by the data of the closure relation as
follows:
\begin{align}
  &\alpha_{\pm}\bigl(\mathcal{E}^{\text{d}\,(+)}(x)\bigr)
  =\mathcal{E}^{\text{d}\,(+)}(x\pm 1)-\mathcal{E}^{\text{d}\,(+)}(x),
  \label{dual:alphapmE}\\
  &R_{-1}\bigl(\mathcal{E}^{\text{d}\,(+)}(x)\bigr)
  R_0\bigl(\mathcal{E}^{\text{d}\,(+)}(x)\bigr)^{-1}
  =-B^{(+)}(x)-D^{(+)}(x),
  \label{dual:Rm1E}\\
  &a^{(+)}\phi^{\text{d}\,(+)}_x(n)
  =-B^{(+)}(x)\phi^{\text{d}\,(+)}_{x+1}(n),\quad
  a^{(-)}\phi^{\text{d}\,(+)}_x(n)
  =-D^{(+)}(x)\phi^{\text{d}\,(+)}_{x-1}(n),
  \label{dual:apmphin}
\end{align}
whereas for the supplementary sector, the order of $x$ is reversed:
\begin{align}
  &\alpha_{\pm}\bigl(\mathcal{E}^{\prime\,\text{d}\,(+)}(x)\bigr)
  =\mathcal{E}^{\prime\,\text{d}\,(+)}(x\mp 1)
  -\mathcal{E}^{\prime\,\text{d}\,(+)}(x),
  \label{dual:alphapmE2}\\
  &R_{-1}\bigl(\mathcal{E}^{\prime\,\text{d}\,(+)}(x)\bigr)
  R_0\bigl(\mathcal{E}^{\prime\,\text{d}\,(+)}(x)\bigr)^{-1}
  =-B^{(-)}(x)-D^{(-)}(x),
  \label{dual:Rm1E2}\\
  &a^{(+)}\phi^{\text{d}\,(-)}_x(n)
  =-D^{(-)}(x)\phi^{\text{d}\,(-)}_{x-1}(n),
  \ \ a^{(-)}\phi^{\text{d}\,(-)}_x(n)
  =-B^{(-)}(x)\phi^{\text{d}\,(-)}_{x+1}(n).
  \label{dual:apmphin2}
\end{align}

The data of $R_i(x)$ are\\
\underline{dual big $q$-Jacobi} :
\begin{align}
  R_1(z)&=(q^{\frac12}-q^{-\frac12})^2(z-aq),\n
  R_0(z)&=(q^{\frac12}-q^{-\frac12})^2(z-aq)^2
  \ \ \Bigl(\Rightarrow\sqrt{R_1(z)^2+4R_0(z)}=(q^{-1}-q)|aq-z|\Bigr),\\
  R_{-1}(z)&=(q^{\frac12}-q^{-\frac12})^2
  \bigl((1+abq)z^2+aq(b+c+a^{-1}c-1-2abq)z\n
  &\phantom{=(q^{\frac12}-q^{-\frac12})^2\bigl(}
  +acq(1-aq)(1-abc^{-1}q)\bigr),\nonumber
\end{align}
\underline{dual big $q$-Laguerre} :
\begin{align}
  R_1(z)&=(q^{\frac12}-q^{-\frac12})^2(z-aq),\n
  R_0(z)&=(q^{\frac12}-q^{-\frac12})^2(z-aq)^2
  \ \ \Bigl(\Rightarrow\sqrt{R_1(z)^2+4R_0(z)}=(q^{-1}-q)|aq-z|\Bigr),\\
  R_{-1}(z)&=(q^{\frac12}-q^{-\frac12})^2
  \bigl(z^2+aq(b+a^{-1}b-1)z+abq(1-aq)\bigr),\nonumber
\end{align}
\underline{dual Al-Salam-Carlitz I} :
\begin{align}
  R_1(z)&=(q^{\frac12}-q^{-\frac12})^2(z-1),\n
  R_0(z)&=(q^{\frac12}-q^{-\frac12})^2(z-1)^2
  \ \ \Bigl(\Rightarrow\sqrt{R_1(z)^2+4R_0(z)}=(q^{-1}-q)|1-z|\Bigr),\\
  R_{-1}(z)&=(q^{\frac12}-q^{-\frac12})^2
  \bigl(z^2+(a-1)z+aq^{-1}\bigr).\nonumber
\end{align}

\subsection{$q$-Meixner and $q$-Charlier}

The polynomials $P_n(\eta)$ \eqref{qMpol} for $q$M and \eqref{qCp} for $q$C,
satisfy $P_n(0)=1$ and the three term recurrence relation has the
following form
\begin{equation}
  \eta P_n(\eta)=A_nP_{n+1}(\eta)-(A_n+C_n)P_n(\eta)+C_nP_{n-1}(\eta),
  \label{3term_P_ii')}
\end{equation}
where $A_n$ and $C_n$ are given in \eqref{3coef} and \eqref{qC3coef}.
{}From this form, we adopt the sinusoidal coordinate for the closure relation
\eqref{closurerel1} as $\eta(x)=q^{-x}-1$, which satisfies $\eta(0)=0$.
This situation is the same as I.

The system has two types of eigenvectors, the ordinary
$\{\phi_n(x)\}$ and the supplementary $\{\phi^{(-)}_n(x)\}$
\eqref{qMori}--\eqref{qMsup} and \eqref{qCsup0}--\eqref{qCsup},
with the corresponding eigenvalues
$\{\mathcal{E}(n)=1-q^n\}$ and $\{\mathcal{E}^{\prime}(n)\}$ 
\eqref{qMsupeig} and \eqref{qCsupeig}, respectively.
The latter is a monotonously decreasing function of $n$ and it lies above the
former $\mathcal{E}^{\prime}(n')>\mathcal{E}(n)$.
We can confirm that \eqref{alphapmE}--\eqref{apmphin} hold for the ordinary
sector, whereas for the supplementary sector the order of $n$ is reversed:
\begin{align}
  &\alpha_{\pm}\bigl(\mathcal{E}^{\prime}(n)\bigr)
  =\mathcal{E}^{\prime}(n\mp 1)-\mathcal{E}^{\prime}(n),
  \label{dual:alphapmE3}\\
  &R_{-1}\bigl(\mathcal{E}^{\prime}(n)\bigr)
  R_0\bigl(\mathcal{E}^{\prime}(n)\bigr)^{-1}
  =A^{(-)}_n+C^{(-)}_n,
  \label{dual:Rm1E3}\\
  &a^{(+)}\phi^{(-)}_n(x)=C^{(-)}_n\phi^{(-)}_{n-1}(x),
  \ \ a^{(-)}\phi^{(-)}_n(x)=A^{(-)}_n\phi^{(-)}_{n+1}(x),
  \label{dual:apmphin3}
\end{align}
where $A^{(-)}_n$ and $C^{(-)}_n$ are given in \eqref{qM:Amn,Cmn} for $q$M
and \eqref{qC3coef} with the replacement $a\to a^{-1}$ for $q$C.

The data of $R_i(x)$ are\\
\underline{$q$-Meixner} :
\begin{align}
  R_1(z)&=(q^{\frac12}-q^{-\frac12})^2(z-1),\n
  R_0(z)&=(q^{\frac12}-q^{-\frac12})^2(z-1)^2
  \ \ \Bigl(\Rightarrow\sqrt{R_1(z)^2+4R_0(z)}=(q^{-1}-q)|1-z|\Bigr),\\
  R_{-1}(z)&=(q^{\frac12}-q^{-\frac12})^2
  \bigl(z^2-(1+bc+c)z+bc-q^{-1}c\bigr),\nonumber
\end{align}
\underline{$q$-Charlier} :
\begin{align}
  R_1(z)&=(q^{\frac12}-q^{-\frac12})^2(z-1),\n
  R_0(z)&=(q^{\frac12}-q^{-\frac12})^2(z-1)^2
  \ \ \Bigl(\Rightarrow\sqrt{R_1(z)^2+4R_0(z)}=(q^{-1}-q)|1-z|\Bigr),\\
  R_{-1}(z)&=(q^{\frac12}-q^{-\frac12})^2
  \bigl(z^2-(1+a)z-q^{-1}a\bigr).\nonumber
\end{align}

\subsection{Discrete $q$-Hermite $\II$}

The polynomials $P_n(\eta)$ \eqref{dqHII:Pn} satisfy
the three term recurrence relation of the following form
\begin{equation}
  \eta P_n(\eta)=A_nP_{n+1}(\eta)+B_nP_n(\eta)+C_nP_{n-1}(\eta),
  \ \ A_n=q^{-n},\ \ B_n=0,\ \ C_n=q^{-n}-1.
  \label{3term_P_iii)}
\end{equation}
{}From this form, we adopt the sinusoidal coordinate for the closure relation
\eqref{closurerel1} as $\eta(x)=\eta^{(\pm)}(x)=\pm cq^x$.
The data of $R_i(z)$ do not depend on $(\pm)$:
\begin{align}
  R_1(z)&=(q^{\frac12}-q^{-\frac12})^2(z-1),\n
  R_0(z)&=(q^{\frac12}-q^{-\frac12})^2(z-1)^2
  \ \ \Bigl(\Rightarrow\sqrt{R_1(z)^2+4R_0(z)}=(q^{-1}-q)|1-z|\Bigr),\\
  R_{-1}(z)&=0.\nonumber
\end{align}
The following relations can be easily verified:
\begin{align}
  &\alpha_{\pm}\bigl(\mathcal{E}(n)\bigr)
  =\mathcal{E}(n\pm 1)-\mathcal{E}(n),
  \label{dqHII:alphapmE}\\
  &R_{-1}\bigl(\mathcal{E}(n)\bigr)R_0\bigl(\mathcal{E}(n)\bigr)^{-1}
  =B_n=0,
  \label{dqHII:Rm1E}\\
  &a^{(+)}\phi^{(\pm)}_n(x)=A_n\phi^{(\pm)}_{n+1}(x),\quad
  a^{(-)}\phi^{(\pm)}_n(x)=C_n\phi^{(\pm)}_{n-1}(x).
  \label{dqHII:apmphin}
\end{align}

\subsection{$q$-Laguerre}

{}From the recurrence relation \eqref{qL3term}, we adopt as the sinusoidal
coordinate for the closure relation $1+\eta(x;\bm{\lambda})$.
The data $R_i(z)$ are
\begin{align}
  R_1(z)&=(q^{\frac12}-q^{-\frac12})^2(z-1),\n
  R_0(z)&=(q^{\frac12}-q^{-\frac12})^2(z-1)^2
  \ \ \Bigl(\Rightarrow\sqrt{R_1(z)^2+4R_0(z)}=(q^{-1}-q)|1-z|\Bigr),\\
  R_{-1}(z)&=-(q^{\frac12}-q^{-\frac12})^2
  \bigl(z^2+(a^{-1}-1)z+a^{-1}q^{-1}\bigr).\nonumber
\end{align}
It is straightforward to verify the fundamental relations
\eqref{alphapmE}--\eqref{apmphin}.

\goodbreak

\end{document}